\documentclass[preprint,12pt]{elsarticle}

\usepackage{amssymb}
\usepackage{amsthm}
\usepackage{mathrsfs}
\usepackage{enumerate}
\usepackage{amsmath}
\usepackage{amsfonts}
\usepackage{amssymb}
\usepackage{multicol}
\usepackage{bm}
\usepackage[table]{xcolor}
\usepackage{caption}
\usepackage{tikz}
\usepackage{graphicx}
\usepackage{mathtools} 
\usepackage{titlesec}    
\usepackage{etoolbox}
\usepackage[colorlinks]{hyperref}
\usepackage[utf8]{inputenc}
\usepackage[english]{babel}
\usepackage{xcolor,soul}
\sethlcolor{yellow}
\usepackage[framemethod=tikz]{mdframed}
\usepackage{floatrow}
\usepackage{multirow}
\makeatletter
\patchcmd{\ttlh@hang}{\parindent\z@}{\parindent\z@\leavevmode}{}{}
\patchcmd{\ttlh@hang}{\noindent}{}{}{}
\makeatother
\usepackage{titletoc}
\usepackage{fancyhdr}
\usepackage{indentfirst}
\usepackage[margin=0.8in,vmargin=2.7cm]{geometry}

\usepackage{mdframed}
\newmdenv[skipabove=4mm,linewidth=1]{MyFrame}   

\numberwithin{equation}{section}

\journal { }

\newtheorem{theorem}{Theorem}[section]

\newtheorem{remark}[theorem]{Remark}

\begin{document}
\begin{frontmatter}

\title{Implicit and semi-implicit second-order time stepping methods for the Richards equation}

\author[addr1,addr2]{Sana Keita\corref{corres}}\ead{skeita@uottawa.ca}  
\cortext[corres]{Corresponding author.}
\author[addr1,addr2]{Abdelaziz Beljadid}\ead{abdelaziz.beljadid@um6p.ma}
\author[addr2]{Yves Bourgault}\ead{ybourg@uottawa.ca}
\address[addr1]{International Water Research Institute, Mohammed VI Polytechnic University,  Morocco} 
\address[addr2]{Department of Mathematics and Statistics, University of Ottawa, Canada \qquad\qquad\qquad\qquad}

\begin{abstract}
This study concerns numerical methods for efficiently solving the Richards equation where different weak formulations and computational techniques are analyzed. The spatial discretizations are based on standard or mixed finite element methods. Different implicit and semi-implicit temporal discretization techniques of second-order accuracy are studied. To obtain a linear system for the semi-implicit schemes, we propose second-order techniques using extrapolation formulas and/or semi-implicit Taylor approximations for the temporal discretization of nonlinear terms. A numerical convergence study and a series of numerical tests are performed to analyze efficiency and robustness of the different schemes. The developed scheme, based on the proposed temporal extrapolation techniques and the mixed formulation involving the saturation and pressure head and using the standard linear Lagrange element, performs better than other schemes based on the saturation and the flux and using the Raviart-Thomas elements. The proposed semi-implicit scheme is a good alternative when implicit schemes meet convergence issues.
\end{abstract}

\end{frontmatter}

\section{Introduction}
The prediction of flows through unsaturated porous media is of great importance in a wide variety of applications in science and engineering. There are plenty of relevant applications in many fields such as hydrology, agriculture, water resources 
management, risk assessment of environmental contaminants, enhanced oil recovery, etc. The Richards equation \cite{Richards1931} describes fluid flows in unsaturated porous medium due to the actions of gravity and capillarity. This equation was derived based on a multiphase flow extension of Darcy’s law and the principle of mass conservation \cite{Bear1972,Richards1931}. Some analytical solutions of the Richards equation are available, but they are limited to very simple geometries and specific initial and boundary conditions. Most practical situations require the use of numerical techniques to produce approximate solutions. The design and analysis of appropriate numerical schemes for the Richards equation is a very challenging task due to the highly nonlinear nature of the equation. Despite intensive research, there is still a strong need for more robust numerical schemes for computing multiphase flows in unsaturated porous media.

Much research has been done on the numerical treatment of the Richards equation. Various classes of approaches are available for the spatial discretization. Standard  finite element methods \cite{Arbogast1993,Forsyth1997,Slodicka2002} and mixed finite element methods \cite{Arbogast1996,Bause2008,Bause2004,Belfort2009,Bergamaschi1999,Farthing2003,ISPop2004,Radu2004,Radu2014,Schneid2004} are commonly used, often combined with mass lumping to ensure nonoscillatory solutions \cite{Belfort2009,Celia1990}. Locally mass-conservative finite volume methods for the Richards equation were proposed in \cite{Eymard1999,Farthing2000,Manzini2004,McBride2006,Forsyth1995}, while others studies have used finite difference methods \cite{Arbogast1993,Celia1990,CLement1994,Haverkamp1977,Miller1998} and a control volume finite element method \cite{Oulhaj2018}.   Another well-known numerical approach used to solve the Richards equation is the Kirchhoff transformation technique \cite{Arbogast1996,Berninger2011,Berninger2015,ISPop2004,Radu2004,Suk2019}. While this technique is very promising, it is limited to specific functions of the capillary pressure \cite{Suk2019}. In terms of temporal discretization techniques,  a proper treatment of the temporal derivative is required for reliable numerical solutions due to the highly nonlinear and stiff nature of the equation. Most of the numerical models employ an implicit or semi-implicit scheme for the temporal discretization.   The main reason is the need of a stable discretization allowing reasonable time steps since, due to the presence of a stiff nonlinear diffusion term, an explicit scheme requires very small time step which makes this scheme impractical to use in terms of efficiency. The backward Euler method is the most commonly used \cite{Arbogast1993,Arbogast1996,Bause2008,Bause2004,Belfort2009,Eymard1999,Forsyth1997,Oulhaj2018,Radu2004,Radu2014,Schneid2004} and it results in a nonlinear system. The conventional approach consists in solving this nonlinear system using an iterative procedure. There are several iterative methods available for the linearization of the resulting nonlinear system. The Newton method, which is quadratically convergent, is often used. However, its convergence occurs only when one starts with an initial guess which is very close to the solution. A method for Richards' equation has been proposed in \cite{Casulli2010}  that exploits the typical behavior of the soil hydraulic functions to define nested iterations that guarantee mass conservation and quadratic convergence of a Newton-type method. The Picard method, which is easy to implement, is also widely used, although it is only linearly convergent. A modified Picard method is used to develop a mass-conservative numerical scheme in \cite{Celia1990}. In \cite{Lott2012}, an acceleration technique is applied to the modified Picard method to provide faster convergence while maintaining the advantages, i.e., the mass conservation and lower memory requirements with respect to Newton's method. A comparison between Picard and Newton iteration approaches and a combination of the modified Picard and Newton methods, as well as $L$-scheme and its combination with  Newton method, can be found in \cite{Lehmann1998,List2016,Paniconi1994}.

Newton and Picard methods can produce good results but these methods entail computational costs associated with having to solve a nonlinear system at each time step. Noniterative methods require a single matrix assembly and factorization per time step and thus are an attractive alternative to traditional iterative methods for the Richards equation \cite{Kavetski2002,Paniconi1991,Ross2003,Zha2013,Zha2016}. Second order time discretization has shown benefits over traditional first order temporal approximation for modeling unsaturated flows through porous media \cite{Farthing2003,Kavetski2002,Paniconi1991,Zha2016}. Furthermore, in our last studies \cite{Keita2021,Keita2020} we showed the advantages of semi-implicit techniques in solving highly nonlinear parabolic equations. This motivates an investigation of second order noniterative approaches to solving the nonlinear Richards equation. This study focuses on finite element methods for the Richards equation where three different weak formulations of the equation are considered. The first one is expressed in terms of the saturation as the only unknown. The second is a mixed formulation based on the saturation and the flux. This mixed formulation is mostly studied and the flux is approximated using the Raviart-Thomas elements \cite{Bause2004,Belfort2009,Bergamaschi1999,Farthing2003,ISPop2004,Radu2014} or the Brezzi-Douglas-Marini element \cite{Bause2008}. The third one is also a mixed formulation, but it consists in approximating simultaneously the saturation and the pressure head by using the standard linear Lagrange element.  All the formulations include a regularization technique to handle degeneracies that could occur in the gradient of the pressure head for some constitutive relationships of the capillary pressure when the saturation reaches the value of one. For each weak formulation, iterative and noniterative second-order methods are studied. We use the second-order backward differentiation formula (BDF2) for the time integration with fixed time step. We note that there exists adaptive BDF2 methods which allow the use of varying time steps \cite{Becker1998,Crouzeix1984,Emmrich2009}. In \cite{Baron2017}, the fully implicit adaptive BDF2 method is applied to Richards' equation.  Our schemes can be modified to allow varying time steps by using the techniques developed in \cite{Baron2017,Becker1998,Crouzeix1984,Emmrich2009}, but this would require modifications in the treatment of the nonlinear terms since some of our schemes are based on semi-implicit BDF2 methods. To obtain a linear system for the noniterative schemes, a second-order extrapolation formula and/or a semi-implicit Taylor expansion of second-order truncation error are used for the nonlinear terms. To the best of our knowledge, the linear mixed scheme based on the saturation and the pressure head, using the standard Lagrange linear element for the spatial discretization and the linearization techniques based on Taylor expansion, has not been explored in previous studies. A main goal is to determine which weak formulation using such temporal discretization is the most competitive among a relatively large set of methods. Here, competitive is to be understood in terms of accuracy versus computational requirements, but also ease of implementation.  Numerical convergence analysis is performed using manufactured and reference solutions to investigate the spatial and temporal order of convergence of each scheme. 

The outline of the paper is as follows. In section~\ref{SectModelEquations}, the Richards equation is introduced. In section~\ref{SectWeakFormulations}, three different weak formulations of the Richards equation including a regularization technique are presented. In section~\ref{SectNumMethods}, several numerical schemes are proposed. In section~\ref{SectNumTests}, convergence rates, accuracy and robustness of the numerical schemes are investigated and compared by means of manufactured and exact solutions and selected test cases. Section~\ref{SectConclusion} provides concluding remarks.

\section{Model equations}
\label{SectModelEquations}
Let $\Omega\subset\mathbb{R}^2$ be an open bounded set with sufficiently smooth boundary $\Gamma$ and $T>0$ a fixed time. The Richards equation is a standard, commonly-used continuum model for describing flow in unsaturated porous media. It can be written in three different forms, namely the ``pressure head form", ``water content form" and ``mixed form" \cite{Celia1990}. The mixed form is given by
\begin{equation}
\dfrac{\partial \theta}{\partial t} - \nabla\cdot[K_s(\bm{x})K_r(\psi)\nabla(\psi+z)]=0\quad\text{in}\quad(0,T)\times\Omega,
\label{RichardsEqMixedForm}
\end{equation}
where $\theta$ is the water content, $\psi$ is the pressure head,  $K_s$ is the saturated hydraulic conductivity, $K_r$ is the water relative permeability, and $z$ is the vertical coordinate and is assumed to be positive upward. The water relative permeability $K_r$ accounts for the effect of partial saturation. While $K_r$ and $\psi$ include physical parameters which depend on the spatial heterogeneity of the medium, for simplicity we write $K_r(\psi)$.  To complete the model formulation of the Richards equation \eqref{RichardsEqMixedForm}, one must specify the constitutive relationship to describe the interdependence between the relative permeability, pressure head and water content. There are several  empirical relationships that are used in modeling, for instance, van Genuchten \cite{vanGenuchten1980}, Brooks-Corey \cite{Brooks1966} and Gardner \cite{Gardner1958} models.

By introducing the effective saturation
\begin{equation}
S=\dfrac{\theta-\theta_r}{\theta_s-\theta_r},
\end{equation}
where $\theta_s$ and $\theta_r$ are the saturated and residual volumetric water contents, respectively, the Richards equation \eqref{RichardsEqMixedForm} can be written as 
\begin{equation}
\phi(\bm{x})\dfrac{\partial S}{\partial t} - \nabla\cdot[K_s(\bm{x})K_r(\psi)\nabla(\psi+z)]=0\quad\text{ in }(0,T)\times\Omega,
\label{RichardsEq}
\end{equation}
where $\phi$ is given by
\begin{equation}
\phi=\theta_s-\theta_r.
\end{equation} 
The water flux $\bm{q}$ is expressed by the extended form of Darcy's law \cite{Philip1969}:
\begin{equation}
\bm{q}=-K_sK_r\nabla(\psi+z).
\label{RichardsEqFlux}
\end{equation}
Equation \eqref{RichardsEq} is to be supplemented by an initial condition
\begin{equation}
S(\bm{x},0)=S_0(\bm{x})
\label{RichardsEqIC}
\end{equation}
and  boundary conditions
\begin{equation}
\bm{q}\cdot\bm{n}=q_{\Gamma}\quad \text{on}\quad \Gamma,
\label{RichardsEqBC}
\end{equation}
where $\bm{n}$ is the outward unit normal vector to the boundary $\Gamma$ of the domain $\Omega$ and $q_{\Gamma}\in\mathbb{R}$.

\section{Weak formulations}
\label{SectWeakFormulations}
For most constitutive relationships used in practice, the pressure head can be expressed as  
\begin{equation}
\psi(S)=h_{cap}(\bm{x}){J}(S),
\label{WaterSectionHead}
\end{equation}
where $h_{cap}:\Omega\rightarrow \mathbb{R}^+$ is the capillary rise function which depends on space for heterogeneous medium and $J:(0,1]\rightarrow (-\infty,0]$ is the Leverett $J$-function \cite{Leverett1941} which is of class $\mathcal{C}^1$ on $(0, 1)$.  The Leverett $J$-function tends to infinity as the effective saturation $S$ goes to zero, and for some constitutive relationships, $J$ is not differentiable at $S=1$, for instance, the constitutive relationship used by Haverkamp et al.\ \cite{Haverkamp1977} (see \eqref{HaverkampJfunction} below for the explicit form of the $J$-function) and the van Genuchten model \cite{vanGenuchten1980} (see \eqref{vanGenuchtenJfunction} below). The computations of solutions may encounter numerical issues when the saturation $S$ approaches the value of zero or one due to the lack of regularity of the constitutive functions. Some $J$-functions are differentiable at $S=1$ such as for the Brooks-Corey model \cite{Brooks1966}. In case of lack of differentiability, mathematical regularization techniques have been proposed in \cite{Schweizer2011,Schweizer2007} to overcome the degeneracy in the pressure head. In this work, we consider the regularization consisting in replacing the derivative $J^{\prime}$ by $J_{\delta}^{\prime}$ which satisfies
\begin{equation}
J_{\delta}^{\prime}(S)=
\left\lbrace
\begin{aligned}
&J^{\prime}(S),\quad \forall S\in(0,1-\delta),\\
&J^{\prime}(1-\delta),\quad \forall S\in[1-\delta, 1],
\end{aligned}
\right.
\label{RegularizedJfunction}
\end{equation}
where $0<\delta\ll1$. The regularized $J_{\delta}^{\prime}$ is a continuous function on the interval $(0,1]$. This is very useful for the convergence of numerical methods involving the computation of the derivative of $J$. The $J$-function in \eqref{WaterSectionHead} may include parameters which depend on the heterogeneity of the medium. Here, the impact of the variability due to heterogeneity is considered dominant by the variation of the intrinsic permeability included in the capillary rise $h_{cap}$ \cite{Leverett1941}. By expanding \eqref{WaterSectionHead}, we get
\begin{equation}
\nabla\psi=J(S)\nabla h_{cap}+h_{cap}J_{\delta}^{\prime}(S)\nabla S,\quad \forall S\in (0,1].
\label{WaterSectionHeadExpanded}
\end{equation}

Substituting \eqref{WaterSectionHeadExpanded} in the Richards equation \eqref{RichardsEq}  leads to
\begin{equation}
\begin{aligned}
&\phi\dfrac{\partial S}{\partial t} - \nabla\cdot[K_sK_r(\psi(S))\big(J(S)\nabla h_{cap}+h_{cap}J_{\delta}^{\prime}(S)\nabla S + \nabla z\big)]=0.
\label{RichardsEqMod}
\end{aligned}
\end{equation}
Multiplying \eqref{RichardsEqMod} by a test function $v$, integrating over the domain $\Omega$ and using  \eqref{RichardsEqBC} give rise to 
\begin{equation}
\begin{aligned}
\int_{\Omega}\phi\dfrac{\partial S}{\partial t}v\,d\bm{x} &+ \int_{\Omega}\big[K_sK_r(\psi(S))\big(J(S)\nabla h_{cap}+h_{cap}J_{\delta}^{\prime}(S)\nabla S+\nabla z\big)\big]\cdot\nabla v\,d\bm{x} \\
&+\int_{\Gamma}q_{\Gamma}v\,d\bm{x}=0,\quad\forall v\in H^1(\Omega).\\
\end{aligned}
\label{RichardsEqWeakForm}
\end{equation}
A weak formulation for \eqref{RichardsEqMod} is: Find $S \in H^1(\Omega)$ such that \eqref{RichardsEqWeakForm} is satisfied.

\begin{remark}
Dirichlet conditions could be imposed on the boundary $\Gamma$. In this case, the solution $S$ is imposed on the boundary $\Gamma$ and the test function $v$ is set to zero  on $\Gamma$. The same idea applies to all the variational formulations below.
\end{remark}

By substituting \eqref{RichardsEqFlux} and  \eqref{WaterSectionHeadExpanded}, the Richards  equation \eqref{RichardsEq}  can be written in the mixed form 
\begin{equation}
\begin{aligned}
&\phi\dfrac{\partial S}{\partial t} + \nabla\cdot\bm{q}=0,\\
&\bm{q}+K_sK_r(\psi)\big(J(S)\nabla h_{cap}+h_{cap}J_{\delta}^{\prime}(S)\nabla S+\nabla z\big)=0.
\end{aligned}
\label{RichardsEqMixed2}
\end{equation}
Multiplying each equation in \eqref{RichardsEqMixed2} by a test function and integrating over the domain lead to 
\begin{equation}
\begin{aligned}
&\int_{\Omega}\phi\dfrac{\partial S}{\partial t}v\,d\bm{x}-\int_{\Omega}\bm{q}\cdot\nabla v\,d\bm{x}+\int_{\Gamma}q_{\Gamma}v\,d\bm{x}=0,\quad\forall v\in H^1(\Omega)\\
&\int_{\Omega}\bm{q}\cdot\bm{w}\,d\bm{x}+\int_{\Omega}\big[K_sK_r(\psi(S))\big(J(S)\nabla h_{cap}+h_{cap}J_{\delta}^{\prime}(S)\nabla S+\nabla z\big)\big]\cdot\bm{w}\,d\bm{x}=0,\quad\forall\bm{w}\in (L^2(\Omega))^2.
\end{aligned}
\label{RichardsEqMixed2WeakForm}
\end{equation} 
A weak formulation for \eqref{RichardsEqMixed2} is: Find $(S,\bm{q}) \in H^1(\Omega)\times H(div;\Omega)$ such that \eqref{RichardsEqMixed2WeakForm} is satisfied, where
\begin{equation}
H(div;\Omega)=\Big\lbrace\bm{w}\in L^2(\Omega)\times L^2(\Omega): \nabla\cdot\bm{w}\in L^2(\Omega)\Big\rbrace.
\label{Hdivspace}
\end{equation}

By using the pressure head relation \eqref{WaterSectionHead}, the Richards equation \eqref{RichardsEq} can also be written in the alternative mixed formulation 
\begin{equation}
\begin{aligned}
&\phi\dfrac{\partial S}{\partial t} - \nabla\cdot[K_sK_r(\psi)\nabla(\psi+z)]=0,\\
&\psi-h_{cap}J(S)=0.
\end{aligned}
\label{RichardsEqMixed1}
\end{equation}
Multiplying each equation in \eqref{RichardsEqMixed1} by a test function and integrating over the domain give
\begin{equation}
\begin{aligned}
&\int_{\Omega}\phi\dfrac{\partial S}{\partial t}v\,d\bm{x} + \int_{\Omega}[K_sK_r(\psi)\big(\nabla\psi+\nabla z\big)]\cdot\nabla v\,d\bm{x}+\int_{\Gamma}q_{\Gamma}v\,d\bm{x}=0,\quad\forall v\in H^1(\Omega),\\
&\int_{\Omega}\psi w\,d\bm{x}-\int_{\Omega}h_{cap}J(S)w\,d\bm{x}=0,\quad\forall w\in H^1(\Omega).\\
\end{aligned}
\label{RichardsEqMixed1WeakForm}
\end{equation}
A weak formulation for \eqref{RichardsEqMixed1} is: Find $(S,\psi) \in H^1(\Omega)\times H^1(\Omega)$ such that \eqref{RichardsEqMixed1WeakForm} is satisfied. The mixed formulation {\eqref{RichardsEqMixed1WeakForm}} is quite simple and leads to a smaller linear system to solve than the mixed formulation \eqref{RichardsEqMixed2WeakForm} as we will see in the next section. It is also advantageous when it comes to the choice of the finite element spaces for the spatial discretization since the standard linear Lagrange element is used instead of the Brezzi-Douglas-Marini element or the Raviart-Thomas elements, which are often used for the mixed formulation \eqref{RichardsEqMixed2WeakForm} and are more difficult to implement.

\section{Numerical methods}
\label{SectNumMethods}
\subsection{Temporal discretization}
\label{SectTemDis}
The time interval $[0,T]$ is discretized as
\begin{equation}
\Delta t=\dfrac{T}{N},\quad t_n=n \Delta t,\quad  n=0,1,2,...,N,
\end{equation} 
where $\Delta t$ is the time step used. We set $S^n\simeq S(t_n)$, $\psi^n\simeq \psi(t_n)$ and $\bm{q}^n\simeq\bm{q}(t_n)$ for $n=0,...,N$.\\
We propose implicit and semi-implicit time stepping discretizations for each of the formulations \eqref{RichardsEqMod}, \eqref{RichardsEqMixed2} and \eqref{RichardsEqMixed1}. Few second order accurate time-stepping schemes were proposed for Richards equation based on the Crank-Nicolson method \cite{Cai2020,Zha2013,Zha2016} and BDF methods \cite{Baron2017,Cumming2011} in their fully implicit form requering Newton or Picard iterations. The Crank-Nicolson method is A-stable but lacks L-stability, and  may lead to non-monotone solutions for larger time steps. Among the simplest methods for solving stiff ODEs, the second-order backward differentiation formula (BDF2) is L-stable. A concise discussion on these concepts and two time-stepping methods are presented in \cite{LeVeque2007}. Consequently, all the methods proposed here approximate the time derivative using the BDF2 formula. We first consider the fully implicit BDF2 method, which is the most stable and accurate, but requires the solution of stiff nonlinear systems at each time step. To avoid the solution of nonlinear systems, semi-implicit (also called IMEX) variants of the BDF2 method are available. These variants are usually based on the extrapolation of nonlinear terms \cite{Akrivis2015,Ascher1995,Ethier2008,Keita2021}, at the expense of a slightly more restrictive time step condition for stability, although with time steps that are still independent from the space step. The extension of these semi-implicit methods to Richards' equation requires some care, since those were so far applied to PDEs with other forms of nonlinear terms. The semi-implicit BDF2 method developed in \cite{Keita2021} for the Cahn-Hilliard equation can be extended for solving Richards' equation. To obtain a linear scheme of second-order time-accuracy for the semi-implicit temporal discretizations, a second-order extrapolation formula is employed for the relative permeability function $K_r(\psi)$ and the regularized derivative $J_{\delta}^{\prime}(S)$, while the Leverett $J$-function $J(S)$ is evaluated using a Taylor approximation with second-order truncation error.

\subsubsection{Formulation \eqref{RichardsEqMod}}
A second-order implicit temporal discretization for \eqref{RichardsEqMod} is given by
\begin{equation}
\begin{aligned}
&\phi\dfrac{3S^{n+1}-4S^{n}+S^{n-1}}{2\Delta t} - \nabla\cdot[K_sK_r(\psi(S^{n+1}))\big(J(S^{n+1})\nabla h_{cap}\\
&\qquad\qquad\qquad\qquad\qquad\qquad\qquad\qquad +h_{cap}J_{\delta}^{\prime}(S^{n+1})\nabla S^{n+1}+\nabla z\big)]=0.
\end{aligned}
\label{RichardsEqImpl}
\tag{$D_{1i}$}
\end{equation}
For a second order semi-implicit temporal discretization of \eqref{RichardsEqMod}, we propose
\begin{equation}
\begin{aligned}
&\phi\dfrac{3S^{n+1}-4S^{n}+S^{n-1}}{2\Delta t} - \nabla\cdot\Big(\big[K_s\big(2K_r(\psi(S^n))-K_r(\psi(S^{n-1}))\big)\big]\\
&\quad\big[\big(J(S^n)+J_{\delta}^{\prime}(S^n)(S^{n+1}-S^n)\big)\nabla h_{cap}+\big(2J_{\delta}^{\prime}(S^n)-J_{\delta}^{\prime}(S^{n-1})\big)h_{cap}\nabla S^{n+1}+\nabla z\big]\Big)=0.
\end{aligned}
\label{RichardsEqSemImp}
\tag{$D_{1s}$}
\end{equation} 

\subsubsection{Formulation \eqref{RichardsEqMixed2}}
An implicit discretization of second-order accuracy for \eqref{RichardsEqMixed2} reads 
\begin{equation}
\begin{aligned}
&\phi\dfrac{3S^{n+1}-4S^{n}+S^{n-1}}{2\Delta t} - \nabla\cdot\bm{q}^{n+1}=0,\\
&\bm{q}^{n+1} + K_sK_r(\psi(S^{n+1}))\big[J(S^{n+1})\nabla h_{cap}+h_{cap}J_{\delta}^{\prime}(S^{n+1})\nabla S^{n+1}+\nabla z\big]=0
\end{aligned}
\label{RichardsEqMixed2Impl}
\tag{$D_{2i}$}
\end{equation}
and a second-order semi-implicit discretization is given by
\begin{equation}
\begin{aligned}
&\phi\dfrac{3S^{n+1}-4S^{n}+S^{n-1}}{2\Delta t} - \nabla\cdot\bm{q}^{n+1}=0,\\
&\bm{q}^{n+1}+\big[K_s\big(2K_r(\psi(S^n))-K_r(\psi(S^{n-1}))\big)\big]\big[\big(J(S^n)+J_{\delta}^{\prime}(S^n)(S^{n+1}-S^n)\big)\nabla h_{cap}\\
&\qquad\qquad\qquad\qquad+\big(2J_{\delta}^{\prime}(S^n)-J_{\delta}^{\prime}(S^{n-1})\big)h_{cap}\nabla S^{n+1}+\nabla z\big]=0.
\end{aligned}
\label{RichardsEqMixed2SemImp}
\tag{$D_{2s}$}
\end{equation}

\subsubsection{Formulation \eqref{RichardsEqMixed1}}
For the formulation \eqref{RichardsEqMixed1}, we consider the second order implicit time integration
\begin{equation}
\begin{aligned}
&\phi\dfrac{3S^{n+1}-4S^{n}+S^{n-1}}{2\Delta t} - \nabla\cdot[K_sK_r(\psi^{n+1})\nabla(\psi^{n+1}+z)]=0,\\
&\psi^{n+1}-h_{cap}J(S^{n+1})=0,
\end{aligned}
\label{RichardsEqMixed1Impl}
\tag{$D_{3i}$}
\end{equation}
and the semi-implicit discretization of second-order accuracy 
\begin{equation}
\begin{aligned}
&\phi\dfrac{3S^{n+1}-4S^{n}+S^{n-1}}{2\Delta t} - \nabla\cdot[K_s\big(2K_r(\psi^n)-K_r(\psi^{n-1})\big)\nabla(\psi^{n+1}+z)]=0,\\
&\psi^{n+1}-h_{cap}\big(J(S^n)+J_{\delta}^{\prime}(S^n)(S^{n+1}-S^n)\big)=0.
\end{aligned}
\label{RichardsEqMixed1Semi}
\tag{$D_{3s}$}
\end{equation}

\subsection{Spatial discretization}
Let $\mathcal{T}_h$ be a triangulation of the domain $\Omega$ into disjoint triangles $\kappa$, with $h_{\kappa}:=diam(\kappa)$ and $h:=\max_{\kappa\in \mathcal{T}_h}h_{\kappa}$, so that 
\begin{equation}
\Omega=\bigcup_{\kappa\in\mathcal{T}_h}\kappa.
\end{equation}
Let $\mathbf{J}$ be the set of nodes of $\mathcal{T}_h$, with coordinates $\lbrace \bm{x}_j\rbrace_{j\in\mathbf{J}}$. We  consider the following finite element spaces 
\begin{align}
\mathcal{V}_h&=\Big\lbrace v_h\in \mathcal{C}^0(\Omega,\mathbb{R}): {v_h}_{|_{\kappa}}\text{ is linear,}\quad\forall \kappa\in \mathcal{T}_h \Big\rbrace,\label{Vhspace}\\
\mathcal{S}_h&=\Big\lbrace v_h\in {\mathcal{V}_h}: 0<v_h(\bm{x}_j)\leqslant 1,\quad\forall \bm{x}_j\in\mathcal{T}_h \Big\rbrace,\label{Shspace}\\
\mathcal{R}_h&=\Big\lbrace \bm{v}_h\in H(div;\Omega):{\bm{v}_h}_{|_{\kappa}}\in RT_0,\quad\forall \kappa\in \mathcal{T}_h   \Big\rbrace,\label{Rhspace}
\end{align}
where $RT_0$ is the Raviart-Thomas finite elements of degree $0$ \cite{Fortin1991}. In what follows, we give a finite element method for each of the discretizations in Section~\ref{SectTemDis}. The following equations are obtained from the ones in Section~\ref{SectTemDis} by integrating over the domain, using the boundary condition \eqref{RichardsEqBC} and replacing the Sobolev spaces by the corresponding finite element spaces.

A finite element approximation for \eqref{RichardsEqImpl} is: Given a suitable approximation of the initial solution $S_h^{0}=\Pi_hS_0$ and a proper initialization for $S_h^{1}\in\mathcal{S}_h$, find $\widehat{S}_h^{n+1}\in {\mathcal{V}}_h$ such that
\begin{equation}
\begin{aligned}
&\int_{\Omega}\phi_h\dfrac{3\widehat{S}_h^{n+1}-4S_h^{n}+S_h^{n-1}}{2\Delta t} v_h\,d\bm{x} \\
&\qquad+ \int_{\Omega}\big[K_sK_r(\psi(\widehat{S}_h^{n+1}))\big(J(\widehat{S}_h^{n+1})\nabla h_{cap}+h_{cap}J_{\delta}^{\prime}(\widehat{S}_h^{n+1})\nabla\widehat{S}_h^{n+1}+\nabla z\big)\big]\cdot\nabla v_h\,d\bm{x} \\
&\qquad+\int_{\Gamma}q_{\Gamma}v_h\,d\bm{x}=0,\quad\forall v_h\in{\mathcal{V}}_h.
\end{aligned}
\label{RichardsEqImplScheme}
\end{equation}
Then, find $S_h^{n+1}\in\mathcal{S}_h$, solution of
\begin{equation*}
\left\lbrace
\begin{aligned}
&S_h^{n+1}=\text{arg }\min_{S_h\in {\mathcal{V}_h}}G(S_h), \text{ where }G(S_h):=\dfrac{1}{2}\Vert S_h-\widehat{S}_h^{n+1} \Vert^2_{L^2(\Omega)},\\
&\text{ subject to }\quad c(S_h):=S_h-1\leqslant 0.\\
\end{aligned}
\right.
\label{P}
\tag{$\mathbf{P}$}
\end{equation*}
Our approach is based on the techniques we developed for the Cahn-Hilliard equation \cite{Keita2020}, consisting in solving an optimization problem after the variational problem to force the discrete solution $S_h$ to satisfy the desired physical property $S\leqslant 1$. The optimization problem (\textbf{P}) can be solved by using the Uzawa algorithm \cite{Allaire2007,Ciarlet1989,Uzawa1958} or a projection technique \cite{Keita2020}, and it takes about 3-4 iterations for the convergence. We refer to \cite{Keita2020} for more discussions on these techniques, their implementation, and computational cost. A finite element approximation for \eqref{RichardsEqSemImp} is: Given a suitable approximation of the initial solution $S_h^{0}=\Pi_hS_0$ and a proper initialization for $S_h^{1}\in\mathcal{S}_h$, find $\widehat{S}_h^{n+1}\in \mathcal{V}_h$ such that
\begin{equation}
\begin{aligned}
&\int_{\Omega}\phi_h\dfrac{3\widehat{S}_h^{n+1}-4S_h^{n}+S_h^{n-1}}{2\Delta t} v_h\,d\bm{x} + \int_{\Omega}\Big(\big[K_s\big(2K_r(\psi(S_h^n))-K_r(\psi(S_h^{n-1}))\big)\big]\\
&\quad\big[\big(J(S_h^n)+J_{\delta}^{\prime}(S_h^n)(\widehat{S}_h^{n+1}-S_h^n)\big)\nabla h_{cap}+\big(2J_{\delta}^{\prime}(S_h^n)-J_{\delta}^{\prime}(S_h^{n-1})\big)h_{cap}\nabla\widehat{S}_h^{n+1}+\nabla z\big]\Big)\cdot\nabla v_h\,d\bm{x} \\
&\qquad\qquad\qquad\qquad\qquad\qquad\qquad+\int_{\Gamma}q_{\Gamma}v_h\,d\bm{x}=0,\quad\forall v_h\in{\mathcal{V}}_h,
\end{aligned}
\label{RichardsEqSemiScheme}
\end{equation}
Then, find $S_h^{n+1}\in \mathcal{S}_h$,  solution of $(\mathbf{P})$. In the following, the schemes \eqref{RichardsEqImplScheme} and $(\mathbf{P})$, and \eqref{RichardsEqSemiScheme} and $(\mathbf{P})$, derived from the formulation \eqref{RichardsEqWeakForm} that uses $S$ as the sole unknown, will be referred to as the implicit $S$-scheme and the semi-implicit $S$-scheme, respectively.

The weak formulation \eqref{RichardsEqMixed2WeakForm}-\eqref{Hdivspace} naturally leads to mixed finite element methods. These methods became popular for modeling flow in porous media due to their capability of handling general irregular grids and allowing simultaneous approximation of water content (or pressure head or saturation) and flux as part of the formulation \cite{Bause2004,Belfort2009,Bergamaschi1999,ISPop2004,Radu2014}. The lowest order Raviart-Thomas elements $(RT_0)$ \cite{Fortin1991} is used in these studies to approximate the flux. We consider the mixed finite element method for \eqref{RichardsEqMixed2Impl}: Given a suitable approximation of the initial solution $S_h^{0}=\Pi_hS_0$ and a proper initialization for $S_h^{1}\in\mathcal{S}_h$, find $(\widehat{S}_h^{n+1},\bm{q}_h^{n+1})\in {\mathcal{V}}_h\times\mathcal{R}_h$ such that
\begin{equation}
\begin{aligned}
&\int_{\Omega}\phi_h\dfrac{3\widehat{S}_h^{n+1}-4S_h^{n}+S_h^{n-1}}{2\Delta t} v_h\,d\bm{x}-\int_{\Omega}\bm{q}_h^{n+1}\cdot \nabla v_h\,d\bm{x}+\int_{\Gamma}q_{\Gamma}v_h\,d\bm{x}=0,\quad\forall v_h\in\mathcal{V}_h,\\
&\int_{\Omega}\bm{q}_h^{n+1}\bm{w}_h\,d\bm{x} + \int_{\Omega}K_sK_r(\psi(\widehat{S}_h^{n+1}))\big[J(\widehat{S}_h^{n+1})\nabla h_{cap}\\
&\qquad\qquad\qquad\qquad\qquad +h_{cap}J_{\delta}^{\prime}(\widehat{S}_h^{n+1})\nabla\widehat{S}_h^{n+1}+\nabla z\big]\cdot\bm{w}_h\,d\bm{x}=0, \quad\forall \bm{w}_h\in\mathcal{R}_h. \\
\end{aligned}
\label{RichardsEqMixed2ImplScheme}
\end{equation}
Then, find $S_h^{n+1}\in \mathcal{S}_h$  solution of $(\mathbf{P})$.\\ 
A mixed finite element approximation for \eqref{RichardsEqMixed2SemImp} is: Given a suitable approximation of the initial solution $S_h^{0}=\Pi_hS_0$ and a proper initialization for $S_h^{1}\in\mathcal{S}_h$, find $(\widehat{S}_h^{n+1},\bm{q}_h^{n+1})\in\mathcal{V}_h\times\mathcal{R}_h$ such that
\begin{equation}
\begin{aligned}
&\int_{\Omega}\phi_h\dfrac{3\widehat{S}_h^{n+1}-4S_h^{n}+S_h^{n-1}}{2\Delta t} v_h\,d\bm{x} - \int_{\Omega}\bm{q}_h^{n+1}\cdot\nabla v_h\,d\bm{x}+\int_{\Gamma}q_{\Gamma}v_h\,d\bm{x}=0,\quad\forall v_h\in\mathcal{V}_h,\\
&\int_{\Omega}\bm{q}_h^{n+1}\cdot\bm{w}_h\,d\bm{x} + \int_{\Omega}\Big(\big[K_s\big(2K_r(\psi(S_h^n))-K_r(\psi(S_h^{n-1}))\big)\big]\big[\big(J(S_h^n)+J_{\delta}^{\prime}(S_h^n)(\widehat{S}_h^{n+1}-S_h^n)\big)\nabla h_{cap}\\
&\qquad+(2J_{\delta}^{\prime}(S_h^n)-J_{\delta}^{\prime}(S_h^{n-1}))h_{cap}\nabla\widehat{S}_h^{n+1}+\nabla z\big]\Big)\cdot\bm{w}_h\,d\bm{x}=0,\quad\forall \bm{w}_h\in\mathcal{R}_h.
\end{aligned}
\label{RichardsEqMixed2SemiScheme}
\end{equation}
Then, find $S_h^{n+1}\in\mathcal{S}_h$, solution of $(\textbf{P})$. The schemes \eqref{RichardsEqMixed2ImplScheme} and $(\textbf{P})$, and \eqref{RichardsEqMixed2SemiScheme} and $(\textbf{P})$, derived from the mixed formulation \eqref{RichardsEqMixed2WeakForm} dealing with $S$ and $\bm{q}$ as unknowns, will be referred to as the implicit $(S,\bm{q})$-scheme and the semi-implicit $(S,\bm{q})$-scheme, respectively.

The weak formulation \eqref{RichardsEqMixed1WeakForm} also leads to mixed finite element methods. But unlike the previous mixed formulation, this one is based on the saturation and the pressure head. This mixed formulation is simpler and more advantageous for the spatial discretization since standard Lagrange linear elements are used instead of the Raviart-Thomas finite elements or the Brezzi-Douglas-Marini finite elements which lead to larger systems to solve, and are therefore more demanding in terms of computational cost. For instance, even with the lowest order Raviart-Thomas elements, the number of degrees of freedom for $\bm{q}$ is equal to the number of edges in the mesh, while with linear Lagrange element the number of unknowns for $\psi$ is equal to the number of triangle vertices and there are roughly half many vertices as edges in a 2-D mesh (if not less). We propose a mixed finite element method  for \eqref{RichardsEqMixed1Impl}: Given a suitable approximation of the initial solution $S_h^{0}=\Pi_hS_0$ and  a proper initialization for $S_h^{1}\in\mathcal{S}_h$, find $(\widehat{S}_h^{n+1},\psi_h^{n+1})\in\mathcal{V}_h\times\mathcal{V}_h$ such that
\begin{equation}
\begin{aligned}
&\int_{\Omega}\phi_h\dfrac{3\widehat{S}_h^{n+1}-4S_h^{n}+S_h^{n-1}}{2\Delta t} v_h\,d\bm{x} + \int_{\Omega}\big[K_sK_r(\psi_h^{n+1})(\nabla\psi_h^{n+1}+\nabla z)\big]\cdot\nabla v_h\,d\bm{x}\\ 
&\qquad\qquad\qquad\qquad\qquad\qquad\qquad+\int_{\Gamma}q_{\Gamma}v_h\,d\bm{x}=0,\quad\forall v_h\in{\mathcal{V}}_h,\\
&\int_{\Omega}\psi_h^{n+1}w_h\,d\bm{x}-\int_{\Omega}h_{cap}J(\widehat{S}^{n+1})w_h\,d\bm{x}=0,\quad\forall w_h\in\mathcal{V}_h.\\
\end{aligned}
\label{RichardsEqMixed1ImplScheme}
\end{equation}
Then, find $S_h^{n+1}\in\mathcal{S}_h$, solution of $(\textbf{P})$.\\
A mixed finite element approximation for \eqref{RichardsEqMixed1Semi} is: Given a suitable approximation of the initial solution $S_h^{0}=\Pi_hS_0$ and a proper initialization for $S_h^{1}\in\mathcal{S}_h$, find $(\widehat{S}_h^{n+1},\psi_h^{n+1})\in\mathcal{V}_h\times\mathcal{V}_h$ such that
\begin{equation}
\begin{aligned}
&\int_{\Omega}\phi_h\dfrac{3\widehat{S}_h^{n+1}-4S_h^{n}+S_h^{n-1}}{2\Delta t} v_h\,d\bm{x} + \int_{\Omega}\big[K_s\big(2K_r(\psi_h^n)-K_r(\psi_h^{n-1})\big)\big(\nabla\psi_h^{n+1}+\nabla z\big)\big]\cdot\nabla v_h\,d\bm{x} \\
&\qquad\qquad\qquad\qquad\qquad\qquad\qquad+\int_{\Gamma}q_{\Gamma}v_h\,d\bm{x}=0,\quad\forall v_h\in{\mathcal{V}}_h,\\
&\int_{\Omega}\psi_h^{n+1}w_h\,d\bm{x}-\int_{\Omega}h_{cap}\big(J(S_h^n)+J_{\delta}^{\prime}(S_h^n)(\widehat{S}_h^{n+1}-S_h^n)\big)w_h\,d\bm{x}=0,\quad\forall w_h\in\mathcal{V}_h.\\
\end{aligned}
\label{RichardsEqMixed1SemiScheme}
\end{equation}
Then, find $S_h^{n+1}\in\mathcal{S}_h$, solution of $(\textbf{P})$. To the best of our knowledge, the linear second-order time stepping mixed finite element method {\eqref{RichardsEqMixed1SemiScheme} for the Richards equation {\eqref{RichardsEq}} has not been explored before. The schemes \eqref{RichardsEqMixed1ImplScheme} and $(\mathbf{P})$, and \eqref{RichardsEqMixed1SemiScheme} and $(\mathbf{P})$  will be referred to as the implicit $(S,\psi)$-scheme and the semi-implicit $(S,\psi)$-scheme, respectively.

\begin{remark}
All the numerical algorithms are two steps schemes and therefore, they require the use of a starting procedure to obtain an approximation of the solution at the first step.
\end{remark}
\begin{remark}
We will use the Newton method for the linearization of the implicit $(S,\psi)$-scheme. The use of this method for the implicit $S$- and $(S,\bm{q})$-schemes requires the computation of the second derivative of the Leverett $J$-function. The resulting $J$-functions from many of the constitutive relationships, in particular  equations \eqref{HaverkampJfunction} and \eqref{vanGenuchtenJfunction} used in this study,  are not sufficiently differentiable (e.g., $J^{\prime\prime}$ is not defined at $S=0$ and $S=1$). This leads to convergence issues when the Newton method is applied. Therefore, we will use the Picard method for the linearization of the implicit $S$- and $(S,\bm{q})$-schemes. The stopping criteria for the  iterative methods is $\Vert u_h^{n,k+1}-u_h^{n,k}\Vert_{L^2(\Omega)} \leqslant \epsilon$, where $0<\epsilon\ll1$ is the chosen tolerance, $u_h^{n,k+1}$ and $u_h^{n,k}$ are the solutions at the iterations $k+1$ and $k$, respectively.
\end{remark}

\begin{remark}
For each of the above schemes, one can easily show the conservation of the variable $\widehat{S}_h$ over time by taking $v_h$ to be constant on the domain $\Omega$ and considering no-flux boundary conditions, i.e., setting $q_{\Gamma}=0$. However, solving directly the additional optimization problem \textbf{(P)} by an Uzawa algorithm or a simple projection technique could  lead to a defect of conservation for the variable $S_h$. We have quantified this error in a previous work on the Cahn-Hilliard equation
(see Figure~1 in \cite{Keita2020}). We also proposed in \cite{Keita2020} new conservative projection techniques for solving the optimization problem while maintaining the conservation. These techniques can be easily applied to the Richards equation to satisfy the conservation. We also compare the time evolution of the total mass of the exact and numerical solutions in the test case in section~\ref{SectionExactSol}. Our tests show a good agreement for the total mass.
\end{remark}

\begin{remark}
The discrete time-derivative introduces a ``distributed'' mass matrix (the term in $\hat{S}_h^{n+1}$) in the global algebraic system for all our finite element methods. As noted in \cite{Celia1990,Karthikeyan2011}, this distributed mass matrix (without lumping) leads to a numerical method that no longer satisfies a discrete maximum principle when the time step is taken smaller than a prescribed threshold value. Without discrete maximum principle, the solution may suddenly become non-monotone with unphysical oscillations from one spatial node to the next. This loss of discrete maximum principle for time steps below a threshold was proven for linear parabolic equations, such as the heat equation \cite{Rank1983}. For Richards equation, the same threshold on the time step holds in numerical tests for the loss of monotonicity of the numerical solutions, though no proof is available even for simple time-stepping schemes such as backward Euler method \cite{Karthikeyan2011}. A common way employed in finite element methods for the Richards equation to recover a discrete maximum principle is the mass lumping technique. An early study on the subject \cite{Celia1990} illustrated that mass lumping can effectively prevent nonoscillatory solutions, by comparing solutions with and without mass lumping on the time derivative. Several forms of mass lumping have been proposed in the standard and mixed finite element methods for unsaturated flow problems \cite{Belfort2009,Cooley1983,Milly1985,Younes2006}. A simple mass lumping technique, based on the row-sum of all the terms in the time derivative, is used in our study and was effective at removing oscillations that may occur near the drainage front. However, proving a discrete maximum principle may turn out to be challenging for our two-step methods.
\end{remark}

\section{Numerical tests}
\label{SectNumTests}
A series of numerical problems is given in this section to test and compare  the performance of the numerical methods discussed in section~\ref{SectNumMethods}. All the test cases are in two dimensions. We use the sparse linear solver UMFPACK which provides a relatively efficient solution.  Computations and results could thus carry over to 3D in a straightforward manner. However, this would require modifications in the implementation of the schemes. All the algorithms for the test cases are implemented using the FreeFem++ software \cite{FreeFem,FreeFemm}. 

\subsection{Manufactured solutions}
In this subsection, we numerically investigate the rates of convergence in space and time of the different schemes. Following the framework presented in \cite{Baron2017,SochalaThesis2008}, we employ a manufactured solution to test the methods. The idea is to use an arbitrary sufficiently differentiable function as exact solution with a source term for the computation of numerical errors. We perform numerical simulations using the computational domain $\Omega=[0\,cm,4\,cm]\times[0\,cm,20\,cm]$ and the final computational time $T=120\,s$. We consider the manufactured solution given by 
\begin{equation}
\psi(x,z,t)=20.4\tanh\big(0.5(z+\frac{t}{12}-15)\big)+c, \quad c\in \mathbb{R},
\label{ManSol}
\end{equation}
where $c$ is a constant satisfying $c\leqslant-20.4$ to ensure the non-positivity of $\psi$.  The Richards equation \eqref{RichardsEq} is then solved with Dirichlet boundary conditions and a source term such that \eqref{ManSol} is a solution, i.e., the closed-form expression \eqref{ManSol} is used to calculate the adequate source term as well as the initial and boundary conditions. We use the following relationships between the relative permeability, pressure head and water content \cite{Haverkamp1977}:
\begin{equation}
\theta=\theta_r+\dfrac{\theta_s-\theta_r}{1+\vert\tilde{\alpha}\psi\vert^{\beta}},\quad K_r=\dfrac{1}{1+\vert\tilde{A}\psi\vert^{\gamma}},
\label{HaverkampRelations}
\end{equation}
where
\begin{equation}
\begin{aligned}
&\theta_s=0.287,\qquad \tilde{\alpha}=0.0271\,cm^{-1},\qquad K_s=9.44\times 10^{-3}\,cm.s^{-1},\qquad \gamma=4.74,  \\
&\theta_r=0.075,\qquad \beta=3.96,\qquad\qquad\qquad \tilde{A}=0.0524\,cm^{-1}.
\end{aligned}
\label{HaverkampParametersValues}
\end{equation}
Using \eqref{HaverkampRelations}, the pressure head $\psi$ can be expressed as in \eqref{WaterSectionHead}, where
\begin{equation}
h_{cap}(\bm{x})=\dfrac{1}{\tilde{\alpha}} \quad\text{and}\quad J(S)=-(S^{-1}-1)^{1/\beta}.
\label{HaverkampJfunction}
\end{equation}
The saturation is given by the relation
\begin{equation}
S=\dfrac{1}{1+\vert\tilde{\alpha}\psi\vert^{\beta}},
\end{equation}
that one uses to calculate the initial and boundary conditions for $S$ from the conditions on $\psi$.

\subsubsection{Numerical example without regularization of the Leverett $J$-function }
For the first test case, we take $c=-41.1$ in \eqref{ManSol}. The profile of the manufactured solutions is shown in Figure~\ref{Fig:ManFacSol1}.
\begin{figure}[!h]
\centering
\begin{tabular}{cc}
\includegraphics[scale=0.55]{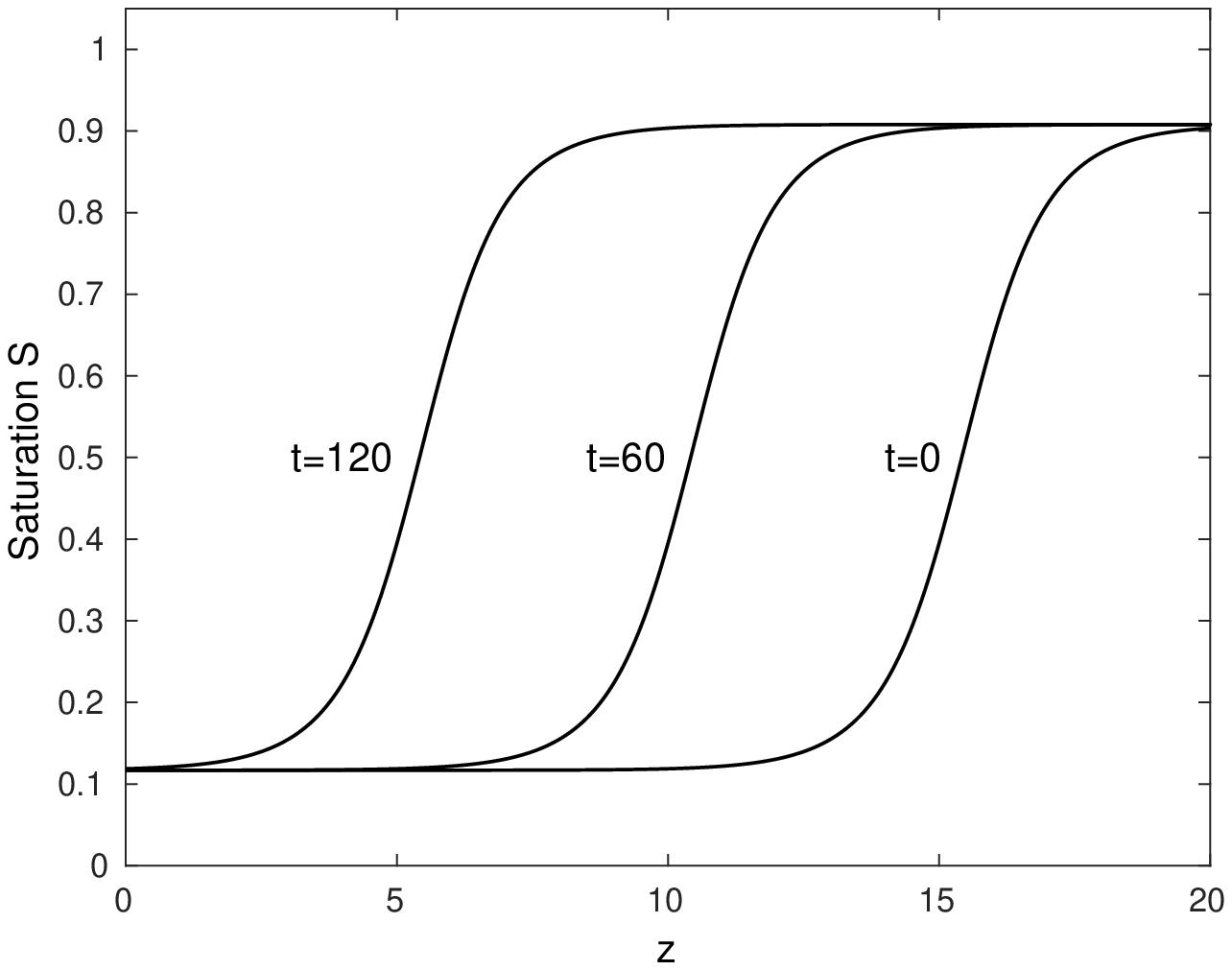}&\includegraphics[scale=0.55]{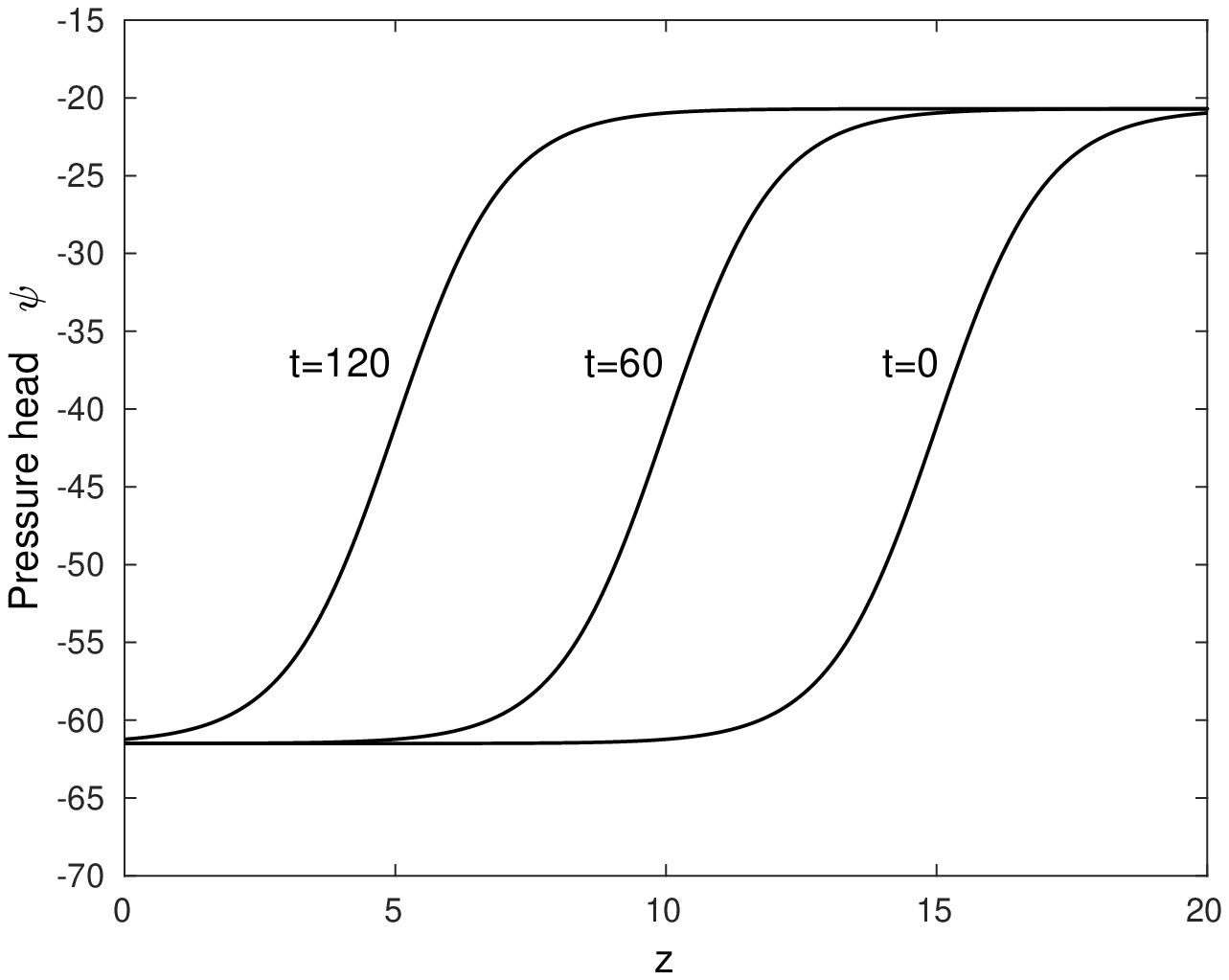}
\end{tabular}
\caption{Test case without regularization of the Leverett $J$-function: Profile of the manufactured solution \eqref{ManSol}  at different times for $c=-41.1$.}
\label{Fig:ManFacSol1}
\end{figure}
The saturation $S$ is far from $1$ for all $t\in[0,T]$. Therefore, $J_{\delta}^{\prime}(S)$ in \eqref{RegularizedJfunction} is exactly given by $J^{\prime}(S)$. Note also that no specific treatment is needed in the numerical discretization of the source term. We use a constant time step $\Delta t=0.2$ for all the schemes. The tolerance for the stopping criteria for the iterative schemes is $\epsilon=10^{-5}$. For each scheme, the $L^2$- and $H^1$-error on the effective saturation $S_h$ and the pressure head $\psi_h$ are computed on different meshes. The numerical results are presented in Figure~\ref{Fig:ErrorEspaceManSolTest1}. We can see that all the schemes have the same order of convergence around $k=1$ for the $H^1$-error on the variables $S_h$ and $\psi_h$. For the $L^2$-error, an order of convergence varying around $k=2$ is obtained on both variables $S_h$ and $\psi_h$ with each scheme. This order is slightly larger with the semi-implicit schemes. There is no significant difference on the errors produced by the different schemes. However, the error on $\psi_h$ is larger than the one on $S_h$. The computational times for each scheme with respect to the mesh refinement are plotted in Figure~\ref{Fig:CPUtimeTest1Graphs}. Here, the semi-implicit schemes are the winner in terms of computational cost since we use the same time step for all the schemes. Larger time step could be used for the implicit schemes as we will see in the next test cases. With respect to the formulations, our numerical simulations show that the $(S,\psi)$-schemes are faster than the $S$-schemes which are faster than the $(S,\bm{q})$-schemes. This shows the efficiency of the mixed formulation based on the saturation and the pressure head with respect to the other formulations.
\begin{figure}[!h]
\centering
\begin{tabular}{cc}
\includegraphics[scale=0.55]{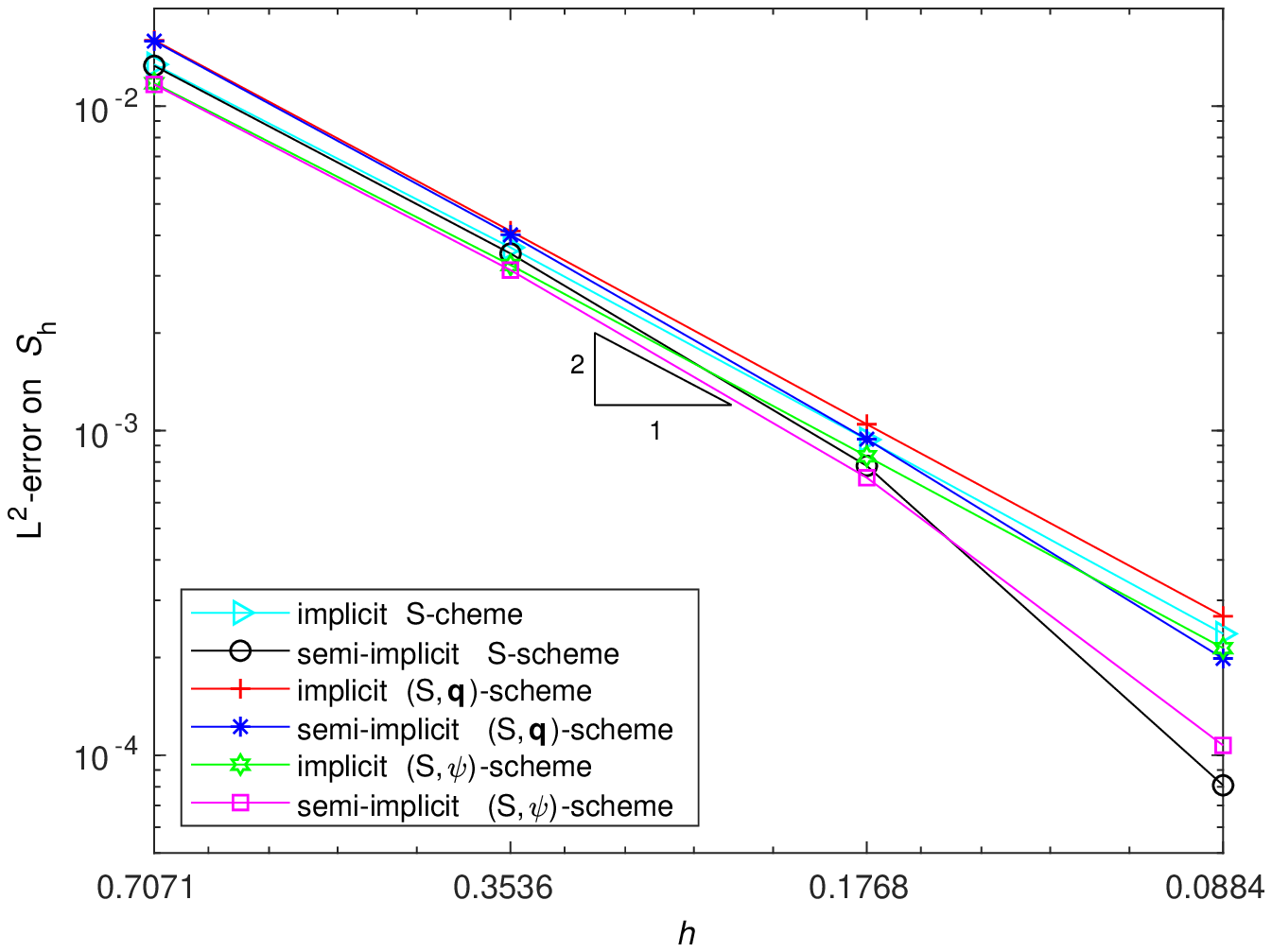}&\includegraphics[scale=0.55]{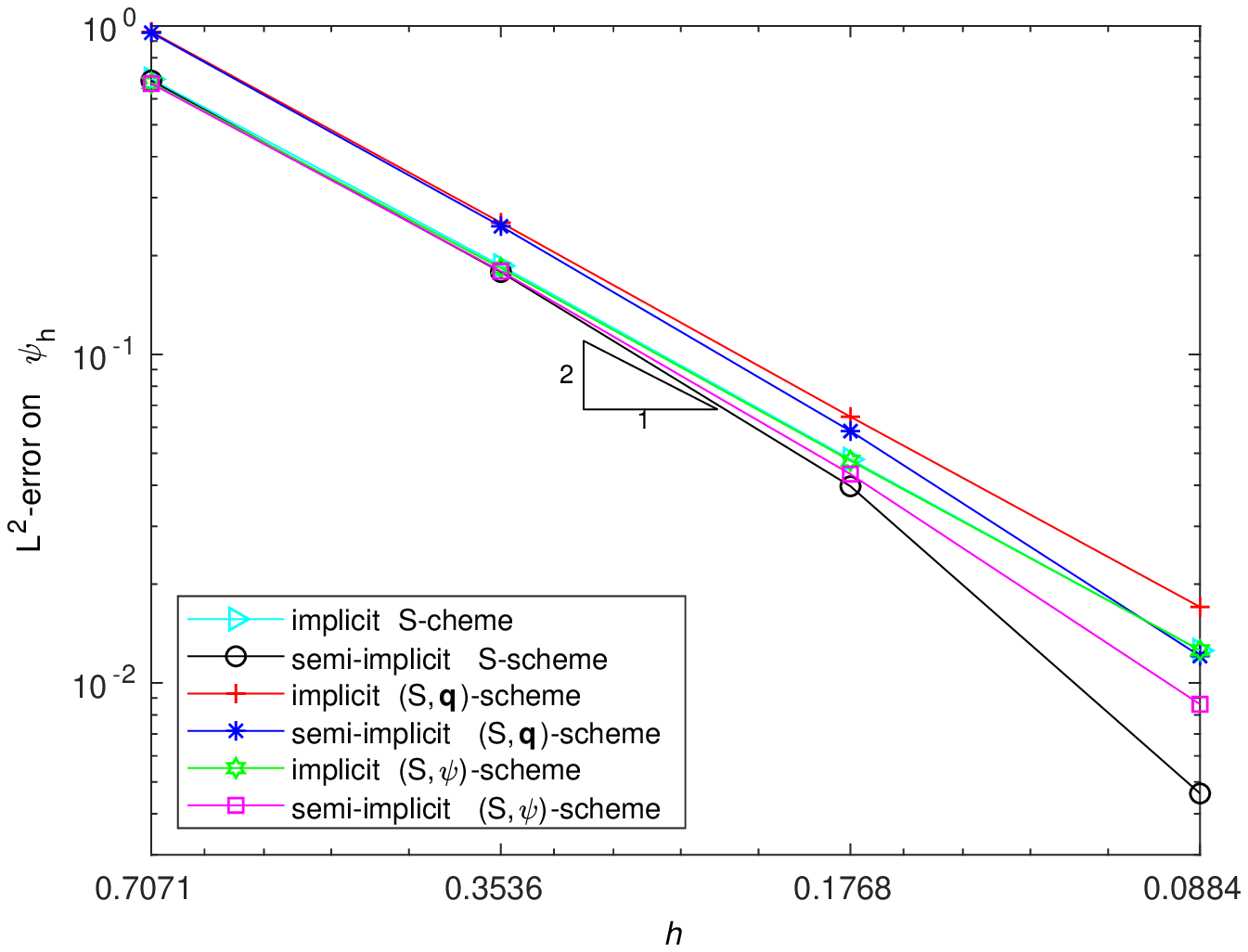}\\
\includegraphics[scale=0.55]{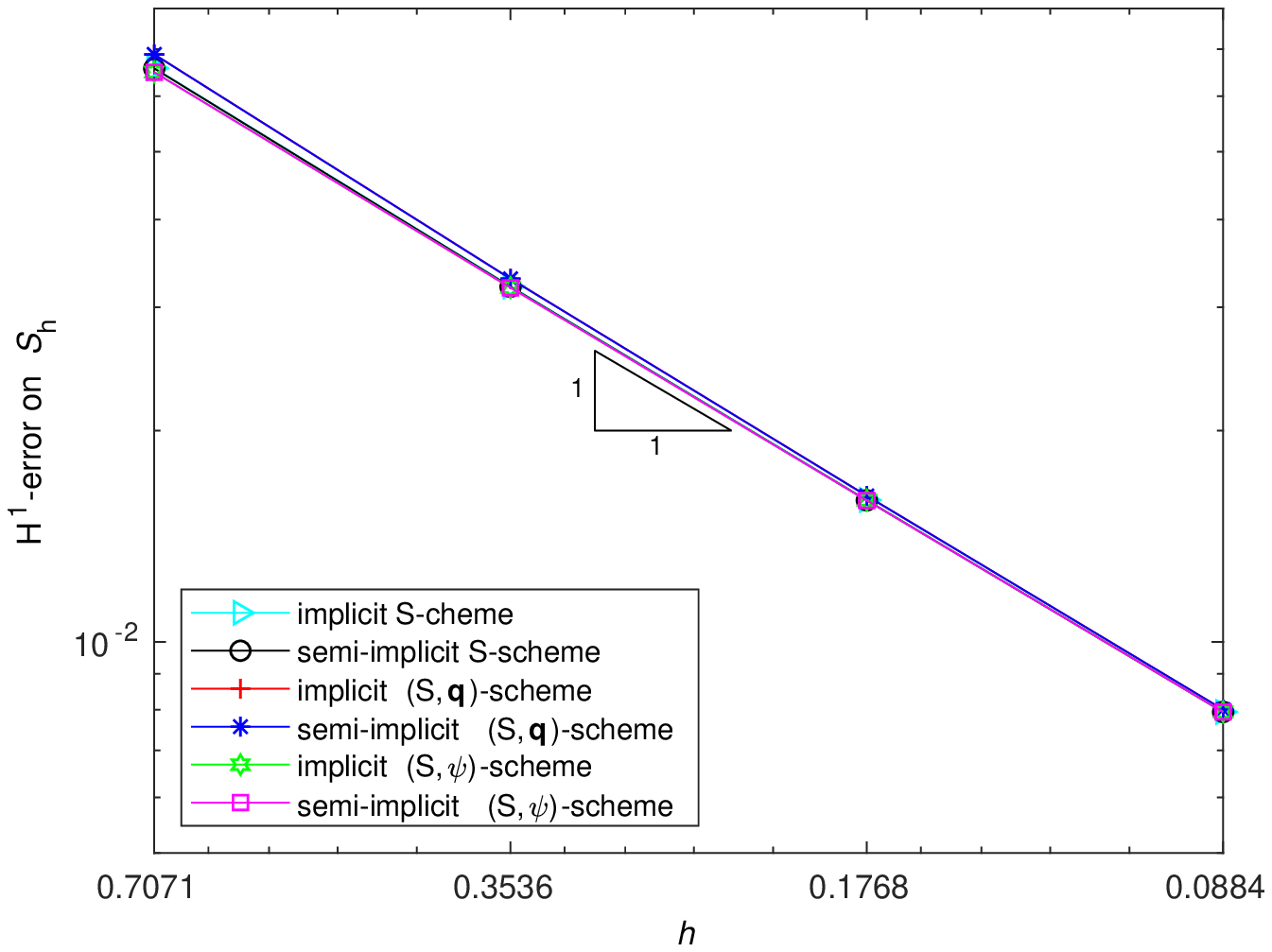}&\includegraphics[scale=0.55]{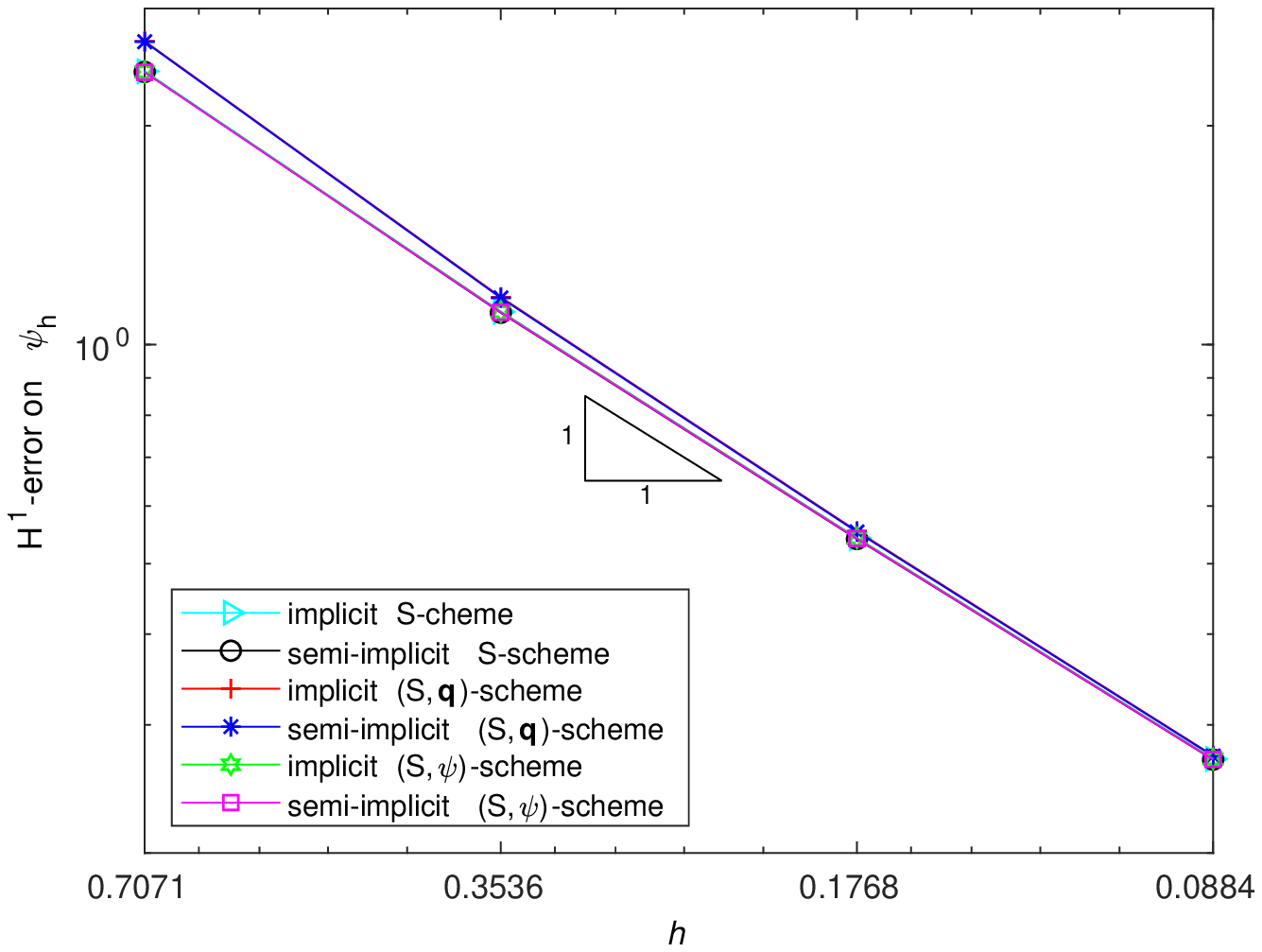}\\
\end{tabular}
\caption{Manufactured solution --- Test case without regularization of the Leverett $J$-function: $L^2$- and $H^1$-error on $S_h$ and $\psi_h$ as a function of the mesh size $h$ for each scheme with respect to the manufactured solution \eqref{ManSol} for $c=-41.1$.}
\label{Fig:ErrorEspaceManSolTest1}
\end{figure}
\begin{figure}[!h]
\centering
\includegraphics[scale=0.55]{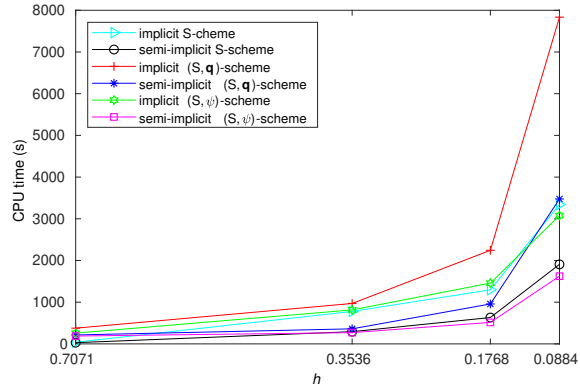}
\caption{Manufactured solution --- Test case without regularization of the Leverett $J$-function: Computational time in seconde (s) as a function of the mesh size $h$ for each scheme.}
\label{Fig:CPUtimeTest1Graphs}
\end{figure}

For the error and convergence analysis in time, we compute a reference solution $(S_h^{\text{ref}}$,$\psi_h^{\text{ref}})$ with the semi-implicit $(S,\psi)$-scheme on a computational grid of size $h=0.176$ ($32\times 160$ elements) by using a small time step $\Delta t= 0.001\,s$ to get a good approximation of the semi-discretized solution of the problem. We then analyze the numerical convergence in time by varying the time step. For each time step, the $L^2$-error on $S_h$ and $\psi_h$ with respect to the reference solution are recorded at the final time $T=120\,s$. Figure~\ref{Fig:TimeOrderTest1Graphs} shows the errors as function of the time step for each scheme. An order of convergence in time varying slightly around $p=2$ is obtained on the variables $S_h$ and $\psi_h$  with each scheme. The schemes produce almost the same error for each variable. However, the error on $\psi_h$ is larger than the error on $S_h$.
\begin{figure}[!h]
\centering
\begin{tabular}{cc}
\includegraphics[scale=0.55]{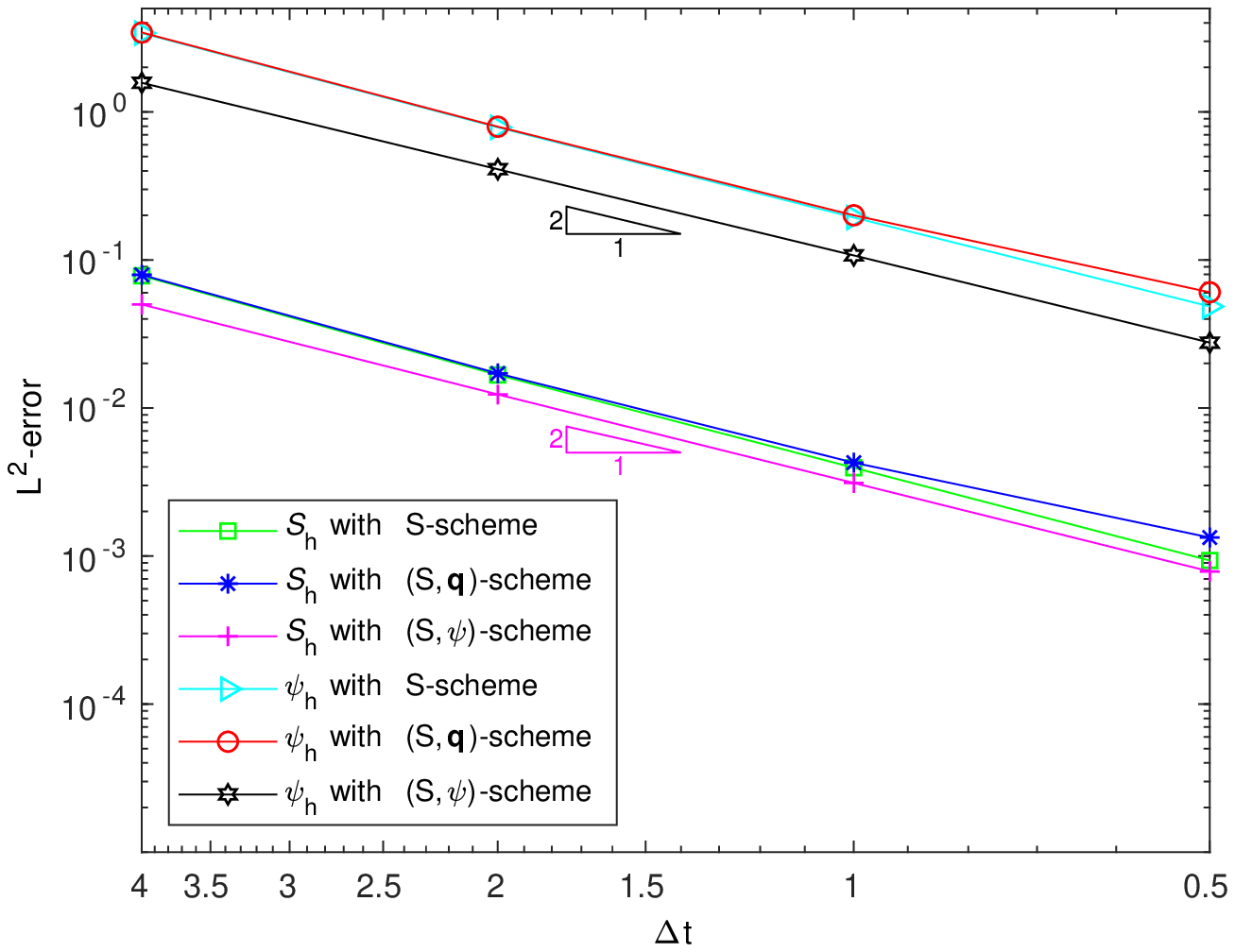}&\includegraphics[scale=0.55]{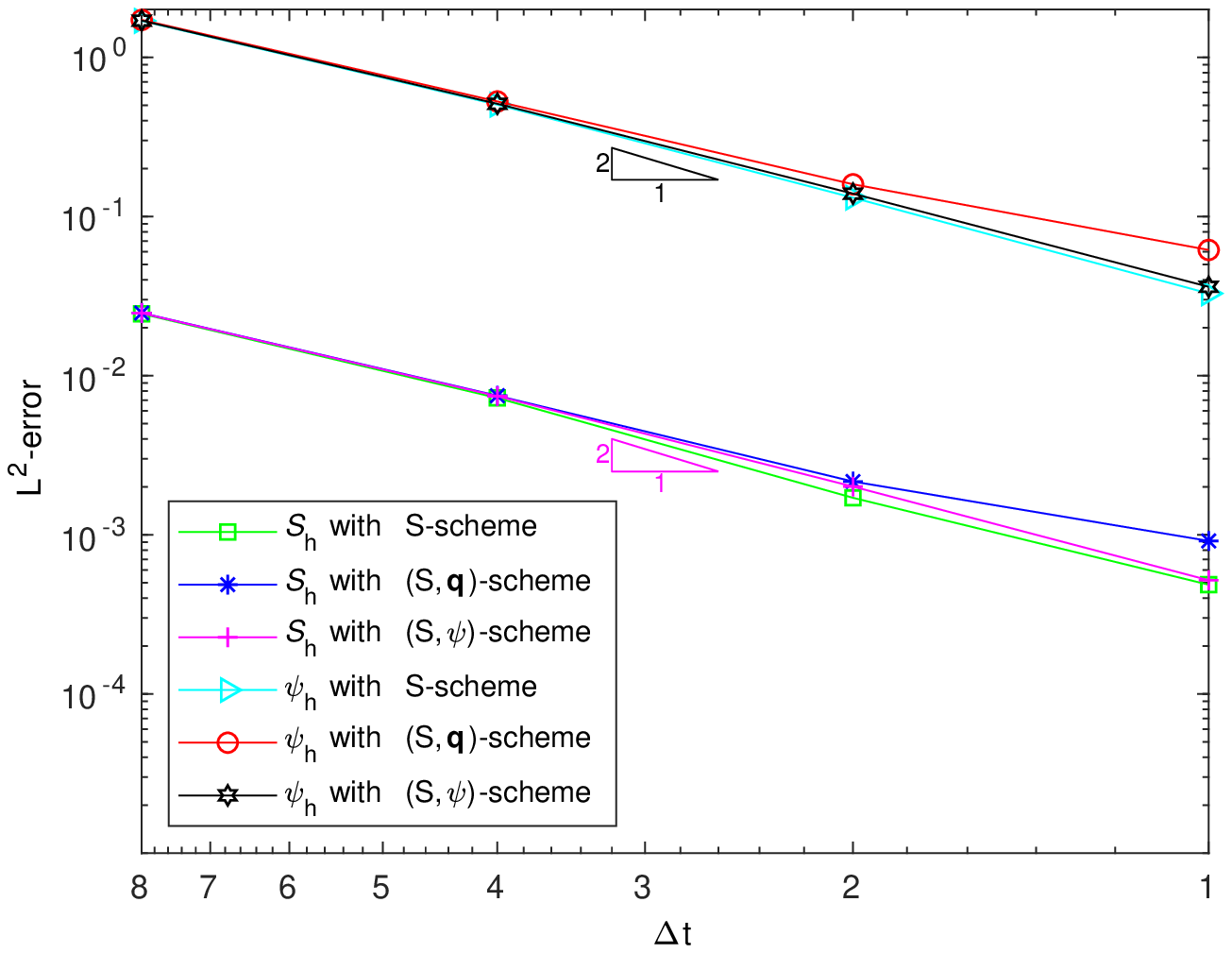}\\
Semi-implicit schemes & Implicit schemes
\end{tabular}
\caption{Manufactured solution --- Test case without regularization of the Leverett $J$-function: $L^2$-error on $S_h$ and $\psi_h$ as a function of the time step $\Delta t$ for each scheme.}
\label{Fig:TimeOrderTest1Graphs}
\end{figure}

Most numerical methods for the Richards equation \eqref{RichardsEqMixedForm} are based on the classical backward Euler method \cite{Arbogast1993,Arbogast1996,Bause2008,Bause2004,Belfort2009,Eymard1999,Forsyth1997,Oulhaj2018,Radu2004,Radu2014,Schneid2004}. For the purpose of comparison of this standard technique with the proposed schemes, we consider the finite element methods for the formulations \eqref{RichardsEqMod} and \eqref{RichardsEqMixed2}, respectively, using the classical backward Euler method for the temporal discretization. These schemes are named backward Euler $S$-scheme for the formulation \eqref{RichardsEqMod} and backward Euler $(S,\bm{q})$-scheme for the formulation \eqref{RichardsEqMixed2}, and are compared against the closest among the proposed schemes, namely the $S$-schemes. Table~\ref{Tab:Comparison} shows the $L^2$-error on $S_h$ and $\psi_h$, and the computational times for the solutions computed with the $S$-schemes, backward Euler $S$-scheme and backward Euler $(S,\bm{q})$-scheme, respectively, on the mesh $32\times 160$ elements with the time step $\Delta t=4\,s$. Our results show that the errors with classical methods are at least seven times larger than the errors with the proposed implicit $S$-scheme. The proposed semi-implicit $S$-scheme produces almost the same errors as the standard methods but with a computational time at least six times smaller. By reducing the time step, the semi-implicit $S$-scheme becomes more accurate and efficient than the standard methods. This confirms the efficiency of the proposed techniques with respect to the standard methods.

\begin{table}
\centering
\begin{tabular}{|c|c|c|c|}
\hline
Methods  & $L^2$-error on $S_h$ & $L^2$-error on $\psi_h$ & CPU(s)\\ \hline
Implicit $S$-scheme &$0.0075831$& $0.53518$ & $237$ \\ \hline
Semi-implicit $S$-scheme & $0.0775579$ & $3.37313$ & $040$ \\ \hline
Backward Euler $S$-scheme & $0.0603825$ & $3.84694$ & $262$\\ \hline
Backward Euler $(S,\bm{q})$-scheme & $0.0595938$ & $3.80333$ & $496$ \\ \hline
\end{tabular}
\caption{Manufactured solution --- Test case without regularization of the Leverett $J$-function: Comparison of errors and computational times between standard and proposed methods.}
\label{Tab:Comparison}
\end{table}

\subsubsection{Numerical example with regularization of the Leverett $J$-function}
For the second test case,  we set $c=-20.4$ in \eqref{ManSol}. The profile of the manufactured solution is shown in Figure~\ref{Fig:ManFacSol2}.
\begin{figure}[!h]
\centering
\begin{tabular}{cc}
\includegraphics[scale=0.55]{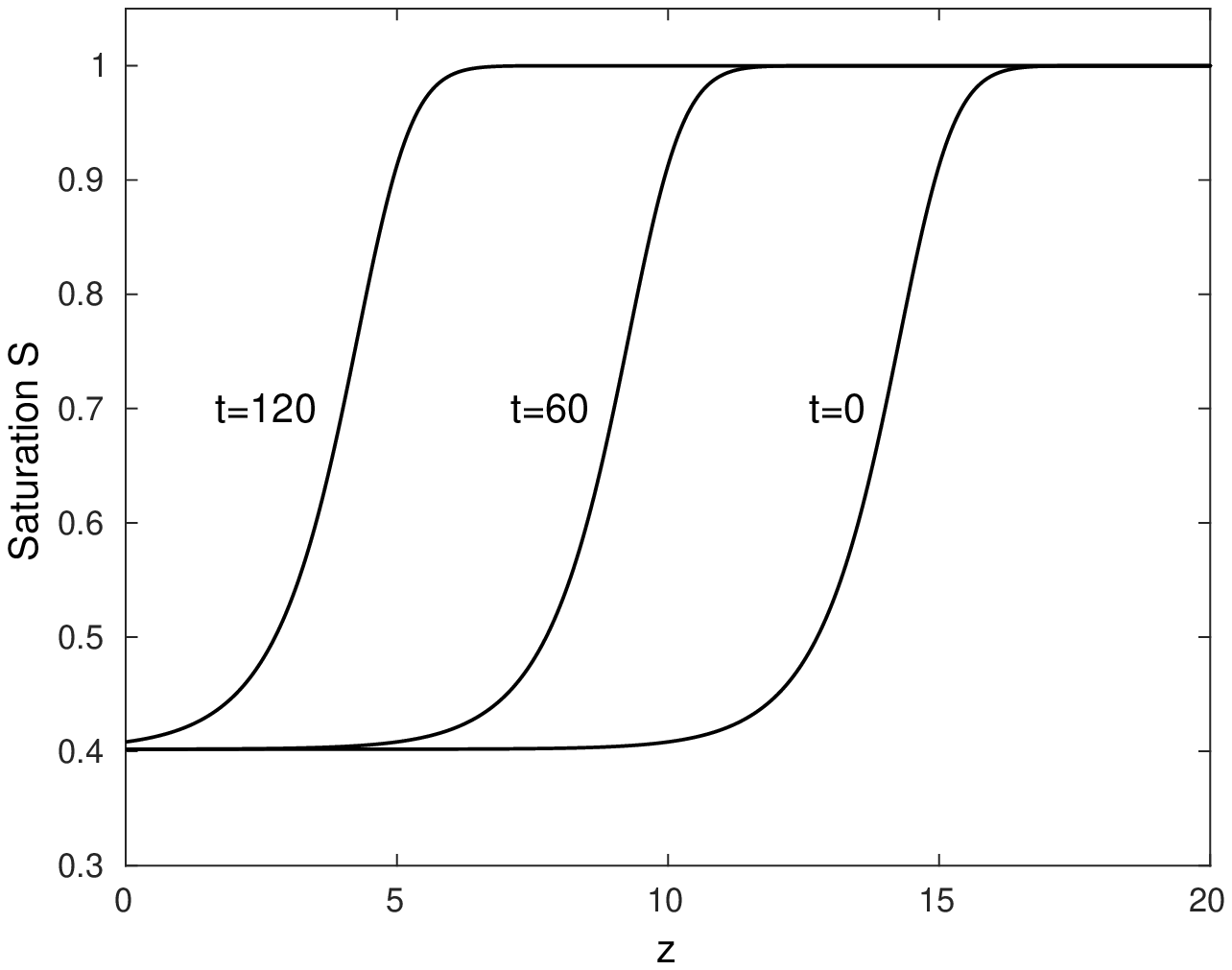}&\includegraphics[scale=0.55]{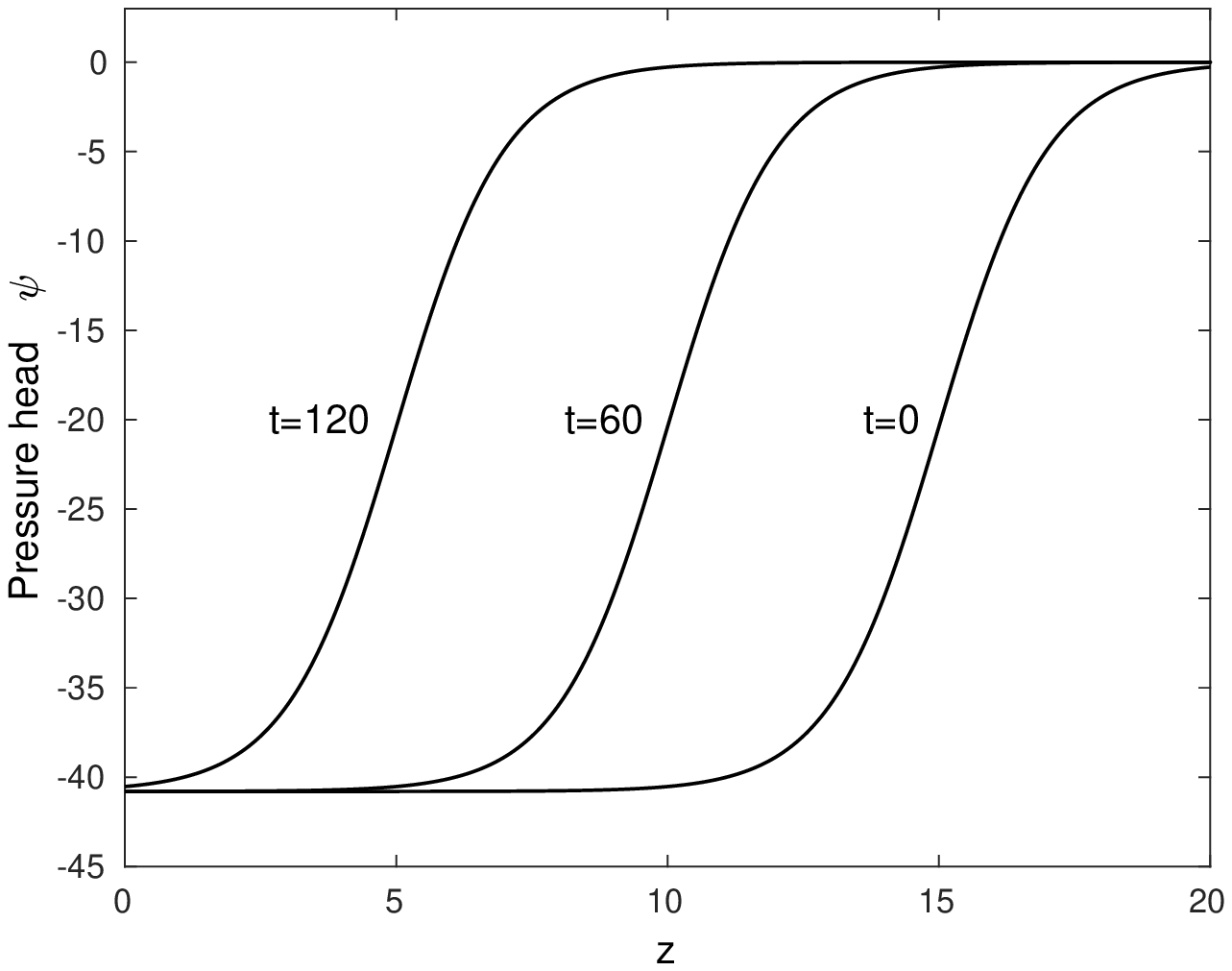}\\
\end{tabular}
\caption{Test case with regularization of the Leverett $J$-function: Profile of the manufactured solution \eqref{ManSol}  at different times for $c=-20.4$.}
\label{Fig:ManFacSol2}
\end{figure}
The purpose of this test case is to verify the convergence behavior of the different schemes with the regularization \eqref{RegularizedJfunction} when the saturation solution $S$ reaches the value of one in some parts of the domain. We have noticed in practice that the $(S,\psi)$-schemes allow smaller values of $\delta$ than the $S$- and $(S,\bm{q})$-schemes. With respect to the temporal discretization, the implicit $S$- and $(S,\bm{q})$-schemes allow smaller values of $\delta$ than their corresponding semi-implicit versions. We take $\delta=10^{-3}$ for the semi-implicit $S$- and $(S,\bm{q})$-schemes, $\delta=10^{-6}$ for the implicit $S$- and $(S,\bm{q})$-schemes, and $\delta=10^{-10}$ for the $(S,\psi)$-schemes. The tolerance for the iterative methods used for the linearization of the implicit schemes is set to $5\times10^{-4}$. For each scheme, the $L^2$- and $H^1$-error on $S_h$ and $\psi_h$ on different meshes are presented in Figure~\ref{Fig:ErrorEspaceManSolTest2} and the computational time for each scheme with respect to the mesh refinement is plotted in Figure~\ref{Fig:CPUtimeTest2Graphs}. 
\begin{figure}[!h]
\centering
\begin{tabular}{cc}
\includegraphics[scale=0.55]{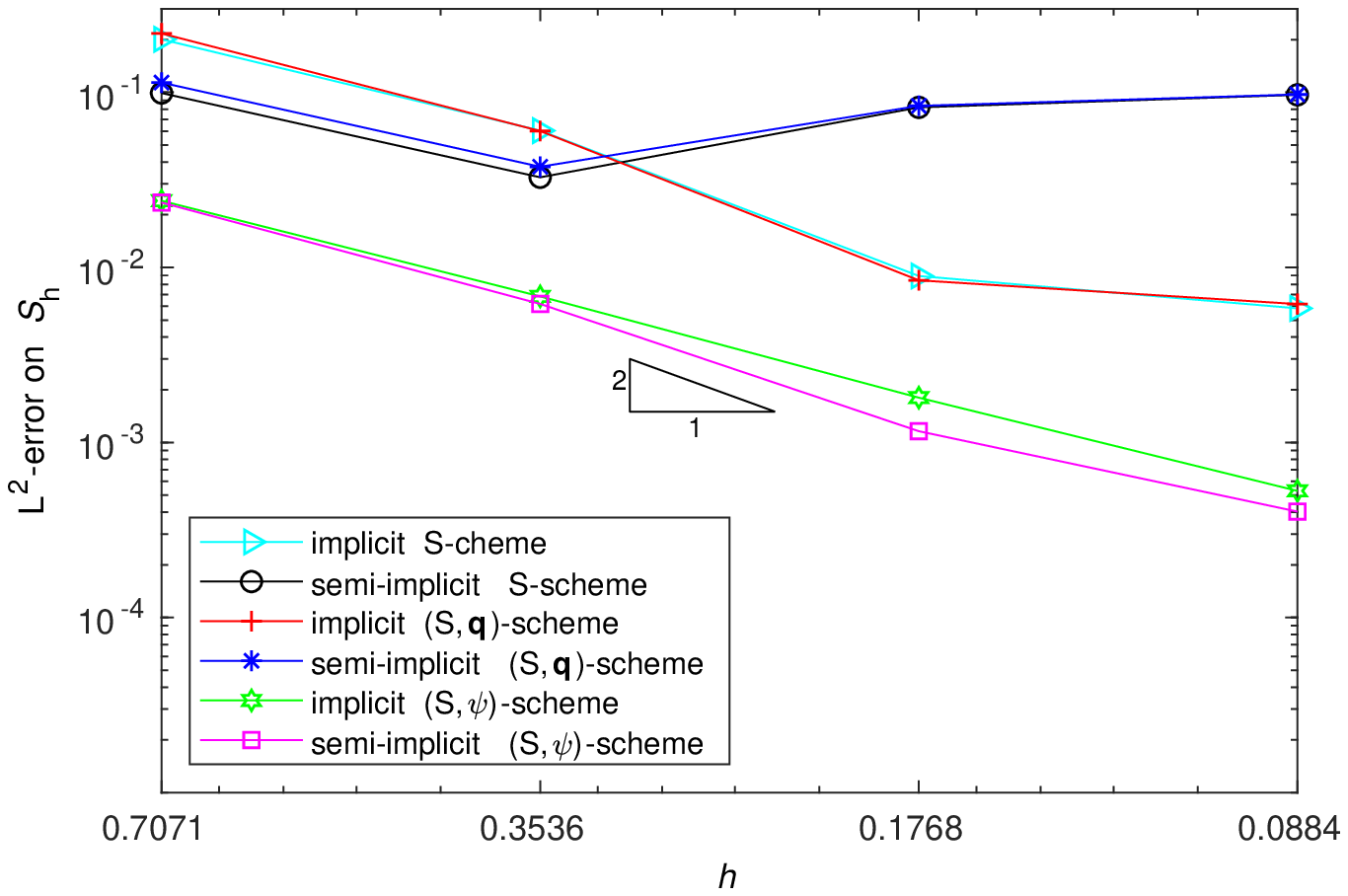}&\includegraphics[scale=0.55]{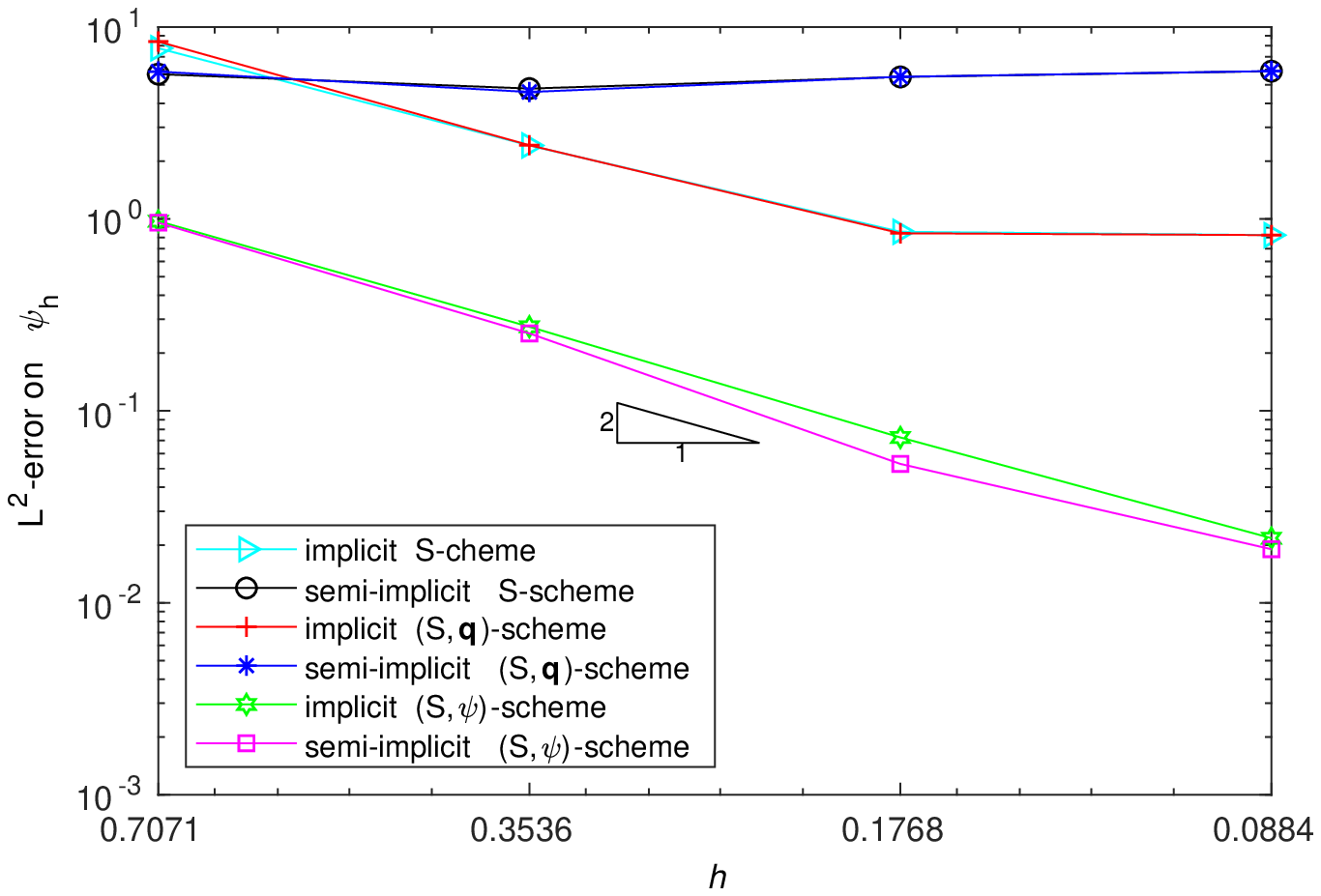}\\
\includegraphics[scale=0.55]{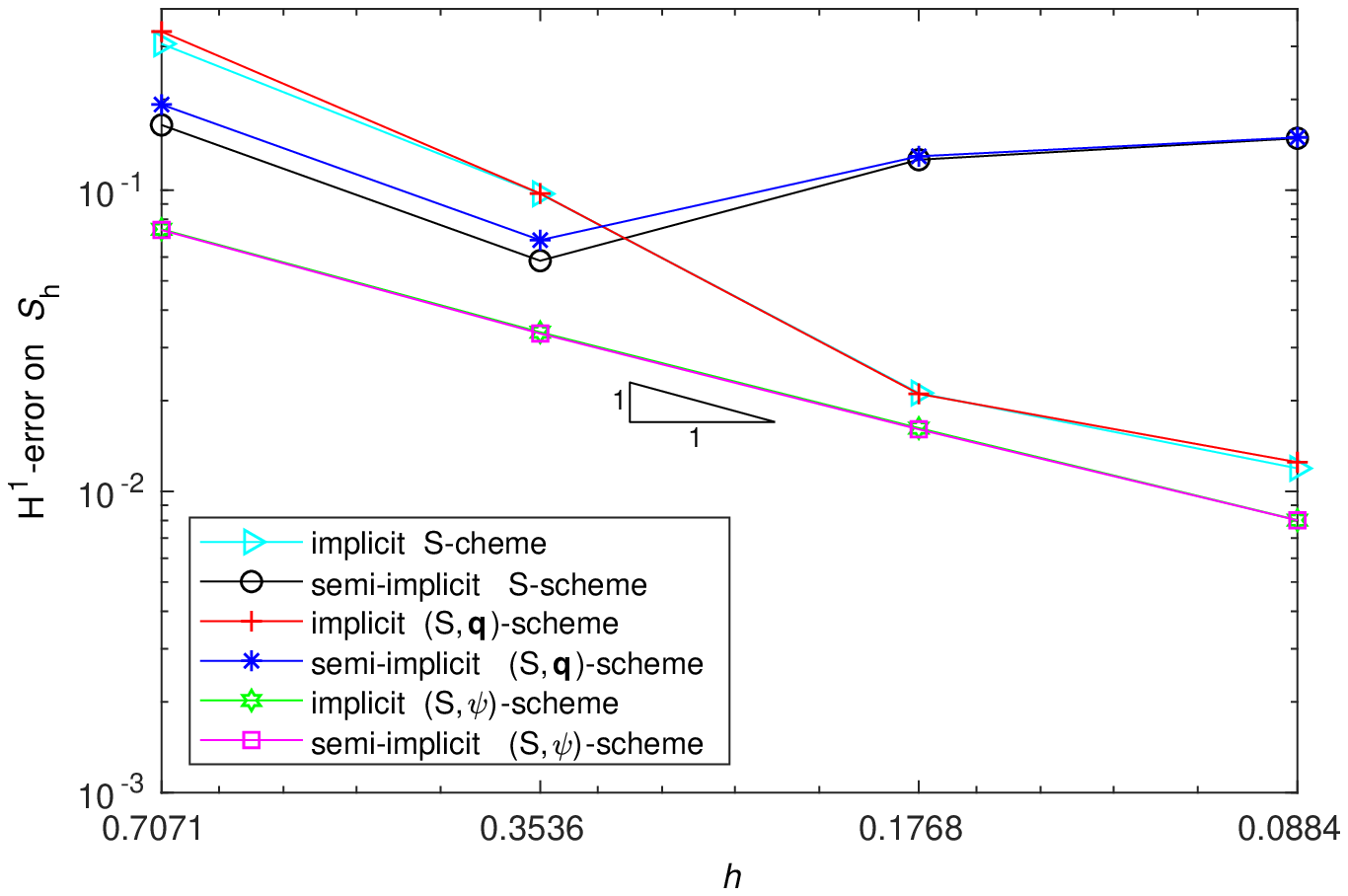}&\includegraphics[scale=0.55]{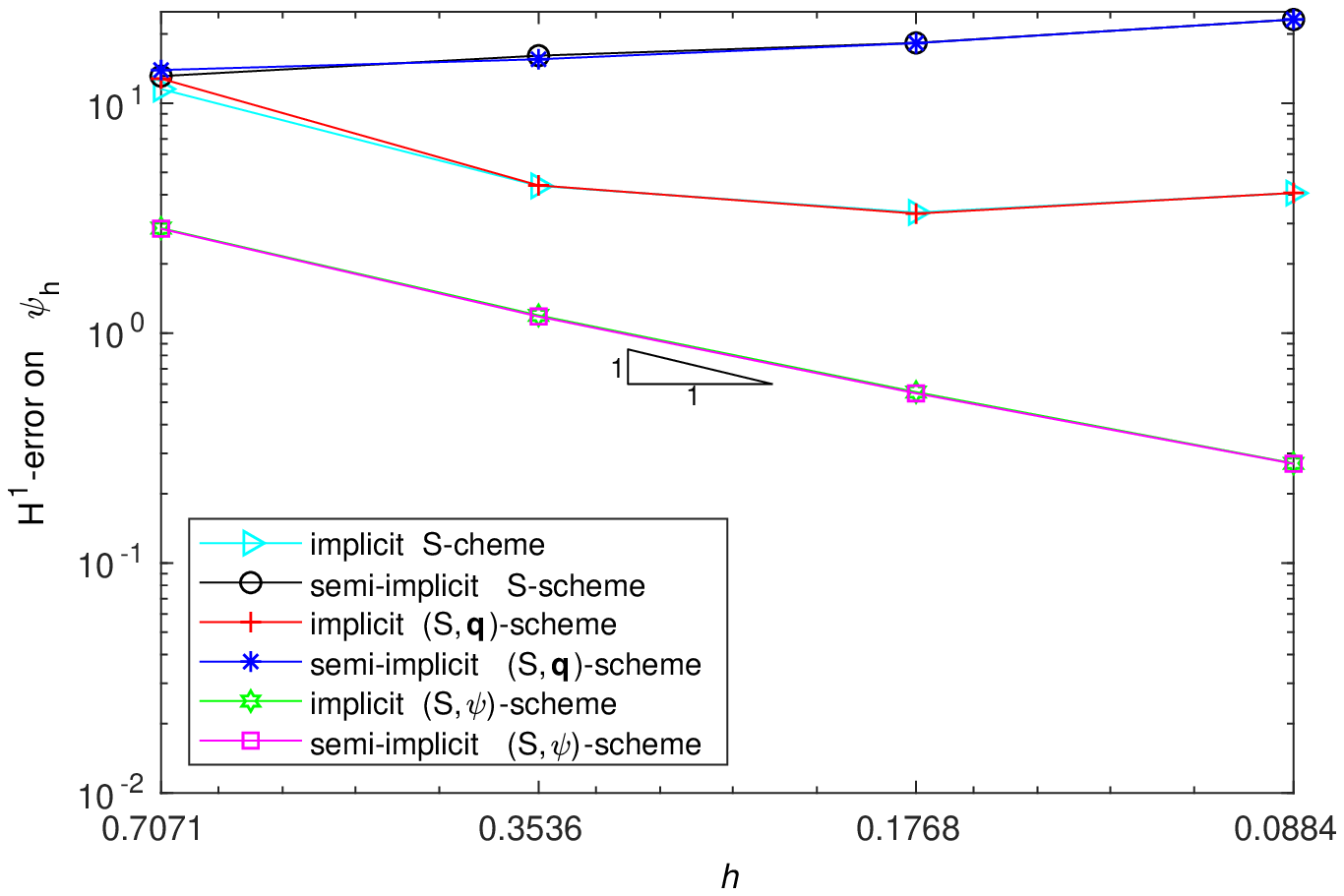}\\
\end{tabular}
\caption{Manufactured solution --- Test case with regularization of the Leverett $J$-function: $L^2$- and $H^1$-error on $S_h$ and  $\psi_h$ as a function of the mesh size $h$ for each scheme with respect to the manufactured solution \eqref{ManSol} for $c=-20.4$.}
\label{Fig:ErrorEspaceManSolTest2}
\end{figure}
\begin{figure}[!h]
\centering
\includegraphics[scale=0.55]{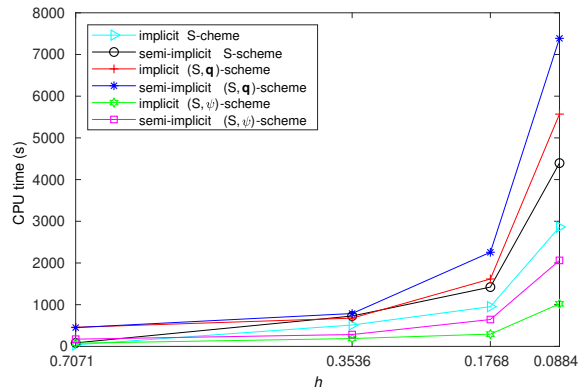}
\caption{Manufactured solution --- Test case with regularization of the Leverett $J$-function: Computational time in seconde (s) as a function of the mesh size $h$ for each scheme.}
\label{Fig:CPUtimeTest2Graphs}
\end{figure}
Numerical results show that the $(S,\psi)$-schemes produce almost the same error for a fixed mesh size, and converge with an order around $k=1$ for the $H^1$-error and an order about $k=2$ for the $L^2$-error on both variables $S_h$ and $\psi_h$. The implicit $S$- and $(S,\bm{q})$-schemes converge slowly with an order less than $1$ on fine meshes, while a lack of convergence is noticed with the semi-implicit $S$- and $(S,\bm{q})$-schemes. This shows again that the linear scheme using the mixed formulation consisting in simultaneously constructing approximations of the saturation and the pressure head with Lagrange linear element performs better in terms of accuracy and efficiency than the schemes using the mixed formulation based on the saturation and the flux using Raviart-Thomas elements. The implicit schemes perform better than the semi-implicit schemes in terms of efficiency. However, the iterative methods used for the linearization of the implicit schemes could not converge for some less regular solution as we will see with the test case in section {\ref{SectInfiltrationHeterogeneous}}.  In such cases, semi-implicit schemes could be an alternative.

We perform numerical tests to obtain the temporal order of convergence of the schemes. To this end, we compute a reference solution $(S_h^{\text{ref}}$,$\psi_h^{\text{ref}})$ with the implicit $(S,\psi)$-scheme using a computational grid of size $h=0.088$ ($64\times 320$ elements) and a time step $\Delta t= 0.1\,s$. We then analyze the numerical convergence in time by varying the time step for each scheme. For each time step, the $L^2$-error on $S_h$ and $\psi_h$ with respect to the reference solution is recorded at the final time $T=120\,s$. Figure~\ref{Fig:TimeOrderTest2Graphs} shows the $L^2$-error on $S_h$ and $\psi_h$ as a function of the time step for each scheme. An order of convergence in time around $p=2$ is obtained with the semi-implicit $(S,\psi)$-scheme and an order around $p=1.8$ is obtained with the implicit $(S,\psi)$-scheme on the variables $S_h$ and $\psi_h$. As for the spatial convergence, a lack of convergence in time is noticed with the $S$- and $(S,\bm{q})$-schemes. 

A second-order time-accuracy and optimal spatial convergence rates are obtained with the $(S,\psi)$-schemes for solutions with regularization, attesting the robustness of the $(S,\psi)$-schemes for solutions in which the saturation reaches the value of one in some parts of the domain. 
\begin{figure}[!h]
\centering
\begin{tabular}{cc}
\includegraphics[scale=0.55]{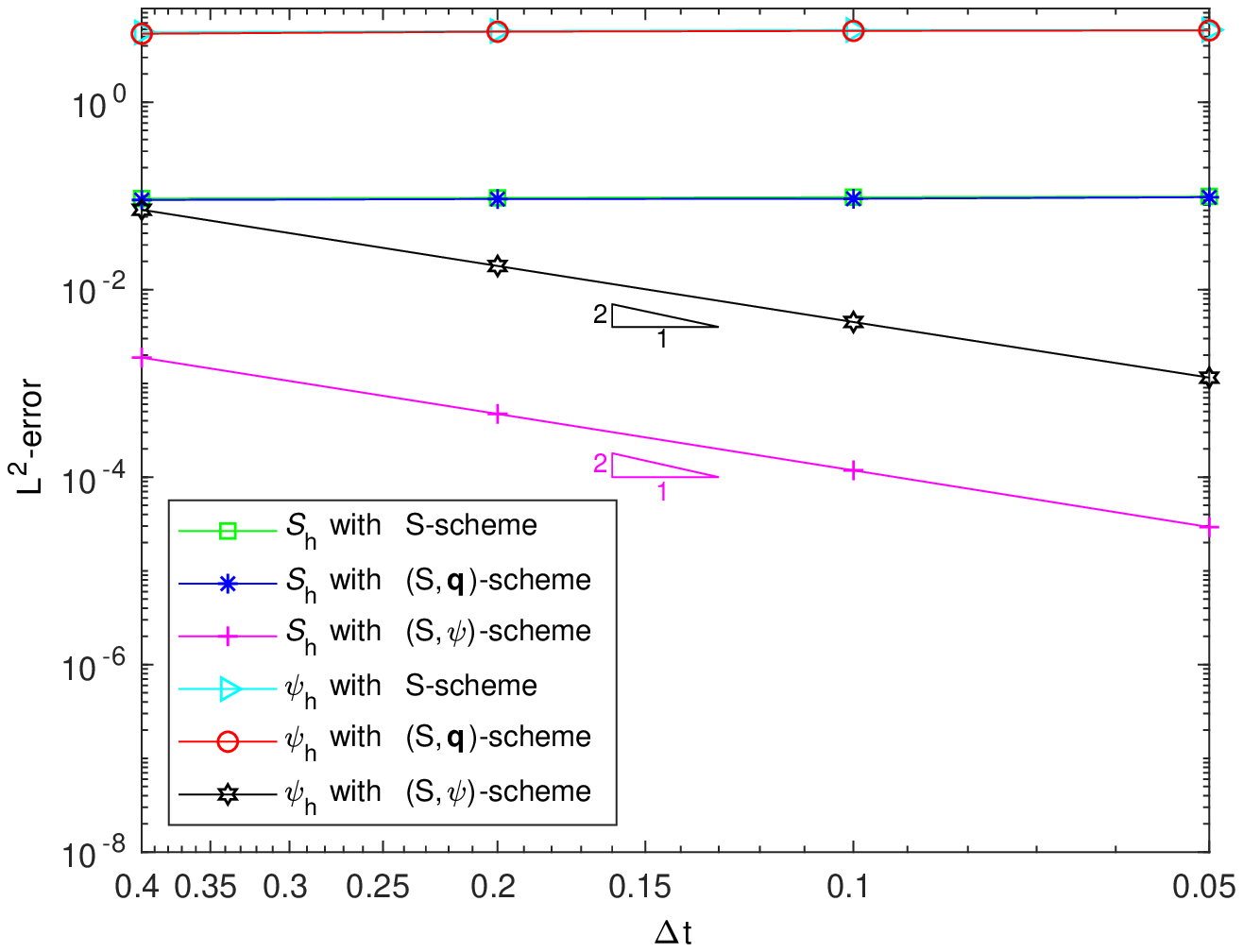}&\includegraphics[scale=0.55]{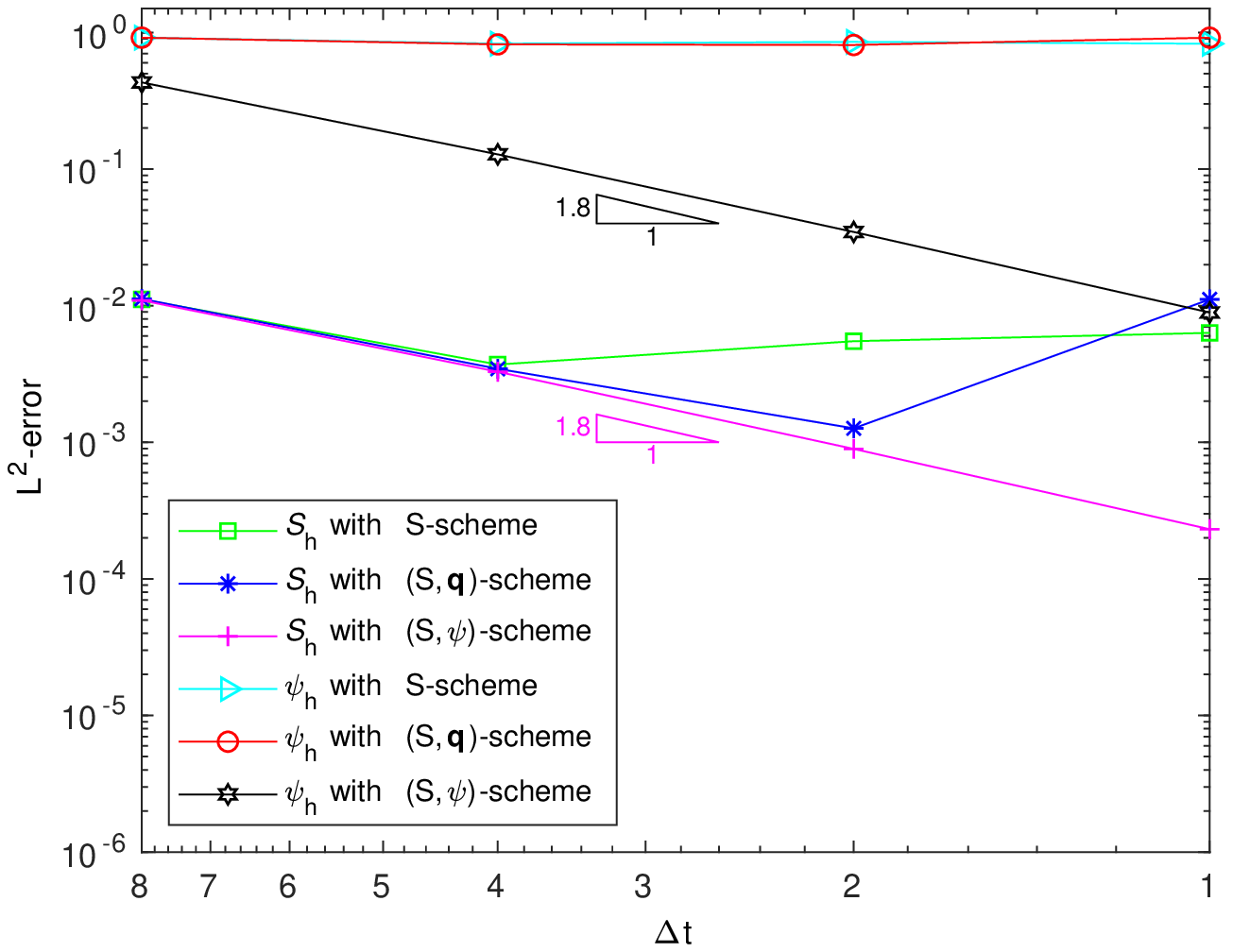}\\
Semi-implicit schemes & Implicit schemes
\end{tabular}
\caption{Manufactured solution --- Test case with regularization of the Leverett $J$-function: $L^2$-error on $S_h$ and $\psi_h$ as a function of the time step $\Delta t$ for each scheme.}
\label{Fig:TimeOrderTest2Graphs}
\end{figure}

\subsection{2-D analytical solution of the Green-Ampt problem}
\label{SectionExactSol}
In this section, we use an analytical solution to \eqref{RichardsEq} for checking the accuracy of our numerical schemes. Closed-form one-dimensional solutions for the Richards equation are available \cite{Broadbridge1988,Fujita1952,Parkin1995,Sander1991,Sander1988}.  Two- and three-dimensional analytical solutions have been derived in \cite{Tracy2006,Tracy2011}. Here, we use a two-dimensional analytical solution proposed in \cite{Tracy2011}.

Figure~\ref{Green-AmptGeometry}\textbf{(i)} shows a 2-D soil sample of dimensions $a\times L$.
\begin{figure}[!h]
\centering
\begin{minipage}[b]{0.45\linewidth}
\begin{tikzpicture}
\draw[black,thick] (0,0) -- (5,0) -- (5,5) -- (0,5) -- (0,0);
\draw[black,thick,-> ] (5.5, 2.8) -> (5.5,5);
\node[ ] at (5.5,2.5) {L};
\draw[black,thick,-> ] (5.5, 2.2) -> (5.5,0);
\draw[black,thick,-> ] (2.3, -0.5) -> (0,-0.5);
\node[ ] at (2.5,-0.5) {a};
\draw[black,thick,-> ] (2.7, -0.5) -> (5,-0.5);
\draw[black,thick,-> ] (1,6) -> (1,5.2);
\draw[black,thick,-> ] (1.5,6) -> (1.5,5.2);
\draw[black,thick,-> ] (2,6) -> (2,5.2);
\draw[black,thick,-> ] (2.5,6) -> (2.5,5.2);
\draw[black,thick,-> ] (3,6) -> (3,5.2);
\draw[black,thick,-> ] (3.5,6) -> (3.5,5.2);
\draw[black,thick,-> ] (4,6) -> (4,5.2);
\node[ ] at (2.5,6.5) {Rainfall};
\node[ ] at (2.5,-1.5) {\textbf{(i)}};
\end{tikzpicture}
\end{minipage}%
\begin{minipage}[b]{0.55\linewidth}
\color{darkgray}
\centering
\includegraphics[scale=0.65]{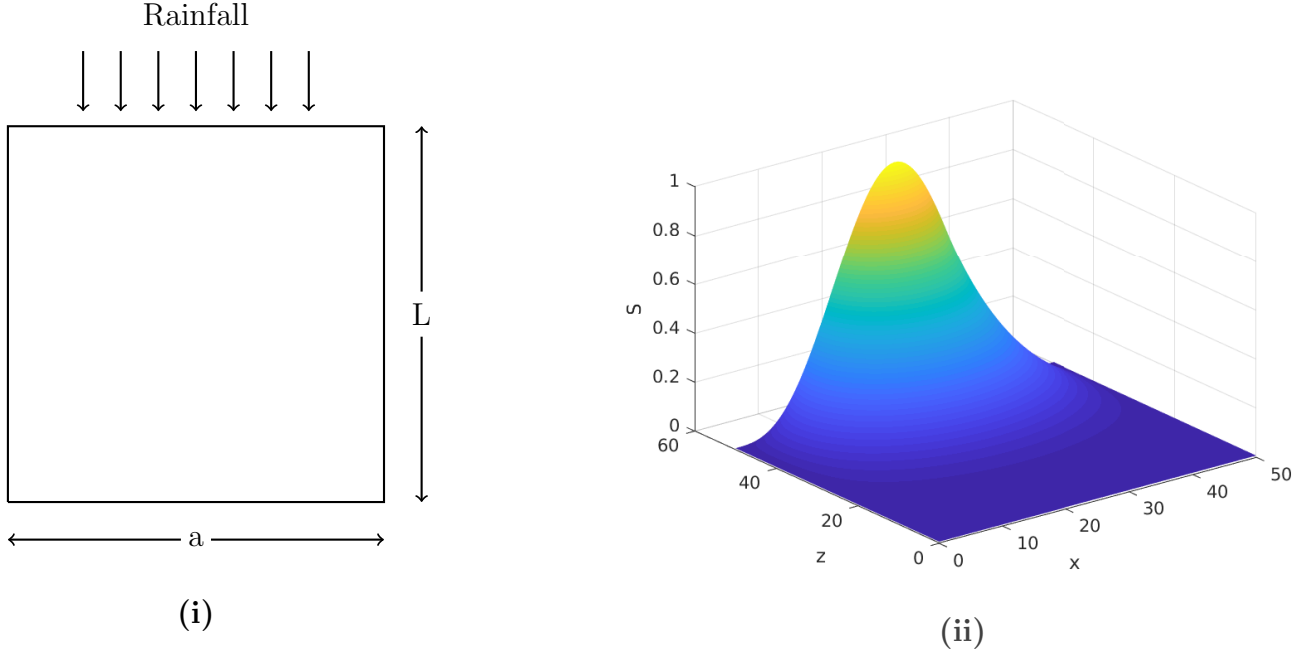}\\
\textbf{(ii)}
\end{minipage}
\caption{\textbf{(i)} Two-dimensional soil sample $[0,a]\times[0,L]$.  \textbf{(ii)} Analytical solution for effective saturation $S$ at time $t=10$ days.}
\label{Green-AmptGeometry}
\end{figure}
The column of soil is initially dry ($\psi=\psi_d$) and the infiltration is considered at the top boundary. This problem is referred to as the Green-Ampt problem in \cite{Tracy2011}. An analytical solution to this problem  using  Gardner formulation \cite{Gardner1958} for the relative hydraulic conductivity $K_r$ and a simple formulation from \cite{Warrick2003} for the water content $\theta$,
\begin{equation}
K_r=e^{(\alpha\psi)},\qquad \theta=\theta_r+(\theta_s-\theta_r)e^{(\alpha\psi)},
\label{GardnerModel}
\end{equation}
is derived in \cite{Tracy2011}, by imposing a specific pressure head at the top that is zero in the middle and tapers rapidly to $\psi_d$ at $x=0$ and $x=a$, and is maintained at $\psi=\psi_d$ along the bottom and vertical sides of the soil. The initial and boundary conditions are \cite{Tracy2011}:
\begin{equation}
\begin{aligned}
&\psi(x,z,0)=\psi_d,\quad (x,z)\in[0,a]\times[0,L],\\
&\psi(x,0,t)=\psi(a,z,t)=\psi(0,z,t)=\psi_d,\quad x\in[0,a], \text{ }z\in[0,L],  \text{ }t\geqslant0,\\
&\psi(x,L,t)=\dfrac{1}{\alpha}\log\Big(\varepsilon+\big(1-\varepsilon\big)\Big[\frac{3}{4}\sin\big(\frac{\pi x}{a}\big)-\frac{1}{4}\sin\big(\frac{3\pi x}{a}\big)\Big]\Big),\quad x\in[0,a], \text{ }z\in[0,L],  \text{ }t\geqslant0
\end{aligned}
\end{equation}
and the analytical solution of this problem is given by \cite{Tracy2011}
\begin{equation}
\psi(x,z,t)=\dfrac{1}{\alpha}\log(\varepsilon+\tilde{\psi}(x,z,t)),
\label{TracyExactSol1}
\end{equation}
where
\begin{equation}
\begin{aligned}
&\tilde{\psi}=(1-\varepsilon)\exp\big(\frac{\alpha}{2}(L-z)\big)\Big\lbrace
\frac{3}{4}\sin\big(\frac{\pi x}{a}\big)\Big[\dfrac{\sinh(\beta_1z)}{\sinh(\beta_1L)}+\frac{2}{Lb}\sum_{k=1}^{\infty}(-1)^k\frac{\lambda_k}{\gamma_{1k}}\sin(\lambda_kz)\exp(-\gamma_{1k}t)\Big] \\
&\qquad\qquad\qquad\qquad\qquad\qquad-\frac{1}{4}\sin\big(\frac{3\pi x}{a}\big)\Big[\dfrac{\sinh(\beta_3z)}{\sinh(\beta_3L)}+\frac{2}{Lb}\sum_{k=1}^{\infty}(-1)^k\frac{\lambda_k}{\gamma_{3k}}\sin(\lambda_kz)\exp(-\gamma_{3k}t)\Big]\Big\rbrace
\end{aligned}
\label{TracyExactSol2}
\end{equation}
with
\begin{equation}
\lambda_k=\dfrac{k\pi}{L},\text{ }
\varepsilon=e^{(\alpha\psi_d)},\text{ }
b=\dfrac{\alpha(\theta_s-\theta_r)}{K_s},\text{ }
\tilde{\lambda}_i=\dfrac{i\pi}{a},\text{ }
\beta_i=\sqrt{\frac{\alpha^2}{4}+\tilde{\lambda}_i^2},\text{ }
\gamma_{ik}=\dfrac{\beta_i^2+\lambda_k^2}{b}.
\end{equation}
Using \eqref{GardnerModel}, the pressure head $\psi$ can be expressed in the form \eqref{WaterSectionHead}, where
\begin{equation}
\quad h_{cap}(\bm{x})=\dfrac{1}{\alpha},\quad J(S)=\log(S).
\label{GardnerJfunction}
\end{equation}
In practice, one should truncate the series in \eqref{TracyExactSol2} to calculate the exact solution. Here, we consider the first $200$ terms of the series.

We solve the equation \eqref{RichardsEq} with the semi-implicit $(S,\psi)$-scheme in the domain represented in Figure~\ref{Green-AmptGeometry}\textbf{(i)} by using the following values for the parameters
\begin{equation}
L=a=50m,\text{ } K_s=0.2m/day,\text{ } \psi_d=-50m,\text{ } \theta_s=0.45,\text{ } \theta_r=0.15,\text{ } \alpha=0.1m^{-1}. 
\end{equation}
The analytical solution for the effective saturation  at time $t=10$ days is represented in Figure~\ref{Green-AmptGeometry}\textbf{(ii)}. The $L^2$- and $H^1$-error norms on the effective saturation $S_h$ and the pressure head $\psi_h$ are shown in Table~\ref{Tab:ErrorExactNumSol}. The contours of the analytical and numerical solutions for the effective saturation  at time $t=10$ days are plotted in Figure~\ref{Fig:ContoursExactNumSol}. The results show that the numerical solution converges to the analytical solution as the mesh is refined. An order of convergence varying around $k=1.8$ for the $L^2$-error and an order varying around $k=0.9$ for the $H^1$-error are obtained on both variables $S_h$ and $\psi_h$. However, the error on the pressure head is larger.
\begin{figure}[!h]
\centering
\begin{tabular}{ccc}
\includegraphics[scale=0.2]{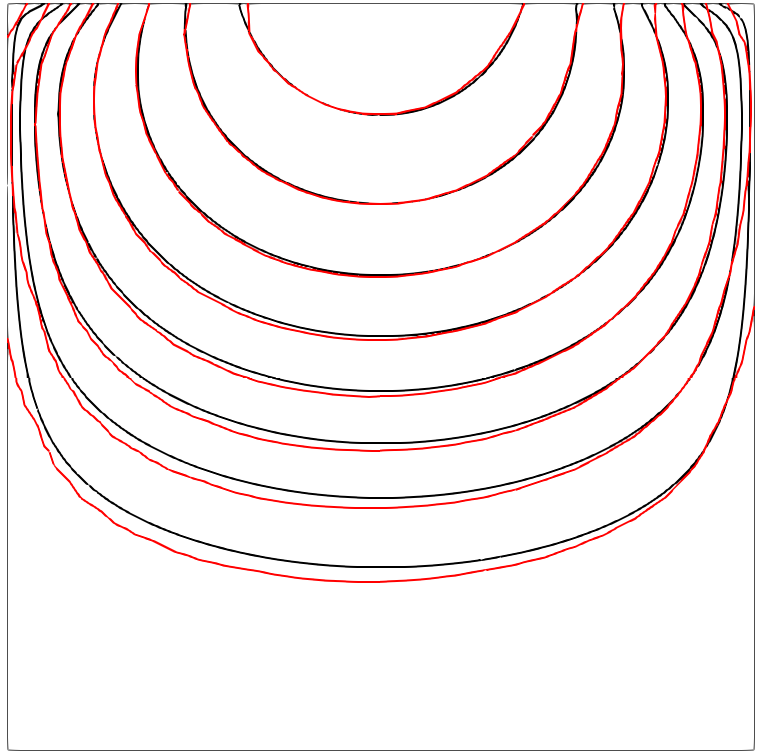}&\includegraphics[scale=0.2]{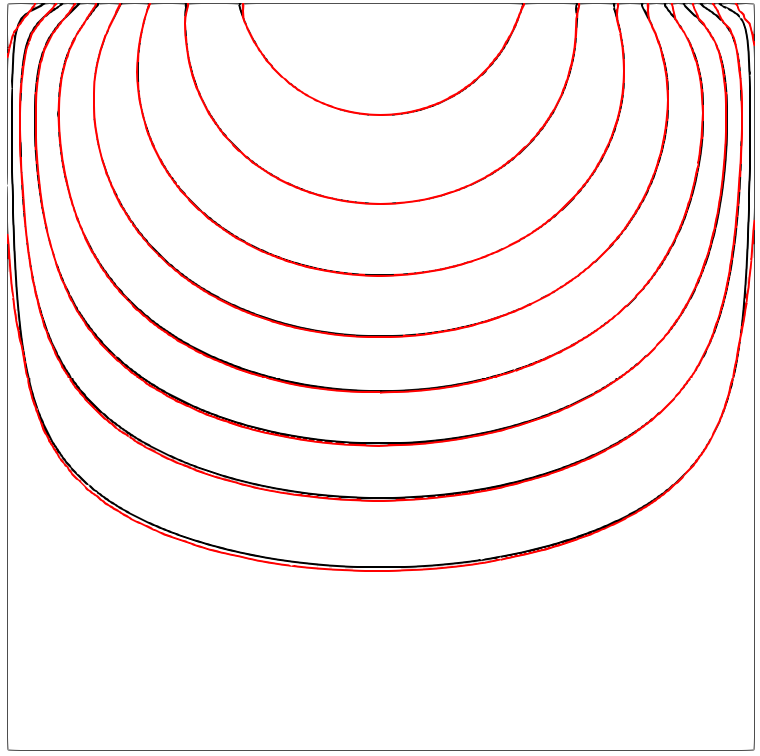}&\includegraphics[scale=0.2]{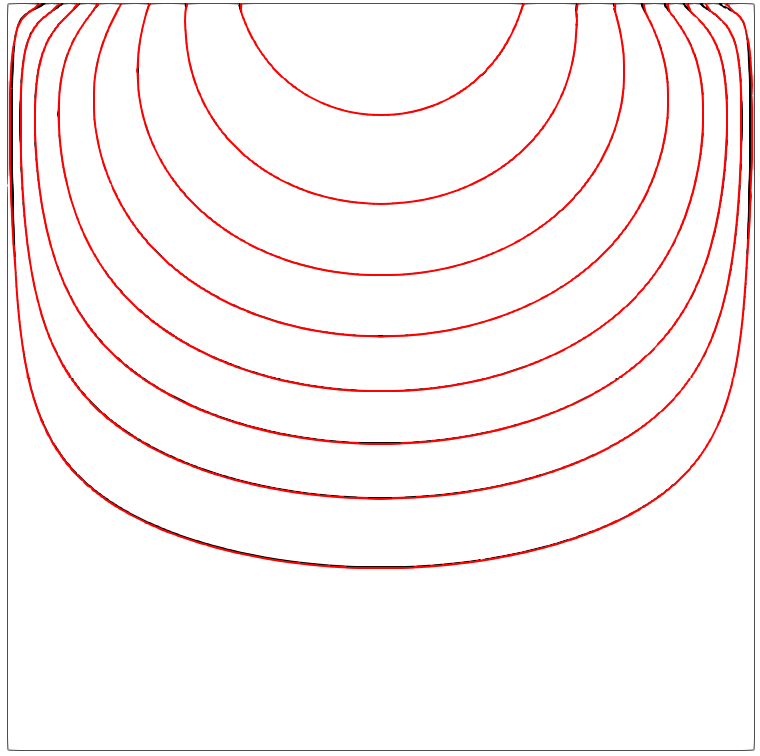}\\
Mesh $25\times 25$&Mesh $50\times 50$&Mesh $100\times 100$
\end{tabular}
\caption{Contours of the analytical (black curves) and numerical (red curves) solutions  on different meshes at time $t=10$ days.}
\label{Fig:ContoursExactNumSol}
\end{figure}
\begin{table}[!h]
\centering
\begin{tabular}{c|c|c|c|c|c|}
\cline{3-6}
\multicolumn{2}{c|}{}&\multicolumn{2}{|c|}{$L^2$-error } & \multicolumn{2}{|c|}{$H^1$-error}  \\ \cline{1-6}
\multicolumn{1}{|c|}{$\Delta t$}&Mesh&$S_h$ & $\psi_h$& $S_h$& $\psi_h$ \\ \cline{1-6}
\multicolumn{1}{|c|}{$0.010$}&$25\times 25$&$0.055429$&$ 26.3803$&$ 0.125187$&$41.3671$ \\ \cline{1-6}
\multicolumn{1}{|c|}{$0.005$}&$50\times 50$&$0.016745$&$8.72881$&$0.057976$&$22.2810$\\ \cline{1-6}
\multicolumn{1}{|c|}{$0.0025$}&$100\times 100$&$ 0.004397$&$2.45371$&$0.027922$&$11.9616$ \\ \cline{1-6}
\multicolumn{1}{|c|}{$0.00125$}&$200\times 200$&$ 0.001182$&$0.54719$&$0.013805$&$6.20522$ \\ \cline{1-6}
\end{tabular}
\caption{Analytical solution: $L^2$- and $H^1$-error on $S_h$ and $\psi_h$ for the solution computed with the semi-implicit $(S,\psi)$-scheme on different meshes.}
\label{Tab:ErrorExactNumSol}
\end{table}

To test the conservation property of the scheme, we compare the total mass in the computational domain of the exact and numerical solutions as time evolves. The time evolution of the total mass of the exact and numerical solutions for the saturation and the relative error are displayed in Figure~\ref{Fig:TotalMassEvolution}. The results show a good agreement for the total mass of the exact and numerical solutions.
\begin{figure}[!h]
\centering
\begin{tabular}{cc}
\includegraphics[scale=0.55]{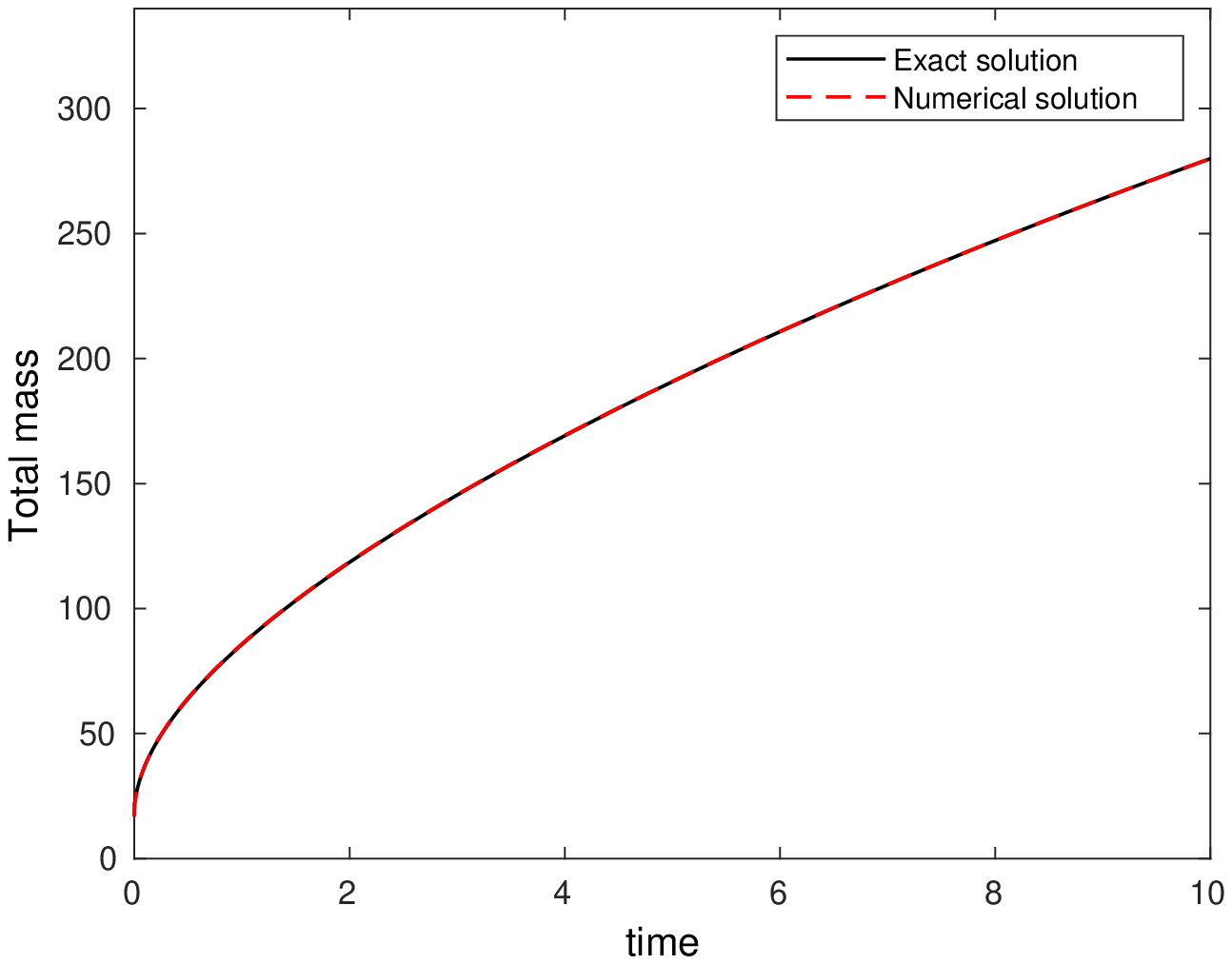}&\includegraphics[scale=0.55]{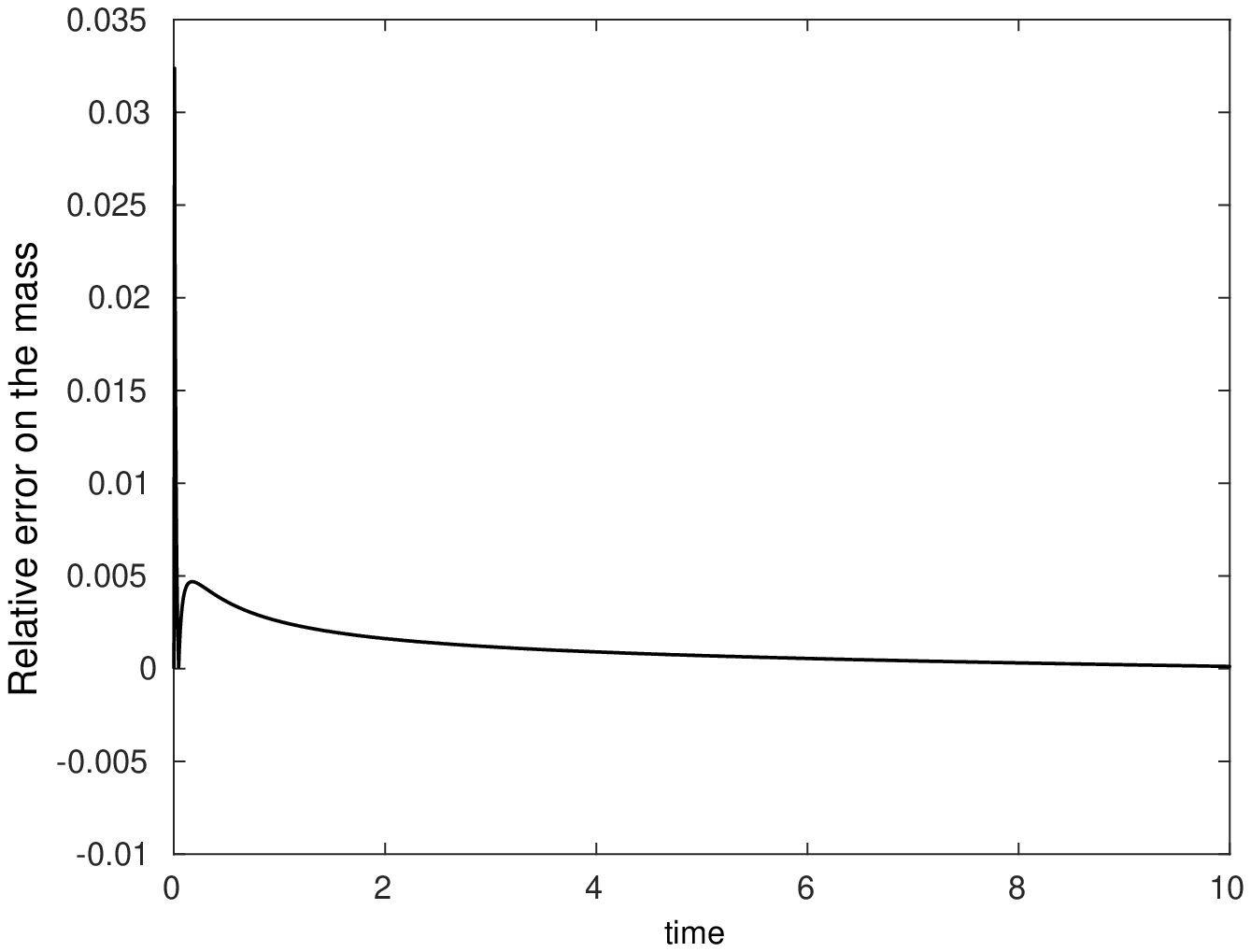}\\
\textbf{(a)}& \textbf{(b)}
\end{tabular}
\caption{ \textbf{(a)} Time-evolution of the total mass in the computational domain of the exact and numerical solutions and \textbf{(b)} relative error on the mesh $100\times 100$ with $\Delta t=0.0025$ day.}
\label{Fig:TotalMassEvolution}
\end{figure}

\subsection{Infiltration in heterogeneous medium}
\label{SectInfiltrationHeterogeneous}
The objective of the next test cases is to investigate the performance of the algorithms for infiltration in heterogeneous medium.
\subsubsection{Curvilinearly layered soil}
This test case was used by Manzini et al.\ \cite{Manzini2004}, where two layers of soil are considered. The capillary pressure is modeled by the van Genuchten constitutive relationship \cite{vanGenuchten1980}: 
\begin{equation}
\theta=\dfrac{\theta_s-\theta_r}{[1+(\alpha\vert\psi\vert)^n]^m}+\theta_r
\label{vanGenuchtenRelation1}
\end{equation}
and  the relative hydraulic conductivity is expressed as follows \cite{Mualem1976,vanGenuchten1980}:
\begin{equation}
K_r=\dfrac{[1-(\alpha\vert\psi\vert)^{n-1}(1+(\alpha\vert\psi\vert)^n)^{-m}]^2}{[1+(\alpha\vert\psi\vert)^n]^{m/2}},\quad m=1-1/n.
\label{vanGenuchtenRelation2}
\end{equation}
Using \eqref{vanGenuchtenRelation1}, the pressure head $\psi$ can be expressed in the form \eqref{WaterSectionHead}, where
\begin{equation}
\quad h_{cap}(\bm{x})=\dfrac{1}{\alpha},\quad J(S)=-(S^{-1/m}-1)^{1/n}
\label{vanGenuchtenJfunction}
\end{equation}
and the saturation can then be expressed as 
\begin{equation}
S=[1+(\alpha\vert\psi\vert)^n]^{-m}.
\end{equation}
The computational domain (see Figure~\ref{CurvilinearlyGeometry}) is the square $[0\,cm,100\,cm]\times[0\,cm,100\,cm]$. 
\begin{figure}[!h]
\centering
\begin{tikzpicture}
\draw[black, thick] (0,0) -- (5,0) -- (5,5) -- (0,5) -- (0,0);
\draw[scale=0.5,domain=0:10,smooth,variable=\x,black] plot ({\x},{10*(0.1*(1-cos(deg(pi*\x/10)))+0.45)});
\node[ ] at (2.5,3.5) {$z\geqslant\xi(x)$};
\node[ ] at (2.5,1.5) {$z\leqslant\xi(x)$};
\node[ ] at (0,-0.3) {0};
\node[ ] at (5,-0.3) {100};
\node[ ] at (-0.4,4.9) {100};
\end{tikzpicture}
\caption{Geometry of two layers of soil separated by a curved interface.}
\label{CurvilinearlyGeometry}
\end{figure}
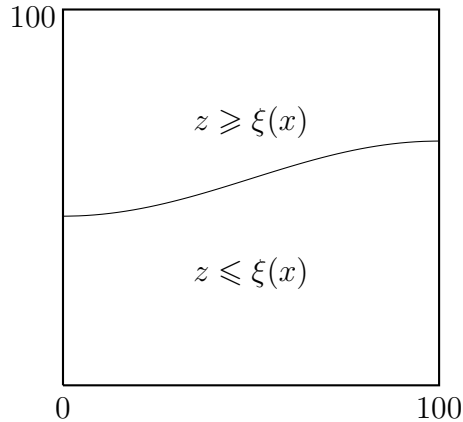
The two different regions of the layered soil are separated by the curved interface 
\begin{equation}
\xi(x)=100\big(0.1(1-\cos(\pi x/100))+0.45\big).
\end{equation}
The parameter values for \eqref{vanGenuchtenRelation1}-\eqref{vanGenuchtenRelation2} are given by 
\begin{equation}
\begin{aligned}
&z\geqslant\xi(x):\quad\theta_s=0.50,\quad \theta_r=0.120,\quad\alpha=0.028\,cm^{-1},\quad n=3.00,\quad K_s=0.25\,cm.h^{-1}, \\
&z\leqslant\xi(x):\quad\theta_s=0.46,\quad \theta_r=0.034,\quad\alpha=0.016\,cm^{-1},\quad n=1.37,\quad K_s=2.00\,cm.h^{-1}. 
\end{aligned}
\label{ManziniTestCaseParameters}
\end{equation}
Homogeneous Dirichlet boundary conditions are set on the top and bottom sides of the domain, while homogeneous Neumann conditions are set on the two vertical sides. We consider the initial condition  $\psi(x,z,t=0)=-z$, the mesh size $h=1\,cm$ and the time step $\Delta t=\frac{1}{60}$ hour. The Newton method applied to the implicit $(S,\psi)$-scheme does not converge with the parameter values \eqref{ManziniTestCaseParameters} because the relative permeability function $K_r$ expressed by \eqref{vanGenuchtenRelation2} is not differentiable at point $\psi=0$ for $n<2$. For physical values of $n$ with $0<n<2$, the implicit $(S,\psi)$-scheme may not converge. In this case, we use the semi-implicit $(S,\psi)$-scheme for the numerical simulations. The time-evolution of the water content is represented in Figure~\ref{Fig:Test4TimeEvolWaterContentMixed1SemiScheme}. Our numerical results agree very well with the ones presented by Manzini et al.\ \cite{Manzini2004}. This shows that linear semi-implicit $(S,\psi)$-scheme is a good alternative when implicit schemes encounter converging difficulties. 
\begin{figure}[!h]
\centering
\begin{tabular}{cccc}
\includegraphics[scale=0.15]{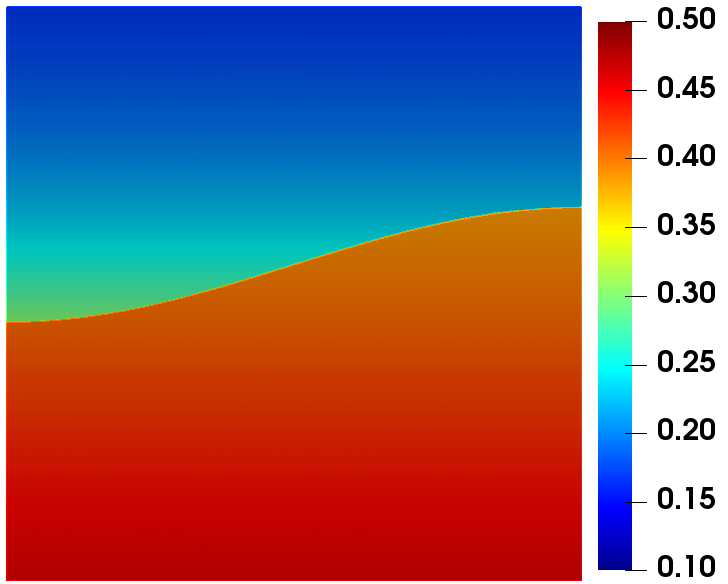}&\includegraphics[scale=0.15]{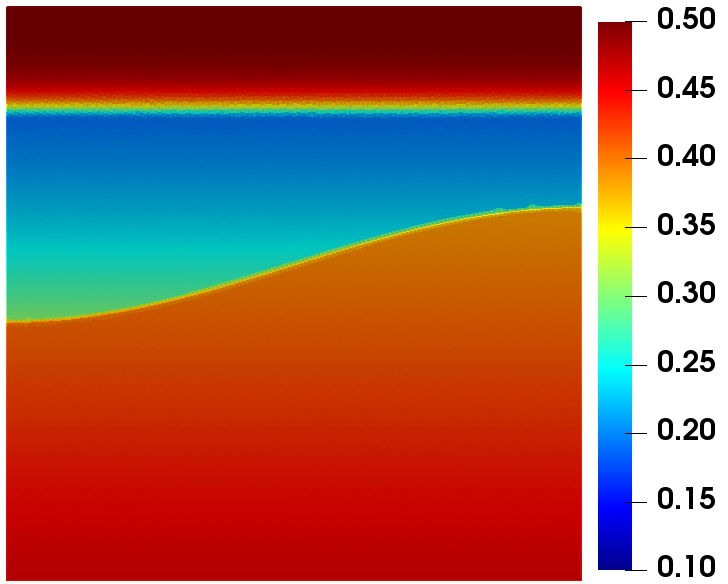}&\includegraphics[scale=0.15]{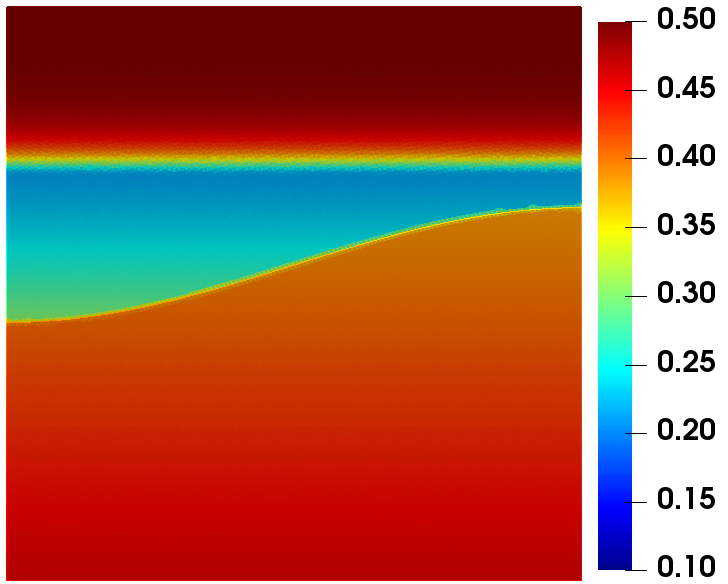}&\includegraphics[scale=0.15]{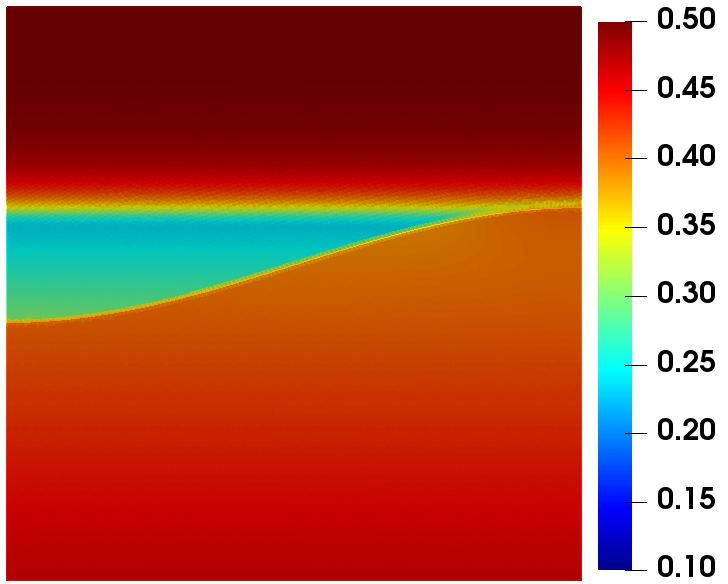}\\
$t=0$&$t=0.25$ day& $t=0.5$ day & $t=0.75$ day\\\\
\includegraphics[scale=0.15]{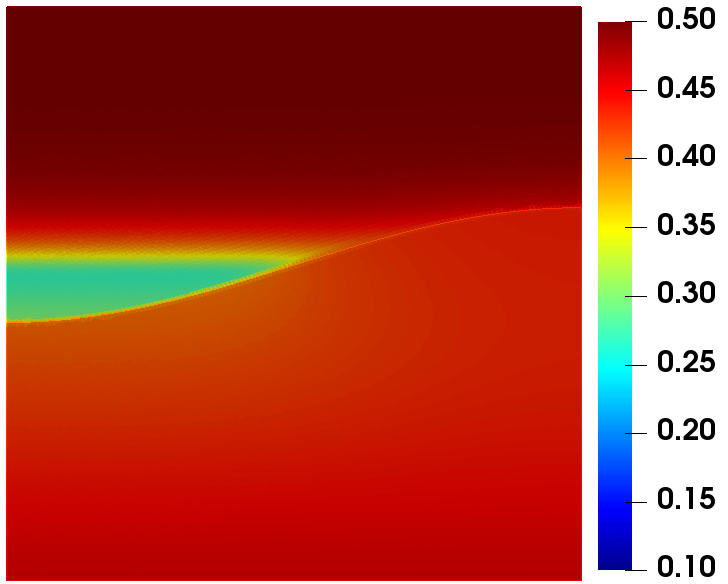}&\includegraphics[scale=0.15]{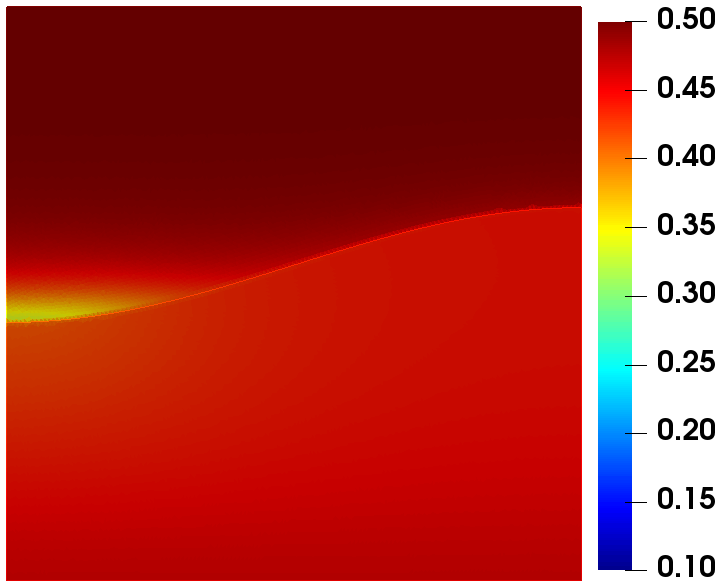}&\includegraphics[scale=0.15]{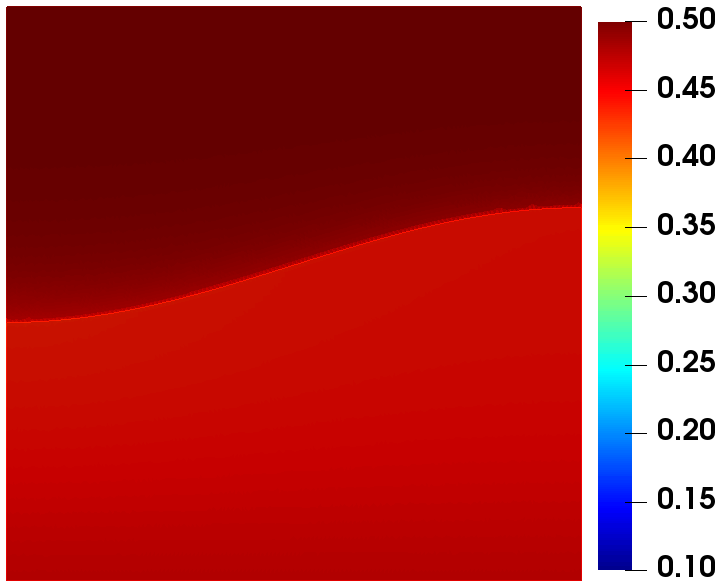}&\includegraphics[scale=0.15]{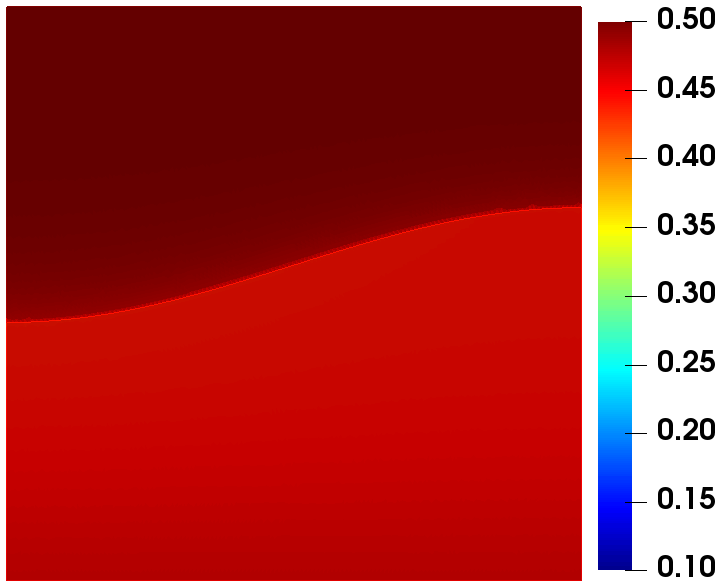}\\
$t=1$ day &$t=1.25$ day& $t=1.5$ day & $t=1.75$ day
\end{tabular}
\caption{Curvilinearly layered soil --- Test case 1: Time-evolution of the water content $\theta_h$ computed with the semi-implicit $(S,\psi)$-scheme.}
\label{Fig:Test4TimeEvolWaterContentMixed1SemiScheme}
\end{figure}

To test the performance of our scheme for large degree of heterogeneity, we reconduct the test case using large difference between the saturated hydraulic conductivities of the layers and we keep the same values for the other parameters. We perform computations for $K_s=2.5\,cm.h^{-1}$ and $K_s=25\,cm.h^{-1}$ for $z\leqslant \xi(x)$, respectively, using different time steps. Numerical results (not shown here) show that the semi-implicit $(S,\psi)$-scheme remains stable in modeling unsaturated flow through highly heterogeneous soils for a wide range of time steps.

To investigate the order of convergence of the semi-implicit $(S,\psi)$-scheme with this non regular solution, we compute a reference solution on a grid of mesh size $h=0.28\,cm$ ($500\times 500$ elements) with a time step $\Delta t=\frac{1}{240}$ hour at time $t=1$ day. The $L^2$- and $H^1$-error norms are computed with respect to this reference solution for different values of the mesh size to evaluate the spatial order of convergence. Figure~\ref{Fig:ErrorManziniMixed1SemiScheme} presents the errors as function of the mesh size. An order of convergence around $k=1.5$ for the $L^2$-error and around $k=1.2$ for the $H^1$-error are observed on the pressure head $\psi_h$ while an order of convergence around $k=0.6$ for the $L^2$-error and around $k=0.1$ for the $H^1$-error are observed on $S_h$. The lack of regularity of this solution deteriorates significantly the order of convergence on the saturation variable as observed in \cite{Traverso2013}. The errors on the pressure head $\psi_h$ are larger than the errors on the saturation $S_h$. However, the magnitude of the difference decreases with the mesh refinement since the orders  for $\psi_h$ are larger than the orders for $S_h$.
\begin{figure}[!h]
\centering
\begin{tabular}{cc}
\includegraphics[scale=0.55]{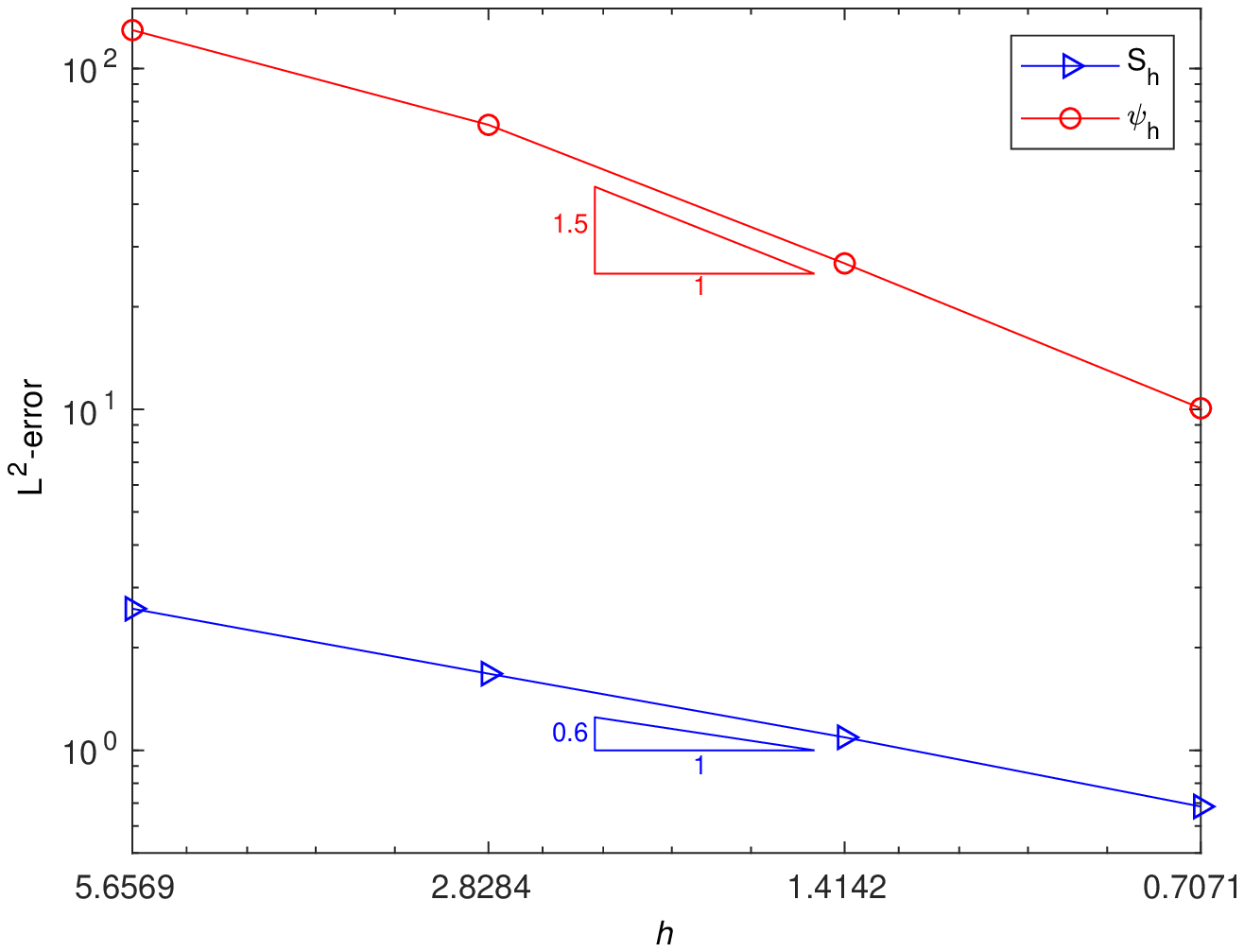}&\includegraphics[scale=0.55]{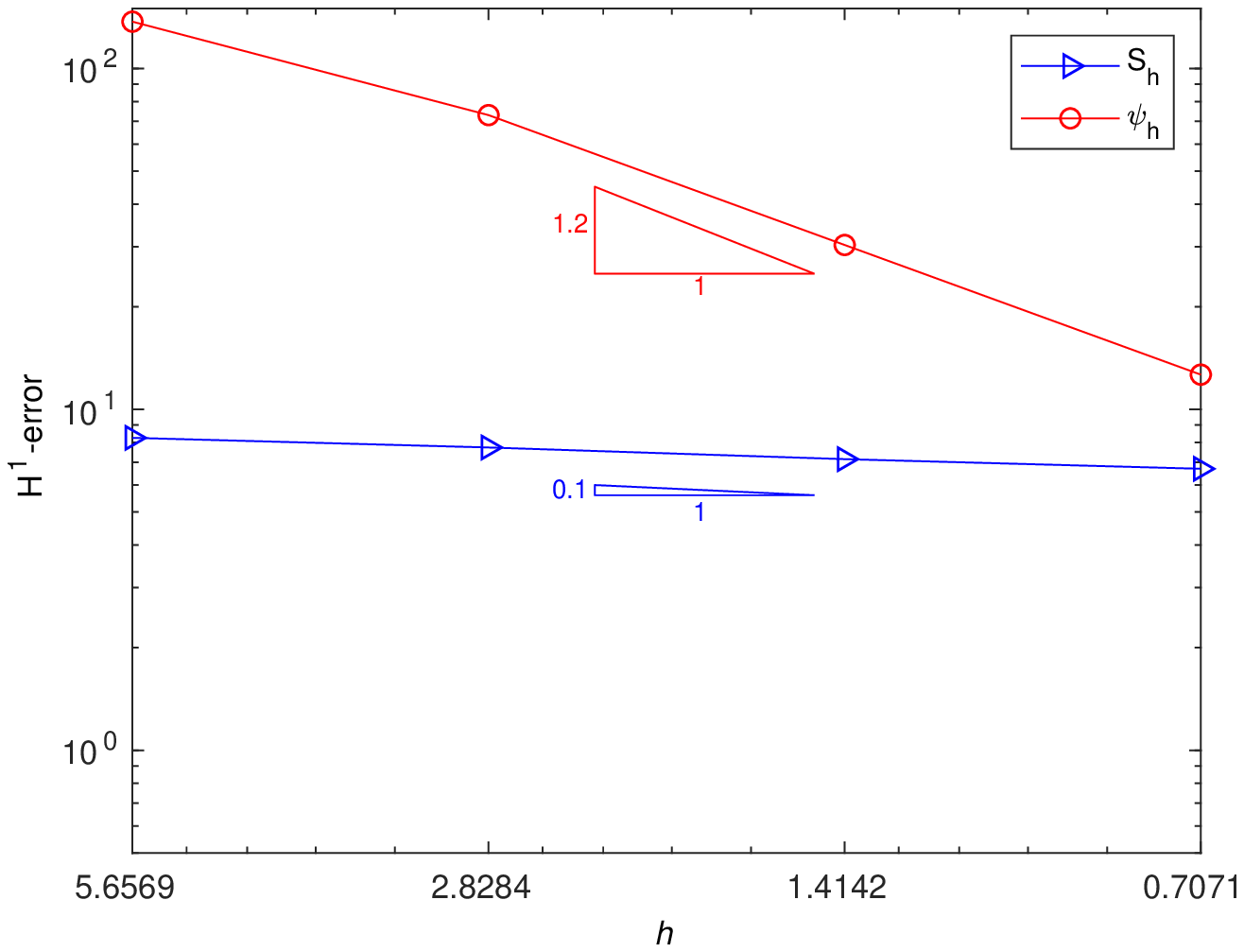}
\end{tabular}
\caption{Curvilinearly layered soil --- Test case 1: $L^2$- and $H^1$-error on $S_h$ and $\psi_h$ as a function of the mesh size $h$  computed with the semi-implicit $(S,\psi)$-scheme with respect to the reference solution.}
\label{Fig:ErrorManziniMixed1SemiScheme}
\end{figure}

Another test case with this geometry was used in \cite{BaronThesis2015} but for two layers of soil having the same characteristics except for the saturated hydraulic conductivity which is $K_s$ at top ($z\geqslant \xi(x)$) and $2K_s$  at bottom ($z\leqslant\xi(x)$). The parameter values for capillary pressure function \eqref{vanGenuchtenRelation1} and the relative permeability function \eqref{vanGenuchtenRelation2} are 
\begin{equation}
\theta_s=0.5,\quad \theta_r=0.12,\quad \alpha=0.02\,cm^{-1},\quad n=3,\quad K_s=6.944\times 10^{-5}cm.s^{-1}.
\label{BaronTestCaseParameters}
\end{equation}
The initial condition is $\psi(x,z,t=0)=-z$. A homogeneous Neumann condition is enforced on the lateral parts of the domain, while a homogeneous Dirichlet one is imposed on top and bottom sides. We reproduce this test case using the implicit $(S,\psi)$-scheme. The mesh size is $h=1\,cm$ and the time step $\Delta t=\frac{1}{36}$ hour. The time-evolution of the relative pressure $\psi_h^t-\psi_h^0$ is represented in Figure~\ref{Fig:OverPressureTestCase}. A good agreement of our numerical results with the ones presented in \cite{BaronThesis2015} is observed.
\begin{figure}[!h]
\begin{tabular}{ccc}
\includegraphics[scale=0.2]{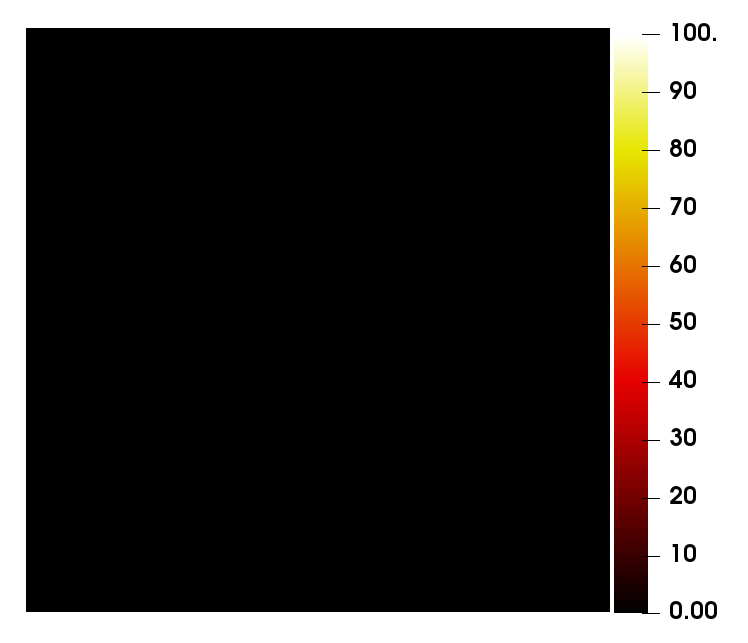}&\includegraphics[scale=0.2]{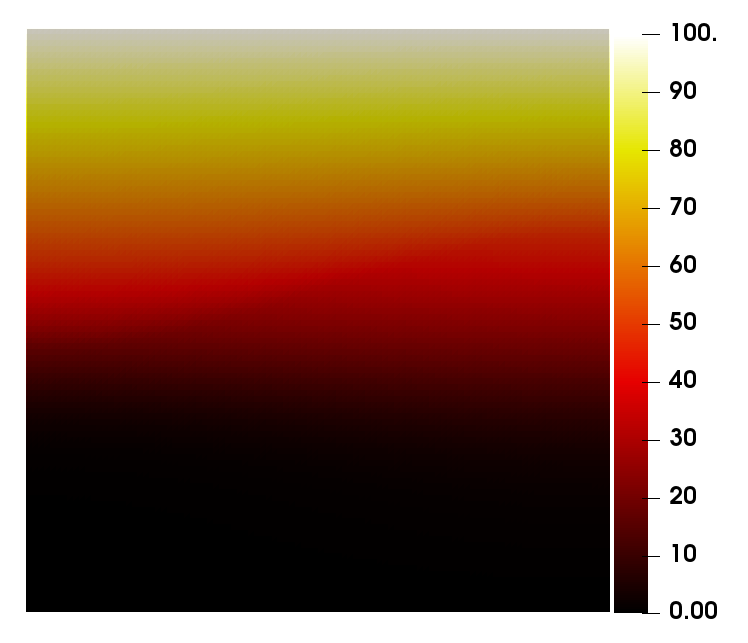}&\includegraphics[scale=0.2]{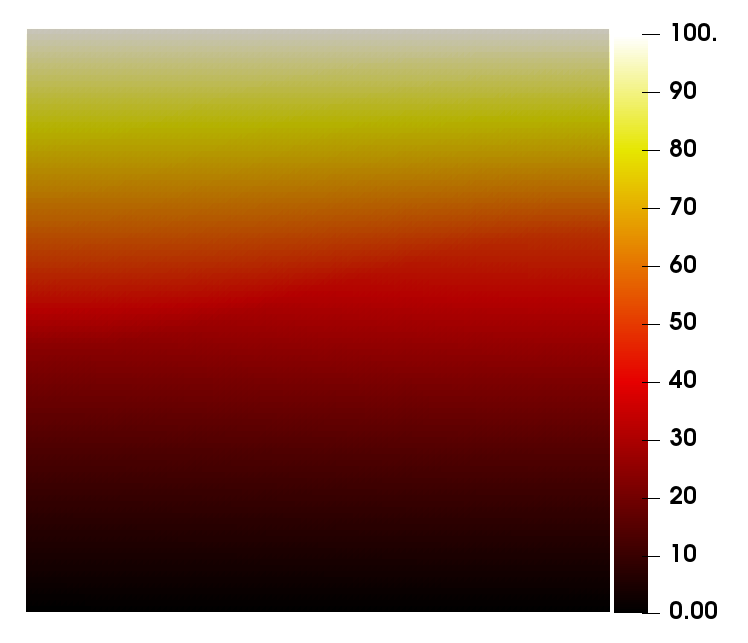}\\
$t=0$&$t=36h$&$t=72h$
\end{tabular}
\caption{Curvilinearly layered soil --- Test case 2: Relative pressure $\psi_h^t-\psi_h^0$ at different times $t$ computed with the implicit $(S,\psi)$-scheme.}
\label{Fig:OverPressureTestCase}
\end{figure}

\subsubsection{Layers of soil of L-shape form}
Here, we perform numerical simulations using the test case presented in {\cite{Baron2017}} for two layers of soil with same characteristics given by the relations {\eqref{vanGenuchtenRelation1}}-{\eqref{vanGenuchtenRelation2}} for the capillary pressure and the relative hydraulic conductivity with the parameter values {\eqref{BaronTestCaseParameters}}, except for the saturated hydraulic conductivity $K_s$ which changes in the medium. The computational domain is the square $[0\,cm,100\,cm]\times[0\,cm,100\,cm]$ with variable saturated hydraulic conductivity as shown in Figure~\ref{L-ShapeGeometry}. The initial condition is $\psi(x,z,t=0)=-z$. A homogeneous Neumann condition is enforced on the lateral parts of the domain, while a homogeneous Dirichlet one is imposed on top and bottom sides. We discretize the domain using the mesh size $h=1\,cm$. The time step is $\Delta t=\frac{1}{180}$ hour. The time-evolution of the effective saturation computed with the implicit $(S,\psi)$-scheme is represented in Figure~\ref{Fig:L-ShapeTestCase}. Good agreement with the results presented in \cite{Baron2017} is observed.
\begin{figure}[!h]
\centering
\begin{tikzpicture}
\draw[black, fill=gray, thick] (0,0) -- (5,0) -- (5,4) -- (2.5,4) -- (2.5,1) -- (0,1) -- (0,0) ;
\draw[black, fill=gray!50, thick] (0,1) -- (2.5,1) -- (2.5,4) -- (2.5,4) -- (5,4) -- (5,5) -- (0,5) -- (0,1) ;
\node[ ] at (1.25,2.5) {$K_s$};
\node[ ] at (3.75,2.5) {$3K_s$};
\node[ ] at (0,-0.3) {0};
\node[ ] at (5,-0.3) {100};
\node[ ] at (-0.4,4.9) {100};
\end{tikzpicture}
\caption{Geometry of two layers of soil of L-shape form.}
\label{L-ShapeGeometry}
\end{figure}
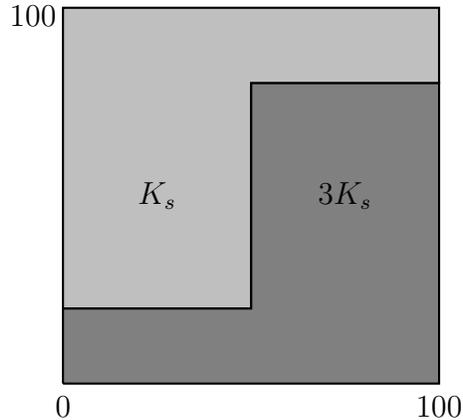
\begin{figure}[!h]
\centering
\begin{tabular}{ccc}
\includegraphics[scale=0.2]{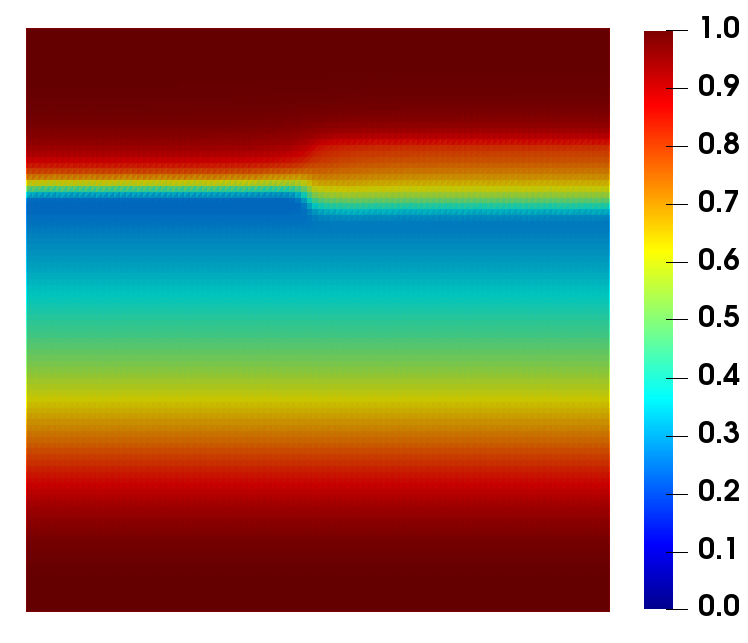}&\includegraphics[scale=0.2]{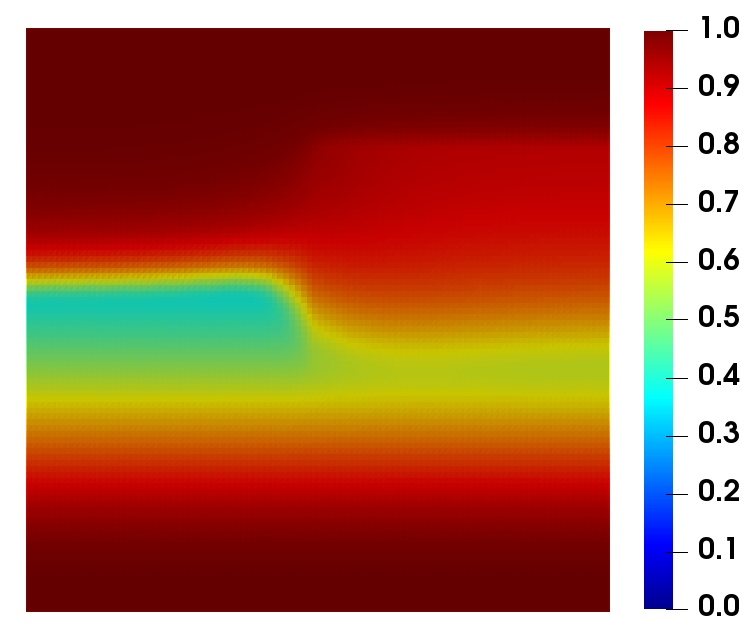}&\includegraphics[scale=0.2]{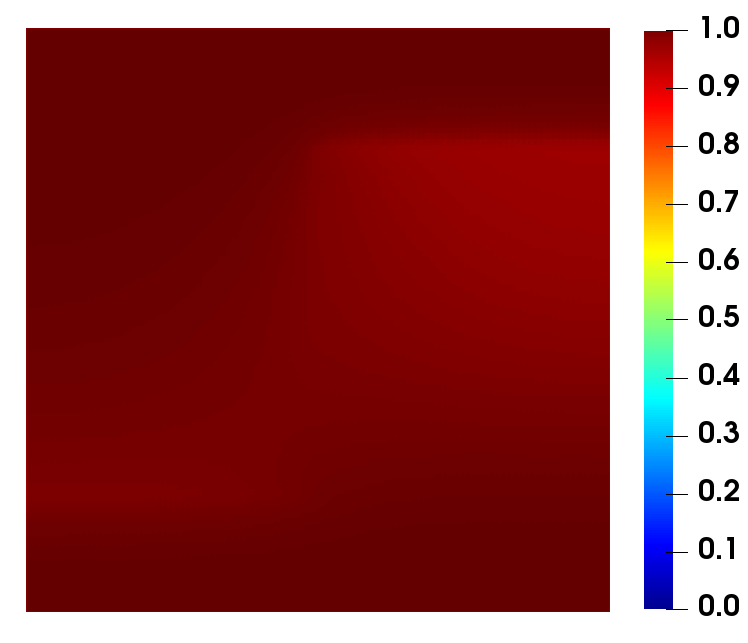}\\
$t=12h$&$t=24h$&$t=48h$
\end{tabular}
\caption{L-shape geometry test case: Saturation $S_h$ at different times computed with the implicit $(S,\psi)$-scheme.}
\label{Fig:L-ShapeTestCase}
\end{figure}

To evaluate the order of convergence of the implicit $(S,\psi)$-scheme with this solution, we compute a reference solution on a grid of mesh size $h=0.28\,cm$ ($500\times 500$ elements) at time $t=24h$. Figure~\ref{Fig:ErrorBaronMixed1ImplScheme} shows the errors as function of the mesh size. A convergence rate around $k=1.8$ for the $L^2$-error and around $k=1.3$ for the $H^1$-error  are observed on each variable. The errors on the pressure head $\psi_h$ are about $70$ times larger than the errors on the saturation $S_h$.
\begin{figure}[!h]
\centering
\begin{tabular}{cc}
\includegraphics[scale=0.55]{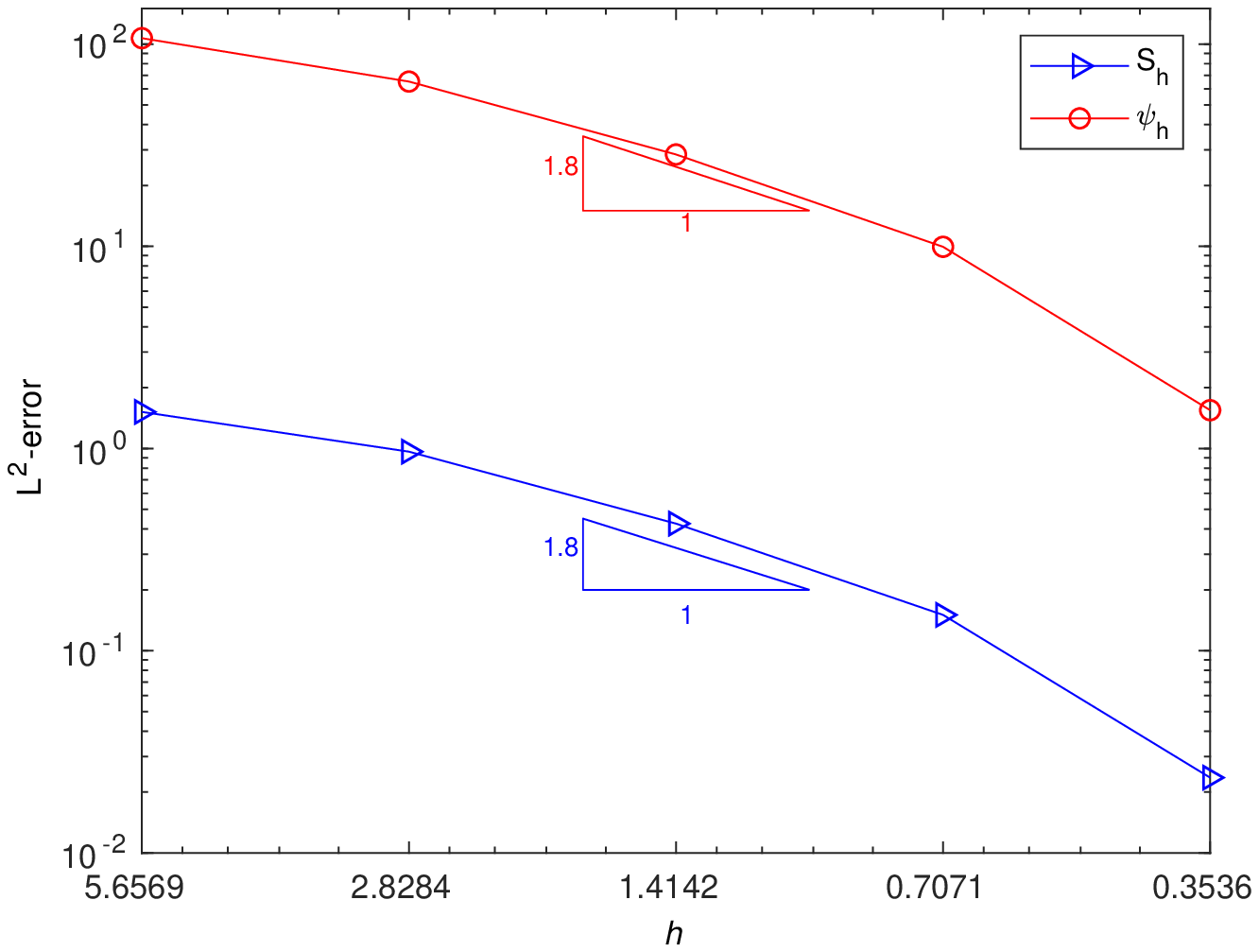}&\includegraphics[scale=0.55]{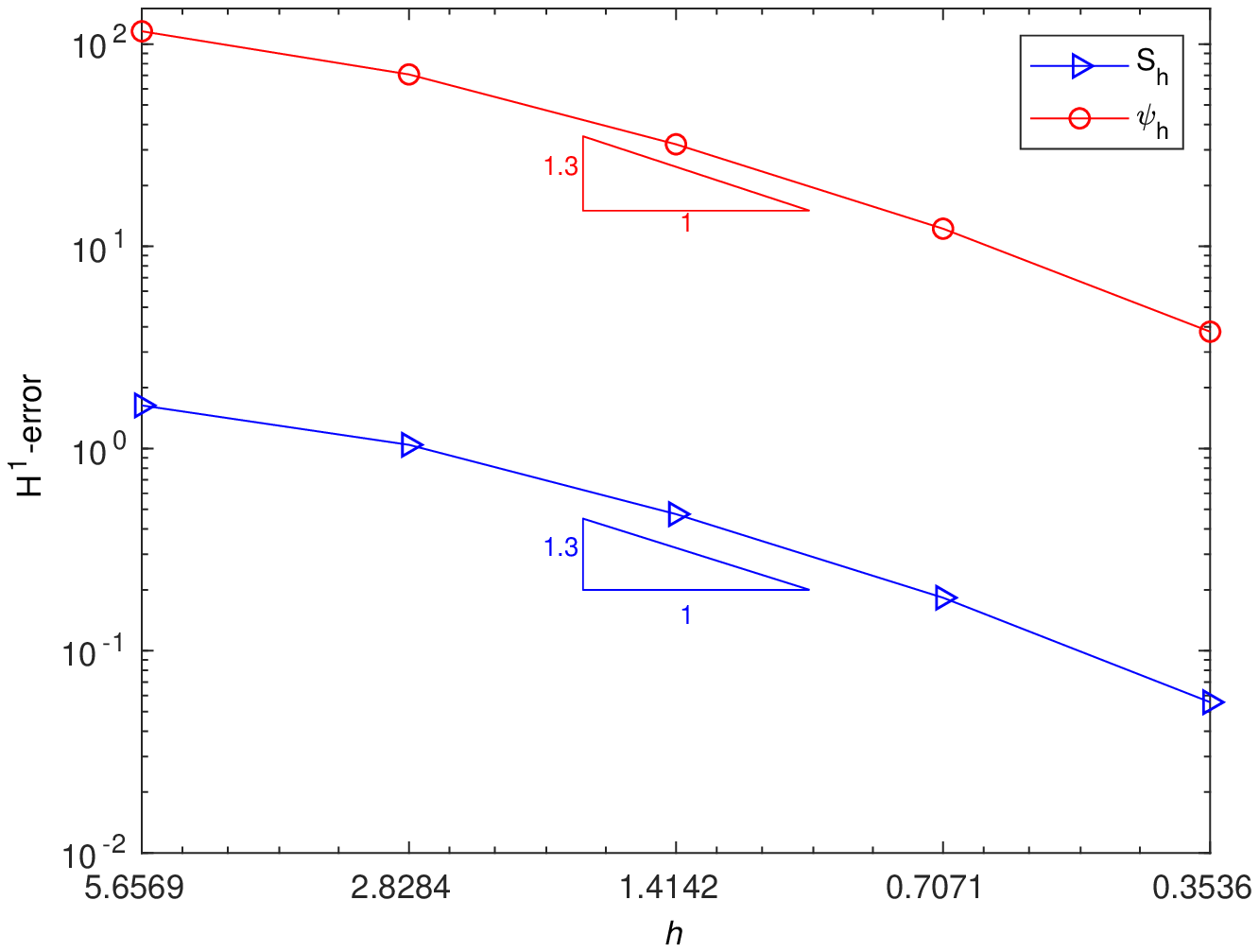}
\end{tabular}
\caption{L-shape geometry test case: $L^2$- and $H^1$-error on $S_h$ and $\psi_h$ as a function of the mesh size $h$ computed with the implicit $(S,\psi)$-scheme with respect to the reference solution.}
\label{Fig:ErrorBaronMixed1ImplScheme}
\end{figure}

\begin{remark}
All the test cases in section~\ref{SectInfiltrationHeterogeneous} have two saturated boundary conditions on both top and bottom sides. The environments tend to full saturation. Therefore, they serve also as an additional check of convergence behavior of the $(S,\psi)$-schemes with the regularization \eqref{RegularizedJfunction} when the saturation reaches the value of one.
\end{remark}

\section{Conclusion}
\label{SectConclusion}
The present study focuses on several iterative and noniterative second-order time stepping standard/or mixed finite element methods to approximate the mixed form of the Richards equation written in three different weak formulations. The first weak formulation uses the saturation as the sole unknown and is discretized using the standard linear Lagrange elements. The second weak formulation uses the saturation and the vector flux as unknowns and is discretized using the linear Lagrange elements for the saturation and the lowest order Raviart-Thomas elements for the vector flux. The third weak formulation is based on the saturation and the pressure head and is discretized using the linear Lagrange element for both variables. The different methods are studied and compared with respect to accuracy and computational requirements using exact solutions and commonly used test cases. Some numerical results are presented involving solutions that require regularization techniques to avoid degeneracy of the capillary pressure function when the saturation reaches the value of one in some regions of the domain.

Our numerical investigations show that the $(S,\psi)$-schemes (derived from the weak formulation using the saturation and the pressure head as unknowns) perform better in terms of accuracy and efficiency than the $(S,\bm{q})$-schemes (derived from the weak formulation using the saturation and the flux head as unknowns) and the $S$-schemes (derived from the weak formulation using the saturation as the sole unknown). The results of our study also indicate that second order accurate time stepping methods are more efficient than the traditional first order accurate schemes.

The implicit $(S,\psi)$-scheme allows the use of larger time steps. However, the iterative procedure may not converge because of the nonlinearity of the soil constitutive laws. The semi-implicit $(S,\psi)$-scheme is linear, as accurate as the implicit $(S,\psi)$-scheme, easy to implement and requires less memory. In summary, our results show that the proposed semi-implicit $(S,\psi)$-scheme is the best among the methods studied and is a good alternative to traditional iterative methods for the Richards equation.

\section*{Acknowledgments}
The first author was supported through a Postdoctoral Fellowship of the UM6P/OCP group of Morocco and a Postdoctoral Fellowship of the Fields Institute. The second author acknowledges funding from  UM6P/OCP group of Morocco. The research of the third author is funded through a Discovery Grant of the Natural Sciences and Engineering Research Council of Canada.

\bibliographystyle{plain} 
\bibliography{biblio}

\begin{thebibliography}{10}

\bibitem{Akrivis2015}
G.~Akrivis.
\newblock Stability of implicit-explicit backward difference formulas for
  nonlinear parabolic equations.
\newblock {\em SIAM J. Numer. Anal.}, 53:464--484, 2015.

\bibitem{Allaire2007}
G.~Allaire.
\newblock {\em Numerical {A}nalysis and {O}ptimization}.
\newblock Oxford University Press, New York, 2007.

\bibitem{Arbogast1993}
T.~Arbogast, M.~Obeyesekere, and M.~F. Wheeler.
\newblock Numerical methods for the simulation of flow in root-soil systems.
\newblock {\em SIAM J. Numer. Anal.}, 30:1677--1702, 1993.

\bibitem{Arbogast1996}
T.~Arbogast and M.~F. Wheeler.
\newblock A nonlinear mixed finite element method for a degenerate parabolic
  equation arising in flow in porous media.
\newblock {\em SIAM J. Numer. Anal.}, 33:1669--1687, 1996.

\bibitem{Ascher1995}
U.~M. Ascher, S.~J. Ruuth, and B.~T.~R. Wetton.
\newblock Implicit-explicit methods for time-dependent partial differential
  equations.
\newblock {\em SIAM J. Numer. Anal.}, 32:797--823, 1995.

\bibitem{BaronThesis2015}
V.~Baron.
\newblock {\em M\'ethodes num\'eriques pour les \'ecoulements en milieu poreux
  : estimations \`a posteriori et strat\'egie d'adaptation.}
\newblock PhD thesis, LMJL-Laboratoire de Math\'ematiques Jean Leray, 2015.

\bibitem{Baron2017}
V.~Baron, Y.~Coudière, and P.~Sochala.
\newblock Adaptive multistep time discretization and linearization based on a
  posteriori error estimates for the {R}ichards equation.
\newblock {\em Appl. Numer. Math.}, 112:104--125, 2017.

\bibitem{Bause2008}
M.~Bause.
\newblock Higher and lowest order mixed finite element approximation of
  subsurface flow problems with solutions of low regularity.
\newblock {\em Adv. Water Resour.}, 31:370--382, 2008.

\bibitem{Bause2004}
M.~Bause and P.~Knabner.
\newblock Computation of variably saturated subsurface flow by adaptive mixed
  hybrid finite element methods.
\newblock {\em Adv. Water Resour.}, 27:565--581, 2004.

\bibitem{Bear1972}
J.~Bear.
\newblock {\em Dynamics of {F}luids in {P}orous {M}edia}.
\newblock American Elsevier Publishing Company, New York, 1972.

\bibitem{Becker1998}
J.~Becker.
\newblock A second order backward difference method with variable steps for a
  parabolic problem.
\newblock {\em BIT Numer. Math.}, 38:644--662, 1998.

\bibitem{Belfort2009}
B.~Belfort, F.~Ramasomanana, A.~Younes, and F.~Lehmann.
\newblock An efficient lumped mixed hybrid finite element formulation for
  variably saturated groundwater flow.
\newblock {\em Vadose Zone J.}, 8:352--362, 2009.

\bibitem{Bergamaschi1999}
L.~Bergamaschi and M.~Putti.
\newblock Mixed finite elements and {N}ewton-type linearizations for the
  solution of {R}ichards' equation.
\newblock {\em Int. J. Numer. Meth. Engng.}, 45:1025--1046, 1999.

\bibitem{Berninger2011}
H.~Berninger, R.~Kornhuber, and O.~Sander.
\newblock Fast and robust numerical solution of the {R}ichards equation in
  homogeneous soil.
\newblock {\em SIAM J. Numer. Anal.}, 49:2576--2597, 2011.

\bibitem{Berninger2015}
H.~Berninger, R.~Kornhuber, and O.~Sander.
\newblock A multidomain discretization of the {R}ichards equation in layered
  soil.
\newblock {\em Comput. Geosci.}, 19:213--232, 2015.

\bibitem{Fortin1991}
F.~Brezzi and M.~Fortin.
\newblock {\em Mixed and {H}ybrid {F}inite {E}lement {M}ethod}.
\newblock Springer, New York, 1991.

\bibitem{Broadbridge1988}
P.~Broadbridge and I.~White.
\newblock Constant rate rainfall infiltration: {A} versatile nonlinear model:
  I. {A}nalytic solution.
\newblock {\em Water Resour. Res.}, 24:145--154, 1988.

\bibitem{Brooks1966}
R.~H. Brooks and A.~T. Corey.
\newblock Properties of porous media affecting fluid flow.
\newblock {\em J. Irrig. Drain. Div. Am. Soc. Civ. Eng.}, 92:61--88, 1966.

\bibitem{Cai2020}
W.~Cai, J.~Wang, and K.~Wang.
\newblock Convergence analysis of {C}rank-{N}icolson {G}alerkin-{G}alerkin
  {FEM}s for miscible displacement in porous media.
\newblock {\em J. Sci. Comput.}, 83:25, 2020.

\bibitem{Casulli2010}
V.~Casulli and P.~Zanolli.
\newblock A nested {N}ewton-type algorithm for finite volume methods solving
  {R}ichards' equation in mixed form.
\newblock {\em SIAM J. Sci. Comput.}, 32:2255--2273, 2010.

\bibitem{Celia1990}
M.~A. Celia, E.~T. Bouloutas, and R.~L. Zarba.
\newblock A general mass-conservative numerical solution for the unsaturated
  flow equation.
\newblock {\em Water Resour. Res.}, 26:1483--1496, 1990.

\bibitem{Ciarlet1989}
P.~G. Ciarlet.
\newblock {\em Introduction to {N}umerical {L}inear {A}lgebra and
  {O}ptimisation}.
\newblock Cambridge University Press, Cambridge, 1989.

\bibitem{CLement1994}
T.~P. Clement, W.~R. Wise, and F.~J. Molz.
\newblock A physically based, two-dimensional, finite-difference algorithm for
  modeling variably saturated flow.
\newblock {\em J. Hydrol.}, 161:71--90, 1994.

\bibitem{Cooley1983}
R.~L. Cooley.
\newblock Some new procedures for numerical solution of variably saturated flow
  problems.
\newblock {\em Water Resour. Res.}, 19:1271--1285, 1983.

\bibitem{Crouzeix1984}
M.~Crouzeix and F.~J. Lisbona.
\newblock The convergence of variable-stepsize, variable-formula, multistep
  methods.
\newblock {\em SIAM J. Numer. Anal.}, 21:512--534, 1984.

\bibitem{Cumming2011}
B.~Cumming, T.~Moroney, and I.~Turner.
\newblock A mass-conservative control volume-finite element method for solving
  {R}ichards' equation in heterogeneous porous media.
\newblock {\em BIT Numer. Math.}, 51:845--864, 2011.

\bibitem{Emmrich2009}
E.~Emmrich.
\newblock Convergence of the variable two-step {BDF} time discretisation of
  nonlinear evolution problems governed by a monotone potential operator.
\newblock {\em BIT Numer. Math.}, 49:297--323, 2009.

\bibitem{Ethier2008}
M.~Ethier and Y.~Bourgault.
\newblock Semi-implicit time-discretization schemes for the bidomain model.
\newblock {\em SIAM J. Numer. Anal.}, 46:2443--2468, 2008.

\bibitem{Eymard1999}
R.~Eymard, M.~Gutnic, and D.~Hilhorst.
\newblock The finite volume method for {R}ichards equation.
\newblock {\em Comput. Geosci.}, 3:259--294, 1999.

\bibitem{Farthing2003}
M.~W. Farthing, C.~E. Kees, and C.~T. Miller.
\newblock Mixed finite element methods and higher order temporal approximations
  for variably saturated groundwater flow.
\newblock {\em Adv. Water Resour.}, 26:373--394, 2003.

\bibitem{Farthing2000}
M.~W. Farthing and C.~T. Miller.
\newblock A comparison of high-resolution, finite-volume, adaptive–stencil
  schemes for simulating advective-dispersive transport.
\newblock {\em Adv. Water Resour.}, 24:29--48, 2000.

\bibitem{Forsyth1997}
P.~A. Forsyth and M.~C. Kropinski.
\newblock Monotonicity considerations for saturated-unsaturated subsurface
  flow.
\newblock {\em SIAM J. Sci. Comput.}, 18:1328--1354, 1997.

\bibitem{Forsyth1995}
P.~A. Forsyth, Y.~S. Wu, and K.~Pruess.
\newblock Robust numerical methods for saturated-unsaturated flow with dry
  initial conditions in heterogeneous media.
\newblock {\em Adv. Water Resour.}, 18:25--38, 1995.

\bibitem{Fujita1952}
H.~Fujita.
\newblock The exact pattern of a concentration-dependent diffusion in a
  semi-infinite medium, {P}art ii.
\newblock {\em Text. Res. J.}, 22:823--827, 1952.

\bibitem{Gardner1958}
W.~R. Gardner.
\newblock Some steady-state solutions of the unsaturated moisture flow equation
  with application to evaporation from a water table.
\newblock {\em Soil Sci.}, 85:228--232, 1958.

\bibitem{Haverkamp1977}
R.~Haverkamp, M.~Vauclin, J.~Touma, P.~J. Wierenga, and G.~Vachaud.
\newblock A comparison of numerical simulation models for one-dimensional
  infiltration.
\newblock {\em Soil Sci. Soc. Am. J.}, 41:285--294, 1977.

\bibitem{FreeFem}
F.~Hecht.
\newblock New development in {F}ree{F}em++.
\newblock {\em J. Numer. Math.}, 20:251--265, 2012.

\bibitem{FreeFemm}
F.~Hecht.
\newblock Free{F}em++: {V}ersion 4.1.
\newblock \url{https://doc.freefem.org/documentation/index.html}, 2019.

\bibitem{Karthikeyan2011}
M.~Karthikeyan, T.-S. Tan, and K.-K. Phoon.
\newblock Numerical oscillation in seepage analysis of unsaturated soils.
\newblock {\em Can. Geotech. J.}, 38:639--651, 2001.

\bibitem{Kavetski2002}
D.~Kavetski, P.~Binning, and S.~W. Sloan.
\newblock Noniterative time stepping schemes with adaptive truncation error
  control for the solution of {R}ichards equation.
\newblock {\em Water Resour. Res.}, 38:29--1--29--10, 2002.

\bibitem{Keita2020}
S.~Keita, A.~Beljadid, and Y.~Bourgault.
\newblock Mass-conservative and positivity preserving second-order
  semi-implicit methods for high-order parabolic equations.
\newblock {\em J. Comput. Phys.}, Revision submitted, arXiv:2010.11913
  [math.NA], 2020.

\bibitem{Keita2021}
S.~Keita, A.~Beljadid, and Y.~Bourgault.
\newblock Efficient second-order semi-implicit finite element method for
  fourth-order nonlinear diffusion equations.
\newblock {\em Comput. Phys. Commun.}, 258:107588, 2021.

\bibitem{Lehmann1998}
F.~Lehmann and P.~Ackerer.
\newblock Comparison of iterative methods for improved solutions of the fluid
  flow equation in partially saturated porous media.
\newblock {\em Transp. Porous Media}, 31:275--292, 1998.

\bibitem{LeVeque2007}
R.~J. LeVeque.
\newblock {\em Finite {D}ifference {M}ethods for {O}rdinary and {P}artial
  {D}ifferential {E}quations}.
\newblock Society for Industrial and Applied Mathematics, 2007.

\bibitem{Leverett1941}
M.~C. Leverett.
\newblock Capillary behavior of porous solids.
\newblock {\em Transaction of AIME}, 142:152--169, 1941.

\bibitem{List2016}
F.~List and F.~A. Radu.
\newblock A study on iterative methods for solving {R}ichards' equation.
\newblock {\em Comput. Geosci.}, 20:341--353, 2016.

\bibitem{Lott2012}
P.~A. Lott, H.~F. Walker, C.~S. Woodward, and U.~M. Yang.
\newblock An accelerated {P}icard method for nonlinear systems related to
  variably saturated flow.
\newblock {\em Adv. Water Resour.}, 38:92--101, 2012.

\bibitem{Manzini2004}
G.~Manzini and S.~Ferraris.
\newblock Mass-conservative finite volume methods on {2-D} unstructured grids
  for the {R}ichards' equation.
\newblock {\em Adv. Water Resour.}, 27:1199--1215, 2004.

\bibitem{McBride2006}
D.~McBride, M.~Cross, N.~Croft, C.~Bennett, and J.~Gebhardt.
\newblock Computational modelling of variably saturated flow in porous media
  with complex three-dimensional geometries.
\newblock {\em Int. J. Numer. Meth. Fluids}, 50:1085--1117, 2006.

\bibitem{Miller1998}
C.~T. Miller, G.~A. Williams, C.~T. Kelley, and M.~D. Tocci.
\newblock Robust solution of {R}ichards' equation for nonuniform porous media.
\newblock {\em Water Resour. Res.}, 34:2599--2610, 1998.

\bibitem{Milly1985}
P.~C.~D. Milly.
\newblock A mass-conservative procedure for time-stepping in models of
  unsaturated flow.
\newblock {\em Adv. Water Resour.}, 8:32--36, 1985.

\bibitem{Mualem1976}
Y.~Mualem.
\newblock A new model for predicting the hydraulic conductivity of unsaturated
  porous media.
\newblock {\em Water Resour. Res.}, 12:513--522, 1976.

\bibitem{Oulhaj2018}
A.~A.~H. Oulhaj, C.~Canc\'es, and C.~Chainais-Hillairet.
\newblock Numerical analysis of a nonlinearly stable and positive control
  volume finite element scheme for {R}ichards equation with anisotropy.
\newblock {\em ESAIM: M2AN}, 52:1533--1567, 2018.

\bibitem{Paniconi1991}
C.~Paniconi, A.~A. Aldama, and E.~F. Wood.
\newblock Numerical evaluation of iterative and noniterative methods for the
  solution of the nonlinear {R}ichards equation.
\newblock {\em Water Resour. Res.}, 27:1147--1163, 1991.

\bibitem{Paniconi1994}
C.~Paniconi and M.~Putti.
\newblock A comparison of {P}icard and {N}ewton iteration in the numerical
  solution of multidimensional variably saturated flow problems.
\newblock {\em Water Resour. Res.}, 30:3357--3374, 1994.

\bibitem{Parkin1995}
G.~W. Parkin, A.~W. Warrick, D.~E. Elrick, and R.~G. Kachanoski.
\newblock Analytical solution for one-dimensional drainage: Water stored in a
  fixed depth.
\newblock {\em Water Resour. Res.}, 31:1267--1271, 1995.

\bibitem{Philip1969}
J.~R. Philip.
\newblock Theory of infiltration.
\newblock {\em Adv. Hydrosci.}, 5:215--296, 1969.

\bibitem{ISPop2004}
I.~S. Pop, F.~Radu, and P.~Knabner.
\newblock Mixed finite elements for the {R}ichards' equation: linearization
  procedure.
\newblock {\em J. Comput. Appl. Math.}, 168:365--373, 2004.

\bibitem{Schweizer2011}
I.~S. Pop and B.~Schweizer.
\newblock Regularization schemes for degenerate {R}ichards equations and
  outflow conditions.
\newblock {\em Math. Models Methods Appl. Sci.}, 21:1685--1712, 2011.

\bibitem{Radu2004}
F.~Radu, I.~S. Pop, and P.~Knabner.
\newblock Order of convergence estimates for an {E}uler implicit, mixed finite
  element discretization of {R}ichards' equation.
\newblock {\em SIAM J. Numer. Anal.}, 42:1452--1478, 2004.

\bibitem{Radu2014}
F.~A. Radu and W.~Wang.
\newblock Convergence analysis for a mixed finite element scheme for flow in
  strictly unsaturated porous media.
\newblock {\em Nonlinear Anal. Real World Appl.}, 15:266--275, 2014.

\bibitem{Rank1983}
E.~Rank, C.~Katz, and H.~Werner.
\newblock On the importance of the discrete maximum principle in transient
  analysis using finite element methods.
\newblock {\em Int. J. Numer. Meth. Engng.}, 19:1771--1782, 1983.

\bibitem{Richards1931}
L.~A. Richards.
\newblock Capillary conduction of liquids through porous mediums.
\newblock {\em Physics}, 1:318--333, 1931.

\bibitem{Ross2003}
P.~J. Ross.
\newblock Modeling soil water and solute transport--fast, simplified numerical
  solutions.
\newblock {\em Agron. J.}, 95:1352--1361, 2003.

\bibitem{Sander1991}
G.~C. Sander, I.~F. Cunning, W.~L. Hogarth, and J.-Y. Parlange.
\newblock Exact solution for nonlinear, nonhysteretic redistribution in
  vertical soil of finite depth.
\newblock {\em Water Resour. Res.}, 27:1529--1536, 1991.

\bibitem{Sander1988}
G.~C. Sander, J.-Y. Parlange, V.~K{\"u}hnel, W.~L. Hogarth, D.~Lockington, and
  J.~P.~J. O'Kane.
\newblock Exact nonlinear solution for constant flux infiltration.
\newblock {\em J. Hydrol.}, 97:341--346, 1988.

\bibitem{Schneid2004}
E.~Schneid, P~Knabner, and F.~A. Radu.
\newblock A priori error estimates for a mixed finite element discretization of
  the {R}ichards' equation.
\newblock {\em Numer. Math.}, 98:353--370, 2004.

\bibitem{Schweizer2007}
B.~Schweizer.
\newblock Regularization of outflow problems in unsaturated porous media with
  dry regions.
\newblock {\em J. Differ. Equ.}, 237:278--306, 2007.

\bibitem{Slodicka2002}
M.~Slodicka.
\newblock A robust and efficient linearization scheme for doubly nonlinear and
  degenerate parabolic problems arising in flow in porous media.
\newblock {\em SIAM J. Sci. Comput.}, 23:1593--1614, 2002.

\bibitem{SochalaThesis2008}
P.~Sochala.
\newblock {\em M\'ethodes num\'eriques pour les \'ecoulements souterrains et
  couplage avec le ruissellement}.
\newblock PhD thesis, \'Ecole Nationale des Ponts et Chauss\'ees, 2008.

\bibitem{Suk2019}
H.~Suk and E.~Park.
\newblock Numerical solution of the {K}irchhoff-transformed {R}ichards equation
  for simulating variably saturated flow in heterogeneous layered porous media.
\newblock {\em J. Hydrol.}, 579:124213, 2019.

\bibitem{Tracy2006}
F.~T. Tracy.
\newblock Clean two- and three-dimensional analytical solutions of {R}ichards'
  equation for testing numerical solvers.
\newblock {\em Water Resour. Res.}, 42:W08503, 2006.

\bibitem{Tracy2011}
F.~T. Tracy.
\newblock Analytical and numerical solutions of {R}ichards' equation with
  discussions on relative hydraulic conductivity.
\newblock In Prof.~Lakshmanan Elango, editor, {\em Hydraulic Conductivity}.
  IntechOpen, 2011.

\bibitem{Traverso2013}
L.~Traverso, T.~N. Phillips, and Y.~Yang.
\newblock Mixed finite element methods for groundwater flow in heterogeneous
  aquifers.
\newblock {\em Comput. Fluids}, 88:60--80, 2013.

\bibitem{Uzawa1958}
H.~Uzawa.
\newblock Iterative methods for concave programming.
\newblock In {\em K. J. Arrow, L. Hurwicz, H. Uzawa (eds.)}, Studies in Linear
  and Non-Linear Programming, pages 154--165. Stanford University Press,
  Stanford, 1958.

\bibitem{vanGenuchten1980}
M.~Th. van Genuchten.
\newblock A closed-form equation for predicting the hydraulic conductivity of
  unsaturated soils.
\newblock {\em Soil Sci. Soc. Am. J.}, 44:892--898, 1980.

\bibitem{Warrick2003}
A.~Warrick.
\newblock {\em Soil {W}ater {D}ynamics}.
\newblock Oxford University Press, New York, 2003.

\bibitem{Younes2006}
A.~Younes, P.~Ackerer, and F.~Lehmann.
\newblock A new mass lumping scheme for the mixed hybrid finite element method.
\newblock {\em Int. J. Numer. Meth. Engng.}, 67:89--107, 2006.

\bibitem{Zha2013}
Y.~Zha, L.~Shi, M.~Ye, and J.~Yang.
\newblock A generalized {R}oss method for two- and three-dimensional variably
  saturated flow.
\newblock {\em Adv. Water Resour.}, 54:67--77, 2013.

\bibitem{Zha2016}
Y.~Zha, M.~Tso, L.~Shi, and J.~Yang.
\newblock Comparison of noniterative algorithms based on different forms of
  {R}ichards' equation.
\newblock {\em Environ. Model. Assess.}, 21:357--370, 2016.

\end{thebibliography}
\end{document}